\documentclass[letterpaper,leqno]{article}
\usepackage{a4, amsxtra,amsmath,amsfonts,amscd, amssymb, mathrsfs}
\usepackage[all]{xy}
\newcommand{\kkk}[1]{}

\newtheorem{thm}{Theorem}[section]
\newtheorem{prop}[thm]{Proposition} 
\newtheorem{cor}[thm]{Corollary}
\newtheorem{lem}[thm]{Lemma}
\newtheorem{ex}[thm]{Example}

\newtheorem{rem}[thm]{Remark}
\newtheorem{art}[thm]{}

\newcommand{\codim}{{\rm codim}}

\newcommand{\Div}{{\rm div}}
\newcommand{\cyc}{{\rm cyc}}
\newcommand{\Max}{{\rm Max}}

\newcommand{\Spec}{{\rm Spec}}
\newcommand{\Spf}{{\rm Spf}}

\newcommand{\ve}{\varepsilon}

\newcommand{\ind}{{\rm ind}}

\newcommand{\Star}{{\rm star}}
\newcommand{\supp}{{\rm supp}}

\newcommand{\val}{{\rm val}}
\newcommand{\Val}{{\rm Val}}
\newcommand{\vol}{{\rm  vol}}

\newcommand{\Acal}{{\mathscr A}}
\newcommand{\Bcal}{{\mathscr B}}
\newcommand{\Ccal}{{\mathscr C}}
\newcommand{\Dcal}{{\mathscr D}}

\newcommand{\Lcal}{{\mathscr L}}
\newcommand{\Mcal}{{\mathscr M}}
\newcommand{\Ocal}{{\mathscr O}}

\newcommand{\Scal}{{\mathscr S}}
\newcommand{\Tcal}{{\mathscr T}}
\newcommand{\Ucal}{{\mathscr U}}
\newcommand{\Xcal}{{\mathscr X}}
\newcommand{\Ycal}{{\mathscr Y}}
\newcommand{\Xfrak}{{\frak X}}
\newcommand{\Xan}{{X_\kdop^{\rm an}}}
\newcommand{\Ufrak}{{\frak U}}

\newcommand{\qdop}{{\mathbb Q}}
\newcommand{\ndop}{{\mathbb N}}
\newcommand{\rdop}{{\mathbb R}}

\newcommand{\kdop}{{\mathbb K}}
\newcommand{\bdop}{{\mathbb B}}
\newcommand{\adop}{{\mathbb A}}
\newcommand{\ldop}{{\mathbb L}}
\newcommand{\zdop}{{\mathbb Z}}

\newcommand{\metr}{\|\hspace{1ex}\|}

\newcommand{\ub}{{\mathbf u}}
\newcommand{\ubb}{{\overline{\mathbf u}}}
\newcommand{\lb}{{\mathbf l}}
\newcommand{\mb}{{\mathbf m}}

\newcommand{\tb}{{\mathbf t}}
\newcommand{\xb}{{\mathbf x}}
\newcommand{\yb}{{\mathbf y}}
\newcommand{\zb}{{\mathbf z}}
\newcommand{\Tor}{{\mathbb G}_m^n}
\newcommand{\rtor}{{\rdop^n/\Lambda}}
\newcommand{\g}{{\frak g}_X}
\newcommand{\gp}{{\frak g}_X^+}
\newcommand{\gh}{\hat{\frak g}_X}
\newcommand{\ghp}{\hat{\frak g}_X^+}

\newcommand{\proof}{\noindent {\bf Proof: \/}}

\newcommand{\qed}{{ \hfill $\square$}}

\newcommand{\Deltabar}{{\overline{\Delta}}}
\newcommand{\Ccalbar}{{\overline{\Ccal}}}
\newcommand{\Omegabar}{{\overline{\Omega}}}
\newcommand{\Sigmabar}{{\overline{\Sigma}}}
\newcommand{\sigmabar}{{\overline{\sigma}}}
\newcommand{\valbar}{{\overline{\val}}}

\newcommand{\Int}{{\rm int}}
\newcommand{\relint}{{\rm relint}}
\newcommand{\card}{{\rm card}}
\title{Tropical varieties \linebreak for non-archimedean analytic spaces}
\date{\today}
\author{Walter Gubler}
\begin{document}
\maketitle


\section{Introduction}

For the whole paper, $\kdop$ denotes an algebraically closed field endowed with a non-trivial non-archimedean complete absolute value $|\phantom{a}|$. The corresponding valuation is $v:=-\log |\phantom{a}|$ with value group $\Gamma := v(\kdop^\times)$. The valuation ring is denoted by $\kdop^\circ$. Note that the residue field $\tilde{\kdop}$ is algebraically closed. In Theorem \ref{Theorem 3}, \S 8 and in the second part of \S 9, we start with a field $K$ endowed with a discrete valuation and we choose $\kdop$ to be the completion of the algebraic closure of the completion of $K$ (see \cite{BGR}, \S 3.4, for these properties of $\kdop$).

On the torus $\Tor$, we always fix coordinates $x_1, \dots , x_n$ and we consider the map
$$\val:\Tor(\kdop) \longrightarrow \rdop^n, \quad \xb \mapsto 
\left( -\log |x_1|, \dots, -\log |x_n| \right).$$ 
For an irreducible closed algebraic subvariety $X$ of $\Tor$ over $\kdop$, the closure of $\val(X)$ in $\rdop^n$ is called a tropical variety. This is the main object of study in tropical algebraic geometry. We refer to \cite{Mi} for a survey of this relatively new area of research. Einsiedler, Kapranov and Lind \cite{EKL} have shown that the tropical variety of $X$ is a connected totally concave $\Gamma$-rational polyhedral set in $\rdop^n$ of pure dimension $\dim(X)$. Here and in the following, the reader is referred to the appendix for the terminology used from convex geometry. In the present paper, the following analytic generalization is given:

\begin{thm} \label{Theorem 1} 
Let $X$ be an irreducible closed analytic subvariety of $\Tor$ over $\kdop$ of dimension $d$. Then the tropical variety associated to $X$ is a connected totally concave locally finite union of $d$-dimensional $\Gamma$-rational polytopes.
\end{thm}

Note that the map $\val$ is continuous with respect to the Berkovich analytic structure on $\Tor$ and therefore the tropical variety is obviously connected and compact. This makes it clear that we benefit a lot by using Berkovich analytic spaces and methods from formal geometry (see \S 2 for a summary). In \S 4, we generalize a result of Mumford to study the special fibre of the analytic subdomain $U_\Delta:= \val^{-1}(\Delta)$ of $\Tor$ associated to a $\Gamma$-rational polytope in $\rdop^n$. This allows us to apply the theory of toric varieties to the reduction of $U_\Delta$. In \S 5, we prove Theorem 1 from the corresponding local case in $U_\Delta$.

The applications will deal with a a totally degenerate abelian variety $A$ over $\kdop$, i.e. $A^{\rm an} = (\Tor)_\kdop^{\rm an}/M$ for a discrete subgroup $M$ of $\Tor(\kdop)$ which is mapped isomorphically onto the complete lattice $\Lambda := \val(M)$ in $\rdop^n$. We get a canonical map $\valbar: A^{\rm an} \rightarrow \rdop^n/\Lambda$ and hence a tropical variety $\valbar(X^{\rm an})$ associated to a closed analytic subvariety $X$ of $A$. In \S 6, we will show that Theorem \ref{Theorem 1} holds also in this framework. This is quite obvious by lifting $X$ to $\Tor$ leading to a periodic tropical variety in $\rdop^n$. As a consequence, we obtain the following dimensionality theorem:

\begin{thm} \label{Theorem 2}
Let $X'$ be a smooth algebraic variety over $\kdop$ with a strictly semistable formal  $\kdop^\circ$-model $\Xcal'$ (see \ref{semistable} and \ref{models}) and let $f: X' \rightarrow A$ be a morphism over $\kdop$. Then the special fibre of $\Xcal'$ has a $\tilde{\kdop}$-rational point contained in at least $1+\dim f(X')$ irreducible components.
\end{thm}
If $X'$ is projective, then we may use a strictly semistable projective $\kdop^\circ$-model. If $X'$ has good reduction at $v$, then $f$ is constant. We will postpone the proof of Theorem \ref{Theorem 2} to the first part of \S 9 where it can be given very neatly and where also a generalization to arbitrary abelian varieties is given. 

In non-archimedean analysis, no good definition is known for the first Chern form of a metrized line bundle. However, Chambert-Loir \cite{Ch} has introduced measures $c_1(\overline{L}_1) \wedge \dots \wedge c_1(\overline{L}_d)$ on the Berkovich analytic space associated to a $d$-dimensional projective variety analogous to the corresponding top dimensional forms in differential geometry. In \S 3, we give a slightly more general approach to these measures using local heights. Our main application of the theory of tropical analytic varieties is the following result:

\begin{thm} \label{Theorem 3}
Let $X$ be a closed subvariety of the abelian variety $A$ over $K$. We assume that $A_\kdop^{\rm an}$ is totally degenerate over $\kdop$ and that $X$ is of pure dimension $d$. Let $\overline{L}_1, \dots , \overline{L}_d$ be ample  line bundles on $A$ endowed with canonical metrics. Then $\mu := (\valbar)_*(c_1(\overline{L}_1|_X) \wedge \cdots \wedge c_1(\overline{L}_d|_X)$ is a strictly positive piecewise Haar measure on the polytopal set $\valbar(X_\kdop^{\rm an})$.
\end{thm}

By a {\it piecewise Haar measure} on the polytopal set $\valbar(X_\kdop^{\rm an})$, we mean that $\valbar(X_\kdop^{\rm an})$ is a union of $d$-dimensional polytopes $\overline{\sigma}$ such that $\mu$ is a real multiple of the relative Lebesgue measure on $\overline{\sigma}$, and {\it stricly positive} means that the multiple is $>0$. 

The proof of Theorem \ref{Theorem 3} will be given in \S 8. It is based on our studies in \S 6 and \S 7 of Mumford's model $\Acal$ of $A$ associated to a $\Gamma$-rational polytopal decomposition of $\rtor$. The main idea is to work with a generic polytopal decomposition $\Ccalbar$ defined over a sufficiently large base extension $\kdop'$ of $\kdop$. Then the irreducible components of the special fibre of the closure of $X$ in $\Acal$ are toric varieties by Theorem \ref{periodic transverse}. Moreover, this holds for the sequence $\frac{1}{m}\Ccal$, $m \in \ndop$, of polytopal decompositions  and this allows us to  compute $\mu$ by Tate's limit argument leading to an explicit expression for $\mu$ in Theorem \ref{explicit expression} and proving Theorem \ref{Theorem 3}. 
In Theorem \ref{explicit generalization}, we prove that $(c_1(\overline{L}_1|_X) \wedge \cdots \wedge c_1(\overline{L}_d|_X)$ itself is induced by an explicit strictly positive Haar measure on the skeleton of a strictly semistable alteration of $X$.

In \cite{Gu4},  Theorem \ref{Theorem 3} is essential to prove the following case of {\it Bogomolov's conjecture} over the function field $F:=k(B)$. Here, $B$ is an integral projective variety over the algebraically closed field $k$ such that $B$ is regular in codimension $1$.  The prime divisors on $B$ are weighted by the degree with respect to a fixed ample class  leading to a theory of heights.

{\bf Bogomolov conjecture (\cite{Gu4}, Theorem 1.1)} {\it
Let $A$ be an abelian variety over the function field $F$ which is totally degenerate at some place $v$ of $F$. Let $X$ be a closed subvariety of $A$ defined over the algebraic closure $\overline{F}$  which is not a translate of an abelian subvariety by a torsion point. For every ample symmetric line bundle $L$ on $A$, there is $\ve > 0$ such that
$$X(\ve):=\{P \in X(\overline{F}) \mid \hat{h}_L(P) \leq \ve\}$$
is not Zariski dense in $X$, where $\hat{h}_L$ denotes the N\'eron--Tate height with respect to $L$.}

\vspace{3mm}

\centerline{\it Terminology}

In $A \subset B$, $A$ may be  equal to $B$. The complement of $A$ in $B$ is denoted by $B \setminus A$ \label{setminus} as we reserve $-$ for algebraic purposes. The zero is included in $\ndop$.

All occuring rings and algebras are commutative with $1$. If $A$ is such a ring, then the group of multiplicative units is denoted by $A^\times$. A variety over a field is a separated reduced scheme of finite type. However,  a formal analytic variety is not necessarily reduced (see \S 2). For the degree of a map $f:X \rightarrow Y$ of irreducible varieties, we use either $\deg(f)$ or $[X:Y]$. The multiplicity of an irreducible component $Y$ of a scheme $S$ is denoted by $m(Y,S)$.

For $\mb \in \zdop^n$, let $\xb^\mb:=x_1^{m_1} \cdots x_n^{m_n}$. The standard scalar product of $\ub,\ub' \in \rdop^n$ is denoted by  $\ub \cdot \ub':=u_1 u_1' + \dots + u_n u_n'$. For the notation used from convex geometry, we refer to the appendix (see also \ref{toric convex geometry} for the periodic case).

\small
The author thanks J. Eckhoff, K. K\"unnemann and F. Oort for precious discussions, and the referee for his suggestions. Part of the research in this paper was done during a $5$ weeks stay at the CRM in Barcelona.
\normalsize

\section{Analytic and formal geometry}

In this section, we gather the results needed from Berkovich spaces and formal geometry.

\begin{art} \rm \label{affinoid algebras}
The completion of $\kdop[x_1, \dots, x_n]$ with respect to the Gauss norm is called the {\it Tate algebra} and is denoted by $\kdop \langle x_1, \dots, x_n \rangle$. It consists of the strictly convergent power series on the closed unit ball $\bdop^n$ in $\kdop^n$.  

A $\kdop$-affinoid algebra $\Acal$ is isomorphic to a quotient $\kdop \langle x_1, \dots, x_n \rangle/I$ and the maximal spectrum $\Max(\Acal)$ is equal to the zero set $Z(I) \subset \bdop^n$ of the ideal $I$. The supremum semi-norm of $\Acal$ on $Z(I)$ is denoted by $|\phantom{a}|_{\rm sup}$. Setting
$$\Acal^\circ:= \{a \in \Acal \mid |a|_{\rm sup} \leq 1 \}, 
\quad \Acal^{\circ \circ}:=\{a \in \Acal \mid |a|_{\rm sup} < 1 \},$$
the {\it residue algebra} is defined by $\tilde{\Acal}:= \Acal^\circ / \Acal^{\circ \circ}$. It is a finitely generated reduced $\tilde{\kdop}$-algebra. For details about affinoid algebras, we refer to \cite{BGR}.
\end{art}

\begin{art} \rm \label{Berkovich spectrum} 
The {\it Berkovich spectrum} $\Mcal(\Acal)$ of a $\kdop$-affinoid algebra $\Acal$ is defined as the set of semi-norms $p$ on $\Acal$ satisfying
$$p(ab)=p(a)p(b), \quad p(1)=1 \quad \text{and} \quad p(a) \leq |a|_{\rm sup}$$
for all $a,b \in \Acal$. It is endowed with the coarsest topology such that the maps $p \mapsto p(a)$ are continuous for all $a \in \Acal$.

The Berkovich spectrum is compact and every $x \in \Max(\Acal)$ gives rise to a semi-norm $a \mapsto |a(x)|$ such that we may view $\Mcal(\Acal)$ as a compactification of $\Max(\Acal)$.
We refer to \cite{Ber} for proofs and more details.

The affine $\tilde{\kdop}$-variety $\Spec(\tilde{\Acal})$ is called the {\it reduction} of $\Mcal(\Acal)$ and the reduction map $p \mapsto \tilde{p}:= \{p<1\}/\Acal^{\circ \circ}$ is surjective. If $\wp$ is a minimal prime ideal of $\tilde{\Acal}$, then there is a unique $p \in \Mcal(\Acal)$ with $\tilde{p}=\wp$ (see \cite{Ber}, Proposition 2.4.4).

\end{art}

\begin{art} \rm \label{subdomains}
An {\it affinoid subdomain} of $X:=\Mcal(\Acal) = \Mcal(\kdop \langle \mathbf x \rangle /I)$ is characterized by an universal property (see \cite{Ber}, 2.2 or \cite{BGR}, 7.2.2). By a theorem of Gerritzen and Grauert (\cite{BGR}, Corollary 7.3.5/3), an affinoid subdomain is a finite union of {\it rational domains}. The latter are defined by
$$X \left(\frac{\mathbf f}{g}\right) := \{x \in X \mid |f_j(x)| \leq |g(x)|, \, j=1,\dots,r\}$$
where $g,f_1, \dots, f_r \in \Acal$ are without common zero. The corresponding affinoid algebra is
$$\Acal \left\langle \frac{\mathbf f}{g}\right\rangle := \kdop \langle \mathbf x, y_1, \dots, y_r \rangle / \langle I, g(\mathbf x) y_j - f_j \mid j=1, \dots ,r \rangle$$
(see \cite{BGR}, Proposition 7.2.3/4). If $g=1$, then $X(\mathbf f)$ is called a {\it Weierstrass domain} in $\Mcal(\Acal)$.
\end{art}

\begin{art} \rm \label{analytic spaces}
An {\it analytic space} $X$ over $\kdop$ is given by an atlas of affinoid subdomains $U=\Mcal(\Acal)$. For the precise definition, we refer to \cite{Ber2}, \S 1. (Note that our definition corresponds to strictly analytic spaces in the notation of \cite{Ber2}, p. 22). The technical difficulty in this definition arises from the fact that the charts $U$ are not open in $X$ but compact. The sheaf of structure $\Ocal_X$ is only defined on the Grothendieck topology of $X$ and it is characterized by $\Ocal_X(U)= \Acal$.
\end{art}

\begin{art} \rm \label{formal analytic varieties}
Let $\Acal$ be a $\kdop$-affinoid algebra. 
A subset $U$ of $\Mcal(\Acal)$ is called {\it formal open} if there is a  open subset $V$ of the reduction $\Spec(\tilde{\Acal})$ such that $U=\pi^{-1}(V)$. The resulting quasi-compact topology on $\Mcal(\Acal)$ is called the {\it formal topology}. Together with the restriction of $\Ocal_{\Mcal(\Acal)}$ to the formal topology, we get a ringed space called a {\it formal affinoid variety} over $\kdop$ and denoted by $\Spf(\Acal)$. By definition, a morphism of affinoid varieties over $\kdop$ is induced by a reverse homomorphism of the corresponding  $\kdop$-affinoid algebras. For details, we refer to \cite{Bo}. 

A {\it formal analytic variety} over $\kdop$ is a $\kdop$-ringed space $\Xfrak$ which has a locally finite open atlas of formal affinoid varieties over $\kdop$.  It has a {\it reduction} $\tilde{\Xfrak}$ and a {\it generic fibre} $\Xfrak^{\rm an}$. If $\Xfrak = \Spf(\Acal)$, then $\tilde{\Xfrak}= \Spec(\tilde{\Acal})$ and $\Xfrak^{\rm an} = \Mcal(\Acal)$. In general, the $\tilde{\kdop}$-variety $\tilde{\Xfrak}$ and the analytic space $\Xfrak^{\rm an}$ are obtained by gluing processes (see \cite{Bo} and \cite{Ber3}, \S 1).

By \ref{Berkovich spectrum}, there is a surjective reduction map $\Xfrak^{\rm an} \rightarrow \tilde{\Xfrak}$, given locally by $p \mapsto \tilde{p}$. For every irreducible component $Y$, there is a unique $\xi_Y \in \Xfrak^{\rm an}$ which reduces to the generic point of $Y$. 
\end{art}

\begin{art} \rm \label{admissible formal schemes}
A $\kdop^\circ$-algebra is called {\it admissible} if it is isomorphic to $\kdop^\circ\langle x_1, \dots , x_n \rangle / I$ for an ideal $I$ and if $A$ has no $\kdop^\circ$-torsion. An {\it admissible formal scheme} $\Xcal$ over $\kdop^\circ$ is a formal scheme which has a locally finite atlas of open subsets isomorphic to $\Spf(A)$ for  admissible $\kdop^\circ$-algebras $A$. The lack of $\kdop^\circ$-torsion is equivalent to flatness over $\kdop^\circ$. 

These spaces are studied in detail by Bosch and L\"utkebohmert (\cite{BL3}, \cite{BL4}) based on results of Raynaud. Note that the locally finiteness condition for the atlas in the definitions of formal analytic varieties and admissible formal schemes does not occur in \cite{Bo} and the above quotes. We need it only to define the generic fibre as an analytic space and it could be omitted working with rigid analytic spaces (see \cite{Gu2}). 

The {\it special fibre} $\tilde{\Xcal}$ of an admissible formal scheme $\Xcal$ over $\kdop^\circ$ is a scheme of locally finite type over $\tilde{\kdop}$ with the same underlying topological space as $\Xcal$ and with $\Ocal_{\tilde{\Xcal}} := \Ocal_{\Xcal} \otimes_{\kdop^\circ} \tilde{\kdop}$.

There is a formal analytic variety $\Xcal^{\rm f-an}$ associated to $\Xcal$. If $\Xcal = \Spf(\Acal)$, then $\Acal:= A \otimes_{\kdop^\circ} \kdop$ is a $\kdop$-affinoid algebra and we set $\Xcal^{\rm f-an}:= \Spf(\Acal)$. In general, $\Xcal^{\rm f-an}$  is obtained by a gluing process. The canonical morphism $(\Xcal^{\rm f-an})\sptilde  \rightarrow \tilde{\Xcal}$ is finite and surjective (see \cite{BL1}, \S 1). 

The  analytic space $\Xcal^{\rm an}:= (\Xcal^{\rm f-an})^{\rm an}$ is called the {\it generic fibre} of $\Xcal$. Similarly as in \ref{formal analytic varieties}, there is a surjective reduction map $\Xcal^{\rm an} \rightarrow \tilde{\Xcal}$. 

If $\Xfrak$ is a formal analytic variety over $\kdop$, we may reverse the above process replacing locally $\Spf(\Acal)$ by $\Spf(\Acal^\circ)$ to get a formal scheme $\Xfrak^{\rm f-sch}$ over $\kdop^\circ$.

All the above constructions are functorial. The functors $\Xcal \rightarrow \Xcal^{\rm f-an}$ and $\Xfrak \rightarrow \Xfrak^{\rm f-sch}$ give an equivalence between the category of admissible formal schemes over $\kdop^\circ$ with reduced fibre and the category of reduced formal analytic varieties over $\kdop$. The reductions are the same, which allows us to flip from one category to the other. For details, see \cite{BL1}, \S 1, and \cite{Gu1}, \S1. 
\end{art}

\begin{art} \rm \label{schemes}
Let $X$ be a scheme  of finite type over a subfield $K$ of $\kdop$. The  analytic space $X^{\rm an}$ over $\kdop$ associated to $X$ is constructed in the following way:  
By using a gluing process, we may assume that $X$ is a closed subscheme of $\adop^n_K$. For $r \in |\kdop^\times|$, the intersection of $X^{\rm an}$ with the closed ball of radius $r$ and center $\mathbf 0$ is defined by the same set of equations as $X$ in $\adop^n_K$. If we glue the balls for $r \to \infty$, then we get $X^{\rm an}$. 
For more details about this functorial construction and the following GAGA theorems, we refer to \cite{Ber}, 3.4.

$X$ is reduced, normal, regular, smooth, $d$-dimensional or connected if and only if $X^{\rm an}$ has the same property. $X$ is separated, resp.\! proper over $K$ if and only if $X^{\rm an}$ is Hausdorff, resp.\! com\-pact. A morphism of schemes of finite type over $K$ is flat, unramified, \'etale, smooth, an open immersion, a closed immersion, dominant, proper, finite if and only if this holds for $\varphi^{\rm an}$. 

Let $\Xcal$ be a flat scheme of finite type over $K^\circ$ with generic fibre $X$ and let $\pi \in \kdop^{\circ \circ}$. Then the associated formal scheme $\hat{\Xcal}$ over $\kdop^\circ$, defined locally by replacing the coordinate ring $A$ by the $\pi$-adic completion of $A \otimes_{K^\circ} \kdop^\circ$, is admissible. Moreover, the special fibre of $\hat{\Xcal}$ is isomorphic to the base change of the special fibre of $\Xcal$ to $\tilde{\kdop}$. Note that $\hat{\Xcal}^{\rm an}$ is an analytic subdomain of $X^{\rm an}$ with $\hat{\Xcal}^{an}(\kdop)$ consisting of the $\kdop^\circ$-integral points of $\Xcal$. If $\Xcal$ is proper over $K^\circ$, then $\hat{\Xcal}^{\rm an}=X^{\rm an}$. For details, we refer to \cite{Gu2}, \S 6). If $\kdop^\circ$ is not a discrete valuation ring, one has to use \cite{Ul}.
\end{art}

\begin{art} \label{etale} \rm 
In the following, we consider an {\it \'etale} morphism $\varphi: \Ycal \rightarrow \Xcal$ of admissible formal schemes over $\kdop^\circ$, i.e. the reduction
$$\varphi_\lambda: (\Ycal, \Ocal_{\Ycal}/\lambda \Ocal_{\Ycal}) \longrightarrow 
(\Xcal, \Ocal_{\Xcal}/\lambda \Ocal_{\Xcal}) $$
is an \'etale morphism of schemes for all $\lambda \in \kdop^{\circ \circ}$. Let $X$, $Y$ be the generic fibres of $\Xcal$ and $\Ycal$. 

For $\tilde{P} \in \tilde{\Xcal}(\tilde{\kdop})$, the {\it formal fibre} $X_+(\tilde{P}) := \{ x \in X \mid \tilde{x} = \tilde{P} \}$ is an open analytic subspace of $X$. Indeed, let $\Spf(A)$ be a formal affine neighbourhood of $\tilde{P}$ in $\Xcal$ and let $f_1, \dots, f_r \in A$ such that $\tilde{P}$ is the only common zero of $\tilde{f_1}, \dots, \tilde{f_r} \in \tilde{A}$  in $\Spec(\tilde{A})$, then 
$$X_+(\tilde{P}) = \{x \in \Mcal(A \otimes_{\kdop^\circ}\kdop) \mid |f_1(x)|<1, \dots , |f_r(x)|<1 \}$$
is an open subdomain of $X$.\end{art}

The following result is a special case of \cite{Ber4}, Lemma 4.4. We give an elementary proof here based on the implicit function theorem.

\begin{prop} \label{implicit function}
Let $\varphi$ be as above and let $\tilde{Q} \in \tilde{\Ycal}(\tilde{\kdop})$ with $\tilde{P}=\tilde{\varphi}(\tilde{Q})$. Then $\varphi$ restricts to an isomorphism $Y_+(\tilde{Q}) \stackrel{\sim}{\rightarrow} X_+(\tilde{P})$ of formal fibres.
\end{prop}

\proof By the local description of \'etale morphisms of schemes, we may assume that $\Xcal=\Spf(A)$ and $\Ycal=\Spf(B)$, where 
$$B=\left(A[t]/ \langle p(t) \rangle \right)_{\{q(t)\}}.$$
Here, $p(t), q(t) \in A[t]$ and the monic polynomial $p(t)$ has the property that the residue class of $\frac{d}{dt} p$ is invertible in $B$ (see \cite{Ber3}, \S 2).
Note that the admissible $\kdop^\circ$-algebra $A$ has the form $\kdop^\circ\langle x_1, \dots, x_n \rangle/ I$. There is $Q \in Y_+(\tilde{Q})(\kdop)$ with reduction $\tilde{Q}$ and hence $P:=\varphi(Q)$ has reduction $\tilde{P}$ (see \cite{BGR}, Theorem 7.1.5/4). 

There is $p(\mathbf x, t) \in  \kdop^\circ\langle\mathbf x\rangle[t]$ with residue class $p(t) \in A[t]$. Clearly, $p(\mathbf x , t)$ has Gauss norm $1$ and $|\frac{\partial}{\partial t} p(Q)|=1$. By the local Eisenstein theorem (\cite{BG}, Theorem 11.5.14), there is a unique formal power series $\xi(\mathbf x) \in \kdop[[\mathbf x]]$ with
$p(\mathbf x, \xi(\mathbf x))=0 $ and $\xi(Q)=0$. 
Moreover, $\xi$ is convergent on $X_+(\tilde{P})$. By Hensel's lemma, we easily deduce that $\mathbf x \mapsto (\mathbf x, \xi(\mathbf x))$ is the inverse of  $\varphi: Y_+(\tilde{Q}) \rightarrow X_+(\tilde{P})$. \qed 

\begin{art} \label{semistable} \rm 
An admissible formal scheme $\Xcal$ over $\kdop^\circ$ is called {\it strictly semistable} if $\Xcal$ is covered by formal open subsets $\Ucal$ with an \'etale morphism 
$$\psi: {\Ucal} \longrightarrow {\Scal} :=\Spf \left( \kdop^\circ \langle x_0, \dots, x_n \rangle / \langle x_0 \dots x_r = \pi \rangle \right)$$
for some $r \leq n$ and $\pi \in \kdop^{\circ \circ}$ (depending on $\Ucal$). 
\end{art}

\begin{prop} \label{semistable properties} 
Let $\Xcal$ be a strictly semistable formal scheme over $\kdop^\circ$.
\begin{itemize}
\item[(a)] The special fibre $\tilde{\Xcal}$ is reduced and hence $\Xcal$ is the formal scheme associated to a formal analytic variety.
\item[(b)] For every $\tilde{P} \in \tilde{\Xcal}(\tilde{\kdop})$, there is a formal open neighbourhood $\Ucal = \Spf(B)$ in $\Xcal$ and an \'etale morphism $\psi$ as in \ref{semistable} such that $\tilde{\psi}(\tilde{P})= \tilde{\mathbf 0}$ and such that every irreducible component of $\tilde{\Ucal}$ passes through $\tilde{P}$.
\item[(c)] Let  $\Ucal$ as in (b) and let $\gamma_j:= \psi^*(x_j)$. The irreducible components of $\tilde{\Ucal}$ are equal to $Y_0, \dots , Y_r$, where $Y_j$ is the zero-scheme of $\tilde{\gamma}_j$. Moreover, the group of units of the $\kdop$-affinoid algebra $B \otimes_{\kdop^\circ} \kdop$ is the direct product of $\kdop^\times  B^\times$ and the free abelian subgroup with basis $\gamma_1, \dots, \gamma_r$.
\end{itemize}
\end{prop}

\proof Property (b) is immediate from strict semistability. For (a) and (c), we may assume $\Ucal = \Xcal$. The special fibre $\tilde{{\Scal}}$ is the zero-scheme of $\tilde{x}_0 \cdots \tilde{x}_r$ in $\adop^{n+1}_{\tilde{\kdop}}$, hence it is reduced and has $r+1$ irreducible components $\tilde{x}_j = \tilde{0}$. Since $\tilde{\psi}$ is \'etale, $\tilde{\Xcal}$ is reduced (\cite{EGA IV}, Proposition 17.5.7) and (a) follows from \ref{admissible formal schemes}. Moreover, $\tilde{\psi}^{-1}\{\tilde{x}_j = \tilde{0}\}$ is smooth and every irreducible component passes through $\tilde{P}$, hence $Y_j$ is irreducible. By flatness of $\tilde{\psi}$, every irreducible component of $\tilde{\Xcal}$ has this form and we get the first part of (c).

We denote the $\kdop$-affinoid algebra $B \otimes_{\kdop^\circ} \kdop$ by ${\cal B}$. Let $u \in {\cal B}^\times$. By Proposition \ref{implicit function}, the formal fibre $X_+(\tilde{P})$ is isomorphic to the open subdomain of $\bdop^n$, given by $|x_j| < 1$ for $j=1, \dots, n$ and $|x_1 \cdots x_r| > |\pi|$. Hence $u$ has the following power series development on $X_+(\tilde{P})$:
$$u|_{X_+(\tilde{P})} =  \sum_{m_1= -\infty}^\infty \dots \sum_{m_r=-\infty}^\infty \sum_{m_{r+1}=0}^\infty  \dots  \sum_{m_n =0}^\infty a_\mb \gamma_1^{m_1} \dots \gamma_n^{m_n}.$$
Since $u$ is a unit, there is a dominant term $a_\mb  \gamma_1^{m_1} \dots \gamma_n^{m_r}$ with $m_{r+1}= \dots = m_n =0$, i.e.
$$|a_\mb \gamma_1^{m_1} \dots \gamma_n^{m_n}| >
|a_{\mathbf s} \gamma_1^{s_1} \dots \gamma_n^{s_n}|$$
on $X_+(\tilde{P})$ for all $\mathbf s \neq \mb$ (use \cite{BGR}, Lemma 9.7.1/1). For $j= 0, \dots, r$, we will compute the multiplicity $m(u,Y_j)$ of $u$ in $Y_j$ (see \cite{Gu2} for the theory of divisors on admissible formal schemes). By compactness of $X$, we easily construct a sequence $P_k \in X_+(\tilde{P})$, convergent to $P_\infty \in X$ and with $|\gamma_i(P_k)| \to 1$ for $i \in \{0, \dots,n\} \setminus \{j\}$. By \cite{Gu3}, Proposition 7.6, $m(u,Y_j)=v(u(P_\infty))$ and hence
$$m(u,Y_j)= \lim_{k \to \infty} v(u(P_k)) = v(a_\mb) + m_j v(\pi).$$
We conclude that the Weil divisor of $u':=u a_\mb^{-1}  \gamma_0^{-m_0} \cdots \gamma_r^{-m_r}$ is zero on $\Xcal$ and hence $|u'(\xi_{Y_j})|=1$. As this holds for all $j=1, \dots, r$ and since $\xi_{Y_0},, \dots, \xi_{Y_r}$ is the Shilov boundary of $X$ (see \cite{Ber}, Proposition 2.4.4), we conclude that $u'$ is a unit in $\cal B^\circ$. By (a) and \ref{admissible formal schemes}, we have $B= \cal B^\circ$ and hence
$${\cal B}^ \times = \kdop^\times B^\times \gamma_1^\zdop \cdots \gamma_r^\zdop.$$
Restriction to the formal fibre $X_+(\tilde{P})$ shows that $\gamma_1, \dots, \gamma_r$ are multiplicatively independent. Moreover, the restriction of an element from $\kdop^\times  B^\times$ to $X_+(\tilde{P})$ has constant dominant term and hence the product of $\kdop^\times  B^\times$ with $\gamma_1^\zdop \cdots \gamma_r^\zdop$ is direct. \qed

\section{Local heights of subvarieties}

In this section, we summarize the theory of local heights of subvarieties. We use the formal and analytic geometry from the previous section. This allows larger flexibility in choosing models and metrics from which we benefit in Section 8.
At the end, we generalize Chambert-Loir's measures associated to metrized line bundles using a different approach through local heights. 
Apart from the definitions in \ref{models} and \ref{formal metrics}, this section will be used only in \S 8 and \S 9.
We consider a proper scheme $X$ over $\kdop$. Note that $X^{\rm an }$ is compact.

\begin{art} \label{models} \rm
A {\it formal $\kdop^\circ$-model} of $X$ is an admissible formal scheme with generic fibre $X^{\rm an}$.

On  analytic spaces, formal analytic varieties and admissible formal schemes, we may define line bundles, sections and Cartier divisors in the usual way. A {\it horizontal cycle} on a formal $\kdop^\circ$-model $\Xcal$ of $X$ is just a cycle on $X$. A {\it vertical cycle} on $\Xcal$ is a cycle on the special fibre $\tilde \Xcal$ with coefficients in $\Gamma$. A {\it cycle on $\Xcal$} is the formal sum of a horizontal and a vertical cycle.
\end{art}

\begin{art} \label{formal metrics} \rm
Let $L$ be a line bundle on $X$. It induces a line bundle $L^{\rm an}$ on $X^{\rm an}$. A {\it formal $\kdop^\circ$-model} of $L$ is a line bundle $\Lcal$ on a formal $\kdop^\circ$-model $\Xcal$ of $X$ with generic fibre $\Lcal^{\rm an}$ equal to $L^{\rm an}$.

A metric $\metr$ on $L^{\rm an}$ is said to be a {\it formal metric} if there is a formal $\kdop^\circ$-model $\Lcal$ of $L$ such that
for every formal trivialization $\Ucal$ of $\Lcal$ and every $s \in \Gamma(\Ucal,\Lcal)$ corresponding to $\gamma \in \Ocal_\Xcal(\Ucal)$, we have 
$\|s(x)\|=|\gamma(x)|$ on $\Ucal^{\rm an}$. 
The formal metric is called {\it semipositive} if the reduction $\tilde \Lcal$ of $\Lcal$ on $\tilde \Xcal$ is numerically effective (see \cite{Kl}). Every line bundle on $X$ has a formal metric (\cite{Gu2}, Corollary 7.7). 

A metric on $L^{\rm an}$ is called a {\it root of a formal metric} if some positive tensor power is a formal metric.
On the space of continuous metrics on $L^{\rm an}$, we use the distance function
$$d(\metr, \metr'):= \sup_{x \in X^{\rm an}} |\log(\metr/\metr')(x)|$$
where $(\metr/\metr')(x)$ is evaluated at the section $1$ of $O_X^{\rm an}=L^{\rm an} \otimes (L^{\rm an})^{-1}$.
\end{art}

\begin{prop} \label{density}
The roots  of formal metrics are dense in the space of continuous metrics on $L^{\rm an}$. In particular, the set of roots of formal metrics on $O_{X}^{\rm an}$ is embedded onto a dense subset of $C(X^{\rm an})$ by the map $\metr \mapsto -\log \metr$.
\end{prop}

\proof This is Theorem 7.12 in \cite{Gu2} holding more generally for compact analytic  spaces. \qed 

\begin{art} \label{pseudo-divisors} \rm
A {\it metrized pseudo-divisor} $\hat{D}$ on $X$ is a quadruple $\hat{D}:=(L, \metr, Y, s)$ where $L$ is a line bundle on $X$, $\metr$ is a metric on $L^{\rm an}$, $Y$ is a closed subset of $X$ and $s$ is a nowhere vanishing section of $L$ on $X \setminus Y$. Then $D:=(L,Y,s)$ is a pseudo-divisor on $X$ (as in \cite{Fu1}, 2.2). The {\it support} $Y$  is denoted by $\supp(D)$ and $O(D):=L$. The most relevant example for applications is the case of an invertible meromorphic section $s$ of a metrized line bundle $(L, \metr)$, where the associated pseudo-divisor $\widehat{\Div}(s)$ is defined by choosing $Y$ as the support of the Cartier divisor $\Div(s)$. Since pseudo-divisors are closed under pull-back, it is much easier to formulate the intersection theory for pseudo-divisors instead of Cartier-divisors.

On a formal $\kdop^\circ$-model $\Xcal$ of $X$, there is a refined intersection theory of formally metrized pseudo-divisors with cycles. It has the properties expected from algebraic geometry (see \cite{Gu2} and \cite{Gu3}). 

For a $t$-dimensional cycle $Z$ on $X$ and formally metrized pseudo-divisors $\hat{D}_0, \dots, \hat{D}_t$ with
\begin{equation} \label{suppcond}
\supp(D_0) \cap \dots \cap \supp(D_t) \cap \supp(Z) = \emptyset,
\end{equation}
there is a {\it local height}
$\lambda(Z):= \lambda_{\hat{D}_0, \dots, \hat{D}_t}(Z)$ 
defined as the intersection number of $\hat{D}_0, \dots, \hat{D}_t$ and $Z$ on a joint formal $\kdop^\circ$-model. In case of a discrete valuation and algebraic $\kdop^\circ$-models, this is the usual intersection product and hence  we get the local heights used in Arakelov geometry. 
\end{art}


\begin{thm} \label{local heights}
Let $\lambda(Z):= \lambda_{\hat{D}_0, \dots, \hat{D}_t}(Z)$ be the local height of a $t$-dimensional cycle $Z$ on $X$ with respect to the formally metrized pseudo-divisors $\hat{D}_0, \dots, \hat{D}_t$ satisfying \eqref{suppcond}.
\begin{itemize}
\item[(a)] $\lambda(Z)$ is multilinear and symmetric in the variables $\hat{D}_0, \dots, \hat{D}_t$, and linear in $Z$.
\item[(b)]
For a proper morphism $\varphi:X' \to X$ and a $t$-dimensional cycle $Z'$ on $X'$, we have 
$$\lambda_{\varphi^*\hat{D}_0, \dots,\varphi^*\hat{D}_t}(Z')
=\lambda_{\hat{D}_0,\dots, \hat{D}_t}(\varphi_*Z').$$
\item[(c)] If $\hat{D}_0$ is the pseudo-divisor of a rational function $f$ on $X$ endowed with the trivial metric and if  $Y$ is a representative of $D_1 \dots D_t.Z \in CH_0 \left(|D_1| \cap \dots \cap |D_t| \cap |Z| \right)$, then
$$\lambda(Z)  =   \log|f(Y)|$$
where the right hand side is defined by linearity with respect to points.
\item[(d)] Let $\lambda'(Z)$ be the local height of $Z$ obtained by replacing the metric $\metr$ of $\hat{D}_0$ by another formal metric $\metr'$ on $O(D_0)$. If the formal metrics of $\hat{D}_1, \dots, \hat{D}_t$ are semipositive and if $Z$ is effective, then 
$$|\lambda(Z) - \lambda'(Z)|
\leq d(\metr, \metr') \deg_{O(D_1), \dots, O(D_t)}(Z).$$
\end{itemize}
\end{thm}
\proof This is proved in \cite{Gu2}, \S 9, in case of Cartier divisors. Using the refined intersection theory for formally metrized pseudo-divisors from \cite{Gu3}, \S 5, this can be proved similarly and is included in \cite{Gu3}, Theorem 10.6. \qed 

\begin{art} \rm \label{approximation metrics}
Formal metrics are closed under tensor product and pull-back. However, canoncial metrics of ample symmetric line bundles on an abelian variety are not formal. Hence it is desirable to extend the local heights to a larger class $\gh$ of metrics keeping these properties and including the canonical metrics. 

The tensor product induces a group law on the isometry classes of metrized line bundles on $X^{\rm an}$ which we denote additively by $+$. 
Let $\g$ be the group of isometry classes of formally metrized line bundles on $X$ and let $\gp$ be the submonoid of classes with semipositive metrics.

The {\it completion} $\ghp$ of $\gp$ is the set of isometry classes of line bundles $(L, \metr)$ on $X$ satisfying the following property: For all $n \in \ndop$, there is a proper surjective morphism $\varphi_n:X_n' \rightarrow X$ and a root of a formal metric $\metr_n$ on $\varphi_n^*(L^{\rm an})$ such that $d_{X_n}(\metr_n,\varphi_n^*\metr)\to 0$. Moreover, $\gh:= \ghp - \ghp$ is called the {\it completion} of $\g$.

By the GAGA principle (\cite{Ul}, Theorem 6.8), every formal metric on a projective scheme over $\kdop$ is induced by an algebraic $\kdop^\circ$-model and hence is a quotient of two very ample metrics. By Chow's lemma, we conclude that $\g \subset \gh$ (see \cite{Gu3}, Proposition 10.5, for details). 

Now it is a formal argument to extend the local heights uniquely to $\gh$-pseudo-divisors such that Theorem \ref{local heights} still holds with $\ghp, \gh$ replacing $\gp$ and $\g$ (see \cite{Gu2}, \S 1). 

It would be easier if we could just work with uniform limits of roots of semipositive metrics instead of $\ghp$. Indeed, this would lead to a satisfactory theory of local heights on projective schemes (see \cite{Gu4}, \S 2). For proper schemes over $\kdop$, the coverings $\varphi_n$ used in the definition of $\ghp$ are necessary to apply Chow's lemma, as we have seen above. 
\end{art}

\begin{art} \rm \label{continuous extension}
It is not possible to extend local heights to all continuous metrics using Proposition  \ref{density} because the continuity property (c) in Theorem \ref{local heights} holds only under semipositivity assumptions. 

However, we can define the local height $\lambda(Z)$ with respect to a continuously metrized pseudo-divisor $\hat{D}_0$ and $\gh$-pseudo-divisors $\hat{D}_1,\dots , \hat{D}_t$ satisfying \eqref{suppcond}. Indeed, by linearity, we may assume that $\hat{D}_1,\dots , \hat{D}_t$ have $\ghp$-metrics and that $Z$ is effective. By Proposition \ref{density}, the metric of $\hat{D}_0$ is limit of formal metrics $\metr_n$ on $O(D_0)$ with corresponding pseudo-divisors $\hat{D}_0^{(n)}$. Then
$$\lambda(Z):= \lim_{n \to \infty} \lambda_{\hat{D}_0^{(n)},\hat{D}_1,\dots , \hat{D}_t}(Z)$$
is well-defined by the extension of Theorem \ref{local heights} to $\gh$. Obviously, Theorem \ref{local heights} still holds for these local heights except the symmetry in (a). Then (c) is true also if $\hat{D}_j$ and not $\hat{D}_0$ is induced by $f$, but (d) only holds if we replace the metric on $O(D_0)$ by another continuous metric.

We apply this to the case $D_0=0$. The generalization of Theorem \ref{local heights}(c) shows that the local height $\lambda_{\hat{D}_0,\dots , \hat{D}_t}(Z)$ depends only on $\hat{D}_0$ and the metrized line bundles $(O(D_j)^{\rm an}, \metr_j)_{j=1,\dots,t}$, but not on the choice of the pseudo-divisors.
\end{art}

\begin{art} \rm \label{Chern forms}
A continuous metric $\metr$ on $O_X^{\rm an}$ is given by $\|1\|:= e^{-g}$ for a continuous function $g$ on $X^{\rm an}$. We denote the metric by $\metr_g$.

Let $\overline{L}_1,\dots, \overline{L}_d$ be $\gh$-metrized line bundles on the $d$-dimensional proper scheme $X$ over $\kdop$. For $j=1, \dots , d$, we choose any pseudo-divisor $D_j$ with $O(D_j)=L_j$, e.g. $D_j = (\overline{L_j},X,1)$. For a continuous function $g$ on $X^{\rm an}$, we consider the continuously metrized pseudo-divisor $\hat{O}^g:=(O_X, \metr_g, \emptyset,1)$. Then we define
$$\int_{X^{\rm an}} g \,c_1(\overline{L}_1) \wedge \dots \wedge c_1(\overline{L}_d)
:= \lambda_{\hat{O}^g, \hat{D}_1, \dots, \hat{D}_d}(X).$$
By \ref{continuous extension}, this is independent of the choice of $D_1, \dots, D_d$ and the generalization of Theorem \ref{local heights}  shows that we get a continuous functional on $C(X^{\rm an})$. By the Riesz representation theorem (\cite{Ru2}, Theorem 6.19), $c_1(\overline{L}_1) \wedge \dots \wedge c_1(\overline{L}_d)$ is a regular Borel measure on $X^{\rm an}$.

These measures were first introduced by Chambert-Loir (see \cite{Ch}) through a slightly different approach and under the additional assumptions that $\kdop$ contains a countable subfield and that $X$ is projective.
\end{art}

\begin{cor} \label{Chern properties}
For $\gh$-metrized line bundles $\overline{L}_1,\dots, \overline{L}_d$ on the $d$-dimensional proper scheme $X$ over $\kdop$, the following properties hold:
\begin{itemize}
\item[(a)] $c_1(\overline{L}_1) \wedge \dots \wedge c_1(\overline{L}_d)$ is multilinear and symmetric in $\overline{L}_1,\dots, \overline{L}_d$.
\item[(b)] If $\varphi:X' \rightarrow X$ is a morphism of integral proper schemes over $\kdop$, then
$$\varphi_* \left( c_1(\varphi^* \overline{L}_1) \wedge \dots \wedge c_1(\varphi^* \overline{L}_d) \right) = \deg(\varphi) c_1(\overline{L}_1) \wedge \dots \wedge c_1(\overline{L}_d).$$
\item[(c)] If the metrics of $\overline{L}_1,\dots, \overline{L}_d$ are in $\ghp$, then
$$\left| \int_{X^{\rm an}} g \,c_1(\overline{L}_1) \wedge \dots \wedge c_1(\overline{L}_d) \right| \leq |g|_{\rm sup} \deg_{L_1, \dots , L_d}(X)$$
for all $g \in C(X^{\rm an})$.
\end{itemize}
\end{cor}

\proof These properties follow immediately from the corresponding properties of the generalization of Theorem \ref{local heights} mentioned in \ref{continuous extension}. \qed

\begin{rem} \rm \label{base change}
Let $\kdop'$ be an algebraically closed extension of $\kdop$ endowed with a complete absolute value extending $|\phantom{a}|$. Obviously, the local heights are invariant under base change to $\kdop'$. If $\pi: X_{\kdop'} \rightarrow X$ denotes the natural projection, then we deduce
$$\pi_* \left( c_1(\pi^* \overline{L}_1) \wedge \dots \wedge c_1(\pi^* \overline{L}_d) \right) = c_1(\overline{L}_1) \wedge \dots \wedge c_1(\overline{L}_d).$$
\end{rem}

\begin{prop} \label{Dirac decomposition} 
Let $\overline{L}_1,\dots, \overline{L}_d$ be formally metrized line bundles on the complete variety $X$ over $\kdop$ of dimension $d$. Then there is a formal $\kdop$-model $\Xcal$ of $X$ with reduced special fibre and for every $j \in \{0, \dots , d\}$ a formal $\kdop^\circ$-model $\Lcal_j$ of $L_j$ on $\Xcal$ inducing the metric of $\overline{L}_j$. For such models, we have always
$$c_1(\overline{L}_1) \wedge \dots \wedge c_1(\overline{L}_d)= 
\sum_Y\deg_{\tilde{\Lcal}_1, \dots, \tilde{\Lcal}_d}(Y) \delta_{\xi_Y},$$
where $Y$ ranges over the irreducible components of $\tilde{\Xcal}$ and $\delta_{\xi_Y}$ is the Dirac measure in $\xi_Y$.
\end{prop}

\proof The existence of such formal $\kdop^\circ$-models follows from \cite{Gu2}, 8.1.  To show equality of the regular Borel measures, it is enough to show that their integrals agree on the  subset of $C(X^{\rm an})$ induced by formal metrics on $O_X^{\rm an}$ (see Proposition \ref{density}). If $\metr_g$ is such a metric with formal model denoted by $\Ocal(g)$, then the section $1$ of $O_X^{\rm an}$ induces a meromorphic section of $\Ocal(g)$. By the very definition of multiplicities (see \cite{Gu2}, \S 3), the corresponding  divisor is vertical and has multiplicity $g(\xi_Y)$ in the irreducible components $Y$ of $\tilde \Xcal$. By definition of divisoral intersections on $\Xcal$ (\cite{Gu2}, \S 4), this leads to the claim. \qed 

\begin{prop} \label{Chern continuity} 
Let $L_1, \dots , L_d$ be line bundles on the $d$-dimensional proper scheme $X$ over $\kdop$. Let $S_j^+$ be the set of $\ghp$-metrics on $L_j^{\rm an}$ endowed with the distance from \ref{formal metrics}. Then  we have a continuous map from $S_1^+ \times \dots \times S_d^+$ to the space of regular Borel measures on $X^{\rm an}$ endowed with the weak topology, given by $(\metr_1, \dots, \metr_d) \mapsto c_1(\overline{L}_1) \wedge \dots \wedge c_1(\overline{L}_d)$. Moreover, $c_1(\overline{L}_1) \wedge \dots \wedge c_1(\overline{L}_d)$ is positive and $X^{\rm an}$ has measure $\deg_{L_1, \dots,L_d}(X)$.
\end{prop}

\proof Recall that the weak topology on the set of regular Borel measures of $X^{\rm an}$ is the coarsest topology such that the map $\mu \mapsto \int f\, \mu$ is continuous for every $f \in C(X^{\rm an})$. By the Riesz representation theorem (\cite{Ru2}, Theorem 6.19), the dual of the Banach space $C(X^{\rm an})$ is isometric to the space $M_{\rm reg}(X^{\rm an})$ of regular Borel measures on $X^{\rm an}$ endowed with the variation as norm. By a standard fact of functional analysis (\cite{Ru1}, Theorem 4.3), we deduce that every closed ball in $M_{\rm reg}(X^{\rm an})$ is compact in the weak topology.

To prove the proposition, we may assume that $X$ is an irreducible variety. Let us consider $\ghp$-metrized line bundles $\overline{L}_1 = (L_1, \metr_1), \dots, (\overline{L}_d, \metr_d)$.

{\bf First step: \/} {\it For $n \in \ndop$, let $\varphi_n:X_n \rightarrow X$ be a generically finite surjective morphism of irreducible complete varieties over  $\kdop$. For $j=1, \dots, d$, let $\metr_{j,n}$ be a $\hat{\frak g}_{X_n}^+$-metric on $\varphi^*L_j$ such that $d_{X_n}(\metr_{j,n},\varphi_n^*\metr_j) \to 0$ for $n \to \infty$. We assume that
$$\left\{\mu_n := (\varphi_n)_*\left( c_1(\varphi_n^*L_1,\metr_{1,n}) \wedge \dots
\wedge  c_1(\varphi_n^*L_d,\metr_{d,n}) \right) \mid n \in \ndop \right\}$$
is bounded in $M_{\rm reg}(X^{\rm an})$. Then $\mu_n$ converges weakly to $\mu:=c_1(\overline{L}_1) \wedge \dots \wedge c_1(\overline{L}_d)$.}

The proof of the first step is by contradiction. By passing to a subsequence and using weak compactness of closed balls, we may assume that $\mu_n$ converges weakly to a regular Borel measure $\mu_\infty \neq \mu$. By Proposition \ref{density}, there is a formal metric $\metr_g$ on $O_X^{\rm an}$ such that 
\begin{equation} \label{neq}
\int g \, d\mu_\infty \neq \int g \, d\mu.
\end{equation}
There is a line bundle $L_0$ on $X$ with $\ghp$-metrics $\metr_{\pm}$ such that $\metr_g = \metr_+/\metr_-$. It is easy to construct non-zero meromophic sections $s_j$ of $L_j$ $(j=0, \dots, d)$ such that $\cap_j \supp(D_j)= \emptyset$ for $D_j:=\Div(s_j)$. Let $\lambda^{\pm}(X)$ be the local heights with respect to $D_0, \dots, D_d$ endowed with $\metr_{\pm}, \metr_1, \dots, \metr_d$ and let $\lambda_n^{\pm}(X_n)$ be the local heights with respect to $\varphi_n^*(D_0),\dots,\varphi_n^*(D_d)$ endowed with $\varphi_n^*\metr_{\pm}, \metr_{1,n}, \dots, \metr_{d,n}$. By Theorem \ref{local heights} generalized to $\ghp$ (see \ref{approximation metrics}), we have
$$\lim_{n \to \infty} \frac{1}{\deg(\varphi_n)} \lambda_n^{\pm}(X_n) = \lambda^{\pm}(X).$$
If we subtract these two formulas, then we get a contradiction to \eqref{neq}. This proves the first step.

If we use Theorem \ref{local heights}(c) for a constant $f$ (and for $\ghp$-metrics as in \ref{approximation metrics}), then we get $\mu(X^{\rm an})=\deg_{L_1, \dots,L_d}(X)$. We claim that $\mu$ is positive. Using Proposition \ref{Dirac decomposition}, this holds for roots of formal metrics. The corresponding measures have norm  $\deg_{L_1, \dots,L_d}(X)$. In general, the metrics used in the approximation process for $\metr_j$ (see the definition of $\ghp$ in \ref{approximation metrics}) may be chosen as in the first step (see \cite{Gu3}, Remark 10.3). Boundedness of $\{\mu_n \mid n \in \ndop \}$ follows from the special case of roots of formal metrics. Then the first step yields positivity of $\mu$. 

Finally, we have to show that $\mu$ is continuous in $(\metr_1, \dots, \metr_d) \in \ghp$. By positivity and $\mu(X^{\rm an})=\deg_{L_1,\dots,L_d}(X)$, every such measure has norm $\deg_{L_1,\dots,L_d}(X)$ and continuity follows also from the first step. \qed

\begin{art} \rm \label{canonical metrics}
Now let $(L,\rho)$ be a rigidified line bundle on the abelian variety $A$ over $\kdop$, i.e. $\rho \in L(\kdop)\setminus \{0\}$. Then there is a {\it canonical metric} $\metr_\rho$ for $(L,\rho)$ which behaves well with respect to tensor product and homomorphic pull-back (see \cite{BG}, Theorem 9.5.7). The construction is by Banach's fixed point theorem and yields that $\metr_\rho$ is a $\ghp$-metric if $L$ is ample and symmetric (see proof of Theorem 9.5.4 in \cite{BG}). 
\end{art}

\begin{rem} \rm \label{cohomological semipositive metrics}
For odd line bundles, we can't be sure that we get $\hat{\frak g}_A$-metrics. Using in $\gp$ cohomologically semipositive metrics instead of the smaller class of semipositive formal metrics, we get a larger class $\gp$ with the same properties. Since canonical metrics of odd line bundles are cohomologically semipositive, the new class $\gh$ includes all canonical metrics on abelian varieties (see \cite{Gu3}, \S 10). All  results of this section hold for this $\gh$ as well.

Let $X$ be a smooth complete variety and assume that one $L_j$ is algebraically equivalent to $0$ endowed with a canonical metric, i.e. pull-back of a canonical metric from the Picard variety. Then Theorem \ref{local heights}, applied to the above metrics, shows that the local height does not depend on the metrics of the other line bundles. In particular, we deduce
$$c_1(\overline{L}_1) \wedge \dots \wedge c_1(\overline{L}_d) = 0.$$
\end{rem}

\begin{art} \label{can measure} \rm 
Let $X$ be a closed subvariety of $A$ of dimension $d$ and let $\overline{L}_1, \dots ,\overline{L}_d$ be canonically metrized line bundles on $A$. Then $\mu := c_1(\overline{L}_1|_X) \wedge \dots \wedge c_1(\overline{L}_d|_X)$ is called a {\it canonical measure} on $X$. Note that the canonical metric is only determined up to $|\kdop^\times|$-multiples by the line bundle. By Theorem \ref{local heights}(c), the canonical measure $\mu$ does not depend on the choice of the canonical metrics. The same argument as in \ref{cohomological semipositive metrics} shows that $\mu=0$ if one line bundle is odd.
\end{art}

\section{Polytopal domains in $\Tor$}

Our goal is to study the formal properties of certain affinoid domains in $\Tor$ associated to polytopes in $\rdop^n$. They are related to Mumford's construction of models of totally degenerate abelian varieties discussed in Section 6. For the terminology used from convex geometry, the reader is referred to the appendix.


Recall that $\Gamma$ is the value group of the valuation $v := - \log |\phantom{a}|$ on $\kdop$.
On $\Tor$, we always fix coordinates $x_1, \dots, x_n$. Then we have a continuous map
$$\val:(\Tor)^{\rm an}_\kdop \longrightarrow \rdop^n, \quad p \mapsto 
\left( -\log p(x_1), \dots, -\log p(x_n) \right).$$

A large part of the following result is contained in \cite{EKL}, 3.1. We give here a different proof which will be used later.

\begin{prop} \label{polyhedral affinoid domains}
Let $\Delta$ be a $\Gamma$-rational polytope. 
Then the set of Laurent series
$$\kdop \langle U_\Delta \rangle := \left \{ \sum_{\mb \in \zdop^n} a_\mb x_1^{m_1} \dots x_n^{m_n} \mid \lim_{|\mb| \to \infty} v(a_\mb) + \mb \cdot \ub =\infty \, \; \forall \ub \in \Delta \right\}$$
is the $\kdop$-affinoid algebra of the Weierstrass domain $U_\Delta := \val^{-1}(\Delta)$ of $(\Tor)^{\rm an}_{\kdop}$. It has supremum norm
\begin{equation} \label{polytopal supnorm}
\left| \sum_{\mb \in \zdop^n} a_\mb \xb^\mb \right|_{\sup} = \sup_{\ub \in \Delta,\mb \in \zdop^n} |a_\mb|e^{-\mb \cdot \ub} = \max_{\text{$\ub$ \rm vertex, $\mb \in \zdop^n$}}|a_\mb|e^{-\mb \cdot \ub}.
\end{equation}
\end{prop}
\proof For $\ub \in \rdop^n$, the polyannulus $\val^{-1}(\ub)$ has multiplicative supremum norm
\begin{equation} \label{u norm}
\left| \sum_{\mb \in \zdop^n} a_\mb \xb^\mb \right|_\ub := \max_{\mb \in \zdop^n} |a_\mb|e^{-\mb \cdot \ub}
\end{equation}
which proves the claim for $\Delta=\{\mathbf 0\}$ (see \cite{BGR}, 6.1.4). In general, $\Delta$ is defined by
$$\Delta = \bigcap_{\mb \in S} \{\ub \in \rdop^n \mid \mb \cdot \ub + v(b_\mb) \geq 0 \}$$
for a finite $S \subset \zdop^n$ and suitable $b_\mb \in \kdop^\times$. We conclude that $U_\Delta$ is the Weierstrass domain in $(\Tor)^{\rm an}_{\kdop}$ given by $|b_\mb \xb^\mb| \leq 1$, $\mb \in S$. By \cite{BGR}, Proposition 6.1.4/2, we deduce that every analytic function $f$ on $U_\Delta$ has a Laurent series expansion $\sum a_\mb \xb^\mb$ on $U_\Delta$. If $\ub$ is a vertex of $\Delta$, then $\val^{-1}(\ub)$ is a Weierstrass domain in $U_\Delta$ and we get $|f|_\ub \leq |f|_{\rm sup}$. Since $\ub \cdot \mb + v(a_\mb)$ achieves its minimum on $\Delta$ always in a vertex $\ub$, we get
$$\sup_{\ub \in \Delta} |f|_\ub \leq |f|_{\rm sup}.$$
By the ultrametric triangle inequality, we have equality  and we deduce easily the claim. \qed

\begin{art} \label{toric variety} \rm 
A $\Tor$-{\it toric variety} over $\tilde{\kdop}$ is a normal variety $Y$ over $\tilde{\kdop}$ with an algebraic $(\Tor)_{\tilde \kdop}$-action  containing a  dense $n$-dimensional orbit. 

The theory of toric varieties will be  very important in the sequel, we refer to \cite{KKMS}, \cite{Fu2} or \cite{Oda} for details. There are bijective correspondences between
\begin{itemize}
\item[(a)] rational polyhedral cones $\sigma$ in $\rdop^n$ which do not contain a linear subspace $\neq \{0\}$,
\item[(b)] finitely generated saturated semigroups $S$ in $\zdop^n$ which generate $\zdop^n$ as a group,
\item[(c)] affine $\Tor$-toric varieties $Y$ over $\tilde{\kdop}$ (up to equivariant isomorphisms).
\end{itemize}
The correspondences are given by $S=\check{\sigma}\cap \zdop^n$ and $Y=\Spec(\tilde{\kdop}[\tilde{\xb}^S])$, where ${\tilde{\xb}}^S := \{{\tilde{\xb}}^\mb \mid \mb \in S\}$ for the coordinates ${\tilde{\xb}}$ on $(\Tor)_{\tilde{\kdop}}$.
\end{art}

\begin{art} \rm \label{torus action}
Let $\Delta$ be still a $\Gamma$-rational polytope in $\rdop^n$. For $\mb \in \zdop^n$, there is $b_\mb \in \kdop^\times$ with $v(b_\mb)=-\min_{\ub \in \Delta} \ub \cdot \mb$. Note that $y_\mb:= b_\mb \xb^\mb$ has supremum norm $1$ on $U_\Delta$. We denote by $\pi:U_\Delta \rightarrow \tilde{U}_\Delta$ the reduction map.

The affinoid torus ${\mathbb T}_1^{\rm an}:= \{p \in (\Tor)^{\rm an}_\kdop \mid p(x_j)=1, \, j=1, \dots,n\}$ acts on $U_\Delta$. Passing to reductions, we get a torus action of $({\mathbb T}_1^{\rm an})\sptilde = (\Tor)_{\tilde{ \kdop}}$ on $\tilde{U}_\Delta$. 
\end{art}

In the framework of  algebraic geometry and for a discrete valuation, the following result is due to Mumford (\cite{Mu}, \S 6). We give here an analytic formulation over $\kdop$ without assuming that the valuation is discrete.

\begin{prop} \label{torus orbits 1}
The following properties hold for $U_\Delta=\val^{-1}(\Delta)$:
\begin{itemize}
\item[(a)] The elements $\left(y_\mb\right)_{\mb \in \zdop^n}$ generate a dense subalgebra of $\kdop \langle U_\Delta \rangle ^\circ$.
\item[(b)] There is a bijective order reversing correspondence between torus orbits $Z$ of $\tilde{U}_\Delta$ and open faces $\tau$ of $\Delta$, given by
$Z_\tau= \pi(\val^{-1}(\tau))$ and $\tau_Z=\val(\pi^{-1}(Z))$.
\item[(c)] $\dim(\tau)+\dim(Z_\tau)=n$.
\item[(d)] If $Y_{\bf u}$ is the irreducible component of $\tilde{\Xcal}$ corresponding to the vertex $\ub= \val(\xi_Y)$ of $\Delta$ by (b), then the natural $(\Tor)_{\tilde{\kdop}}$-action of $\tilde{U}_\Delta$ makes $Y_\ub$ into an affine toric variety. The corresponding rational polyhedral cone is generated by $\Delta - \ub$.
\item[(e)] If $\Delta'$ is also a $\Gamma$-rational polytope with $\Delta' \subset \Delta$, then the canonical morphism $U_{\Delta'} \rightarrow U_\Delta$ induces an open immersion of the reductions if and only if $\Delta'$ is a closed face of $\Delta$.
\end{itemize}
\end{prop}

\proof Property (a) follows easily from \eqref{polytopal supnorm}. For every $\ub \in \Delta$, $|\phantom{a}|_\ub$ from \eqref{u norm} restricts to a multiplicative norm
on $\kdop \langle U_\Delta \rangle$ which is bounded by the supremum norm. Hence it may be viewed as a point $\xi_\ub \in U_\Delta$.

The {\it Shilov boundary} is the unique minimal set $\Theta$ of $U_\Delta$ such that every $f \in \kdop \langle U_\Delta \rangle$ has its maximum in $\Theta$. By \cite{Ber}, Proposition 2.4.4, $\Theta$ is the set of $\xi_Y \in U_\Delta$ corresponding to the irreducible components $Y$ of $\tilde{U}_\Delta$ by \ref{formal analytic varieties}. Using \eqref{polytopal supnorm}, we get $\Theta= \{|\phantom{a}|_\ub \mid \text{$\ub$ vertex of $\Delta$}\}$. Note that the vertex $\ub$ corresponding to $Y$ is given by $\ub=\val(\xi_Y)$. By definition of $\xi_Y$, we have 
$$\tilde{\kdop}[Y] = \kdop \langle U_\Delta \rangle^\circ / \{|\phantom{a}|_\ub < 1\}.$$
To prove (d), it is enough to show that $\tilde{\kdop}[Y]$ is isomorphic to $\tilde{\kdop}[{\tilde{\yb}}^{\check{\sigma} \cap \zdop^n}]$ for the cone $\sigma$ generated by $\Delta - \ub$. By a change of coordinates, we may assume that $\ub = \mathbf 0$. For $S:= \{ \mb \in \zdop^n \mid v(b_\mb)=0 \}$, (a) yields that $\tilde{\kdop}[Y]$ is generated by $(\tilde{y}_{\mb})_{\mb \in S}$. By construction, we have $S= \check{\sigma} \cap \zdop^n$. It remains to show that a relation 
$$\tilde{p}(\tilde{y}_{\mb_1}, \dots , \tilde{y}_{\mb_r}) = \tilde{0} \in \tilde{\kdop}[Y]$$
comes from a relation in $\tilde{\kdop}[{\tilde{\yb}}^S]$. For simplicity, we may assume $b_\mb =1$ for all $\mb \in S$. Then we have $|p(y_{\mb_1}, \dots, y_{\mb_r})|<1$ on $U_\Delta$. Repacing $y_\mb$ by $\xb^\mb$ and collecting terms of the same degree, we get the desired relation in $\tilde{\kdop}[{\tilde{\yb}}^S]$. This proves (d).

To prove (b), let $\tau$ be an open face of $\Delta$. There is $I \subset \zdop^n$ such that $\tau$ is given by 
\begin{equation} \label{open face}
\ub \cdot \mb + v(b_\mb) 
\begin{cases}=0 & \text{if $\mb \in I$},\\
>0 & \text{if $\mb \in \zdop^n \setminus I$}.
\end{cases}
\end{equation}
Then $\tilde{x} \in Z_\tau:=\pi(\val^{-1}(\tau))$ if and only if $\tilde{y}_\mb(\tilde{x}) \neq \tilde{0}$ for $\mb \in I$ and $\tilde{y}_\mb(\tilde{x}) = \tilde{0}$ for $\mb \not \in I$. 

We prove first that $Z_\tau$ is a torus orbit. We choose a vertex $\ub$ of $\overline{\tau}$ with associated toric variety $Y_\ub$. By a change of coordinates, we may assume again that $\ub=\mathbf 0$. It follows from the proof of (d) that $Y_\ub$ is given by the equations $\tilde{y}_\mb= \tilde{0}$ for $\mb \not \in S$. Since $I \subset S$, we conclude that $Z_\tau \subset Y_\ub$. Then $Z_\tau$ is given in $Y_\ub$ by the equations $\tilde{y}_\mb=\tilde{0}$ for $\mb \in S \setminus I$ and the inequalities $\tilde{y}_\mb \neq \tilde{0}$ for $\mb \in I$ and hence $Z_\tau$ is a torus orbit (see \cite{Fu2}, 3.1).

Since $\pi$ is surjective, all torus orbits have  this form. The above characterization of $Z_\tau$ shows that we get a bijective correspondence in (b). Moreover, $\tau \subset \val(\pi^{-1}(Z_\tau))$ is also clear and equality follows as the torus orbits (resp. open faces) form a disjoint covering of  $\tilde{U}_\Delta$ (resp. $\Delta$).

Finally, (c) and (e) follow from the theory of toric varieties applied to $Y_\ub$ for a vertex $\ub$ of $\overline{\tau}$ (resp. $\Delta'$). \qed

\begin{cor}  \label{u points}
For $\ub \in \Delta$, let $\xi_\ub$ be the point of $U_\Delta$ defined by \eqref{u norm}. Then $\tilde{\xi}_\ub$ is the generic point $\zeta$ of the torus orbit $Z_\tau$ associated to the unique open face $\tau$ containing $\ub$. Moreover, $\xi_\ub$ is the Shilov boundary of $\val^{-1}(\ub)$.
\end{cor}
\proof By Proposition \ref{torus orbits 1}(b), we have $\tilde{\xi}_\ub \in Z_\tau$. Let $f = \sum a_\mb \xb^\mb \in \kdop \langle U_\Delta \rangle^\circ$ with $\tilde{f}(\tilde{\xi}_\ub)= \tilde{0}$. Note that \eqref{open face} determines $I$ and describes the open face $\tau$. Using $f \in \{|\phantom{a}|_\ub < 1 \}$,  we have $|a_\mb | < e^{\mb \cdot \ub} = |b_\mb|$ for all $\mb \in I$. By the description of $Z_\tau$ in the proof of Proposition \ref{torus orbits 1}, $\tilde{f}$ vanishes identically on $Z_\tau$ and hence $\tilde{\xi}_\ub = \zeta$. We have also seen that the vertices of $\Delta$ correspond to the Shilov boundary of $U_\Delta$. Using $\{\ub\}$ instead of $\Delta$, we get the last claim. \qed

\begin{art} \label{globalization} \rm
We globalize our considerations. Let $\Ccal$ be a $\Gamma$-rational polytopal complex in $\rdop^n$. By Proposition \ref{torus orbits 1}, it is easy to deduce that $(U_\Delta)_{\Delta \in \Ccal}$ is a formal analytic atlas of $U = \cup_{\Delta \in \Ccal} U_\Delta$. The associated admissible formal scheme $\Ucal$ over $\kdop^\circ$ (see \ref{admissible formal schemes}) has a formal open affine atlas
$$\Ucal_\Delta := \Spf \left( \kdop \langle U_\Delta \rangle^\circ \right), \quad \Delta \in \Ccal.$$
Clearly, the affinoid torus ${\mathbb T}_1^{\rm an}$ acts on $U$, ${\mathbb T}_1 := \Spf( \kdop \langle \val^{-1}(\mathbf 0) \rangle^\circ )$ acts on $\Ucal$ and $(\Tor)_{\tilde{\kdop}}$ acts on the special fibre $\tilde{\Ucal}$. Proposition \ref{torus orbits 1} yields the following result:
\end{art}

\begin{prop} \label{torus orbits global}
Under the hypothesis above, we have the following properties of $\Ucal$:
\begin{itemize}
\item[(a)] There is a bijective correspondence between torus orbits of $\tilde{\Ucal}$ and $\{\relint(\Delta) \mid \Delta \in \Ccal\}$.
\item[(b)] The irreducible components of $\tilde{\Ucal}$ are in bijective correspondence with vertices of $\Ccal$.
\item[(c)] If $Y_\ub$ is the irreducible component of $\tilde{\Ucal}$ corresponding to the vertex $\ub$, then $Y_\ub$ is a toric variety with fan given by the cones $\sigma$ in $\rdop^n$ which are generated by $\Delta - \ub$ for $\Delta \in \Star(\ub)$.
\item[(d)] For $\Delta, \Delta' \in \Ccal$, $\Ucal_{\Delta'}$ is an open subset of $\Ucal_\Delta$ if and only if $\Delta'$ is a closed face of $\Delta$.
\end{itemize}
\end{prop}

\begin{rem} \label{polytopal toric variety} \rm 
Recall from the appendix that $\relint(\Delta)$ denotes the relative interior of $\Delta$. Every $\Delta \in \Ccal$ induces a toric variety $Y_\Delta$, given as the closure of the torus orbit associated to $\relint(\Delta)$. Let $\ldop_\Delta$ be the linear space in $\rdop^n$ generated by $\Delta - \ub$, $\ub \in \Delta$, and let $N_\Delta := \ldop_\Delta \cap \zdop^n$. Then the subtorus $H_\Delta$ of $(\Tor)_{\tilde{\kdop}}$, given by $H_\Delta(\tilde{\kdop}) = N_\Delta \otimes_\zdop \tilde{\kdop}^\times$, acts trivially on $Y_\Delta$ and hence $Y_\Delta$ is a toric variety with respect to the torus $T_\Delta = (\Tor)_{\tilde{\kdop}}/H_\Delta$. If we project the cones in $\rdop^n$ generated by some $\Delta' - \ub$, $\Delta' \in \Star(\Delta)$, to $\rdop^n/ \ldop_\Delta$, then we get the fan associated to $Y_\Delta$. For details, we refer to \cite{Fu2}, 3.1. 
\end{rem}
\begin{art} \label{subdivision} \rm
Let $\Dcal$ be also a $\Gamma$-rational  polytopal complex in $\rdop^n$ which subdivides $\Ccal$. Then the atlas $(U_\sigma)_{\sigma \in \Dcal}$ yields a formal analytic structure on $U$ which is finer than $\Ucal^{\rm f-an}$. We denote the associated formal scheme over $\kdop^\circ$ by $\Ucal'$ and we get a morphism $\Ucal' \rightarrow \Ucal$.
\end{art}

\begin{prop} \label{toric and etale}
With $\Ucal$ and $\Ucal'$ as in \ref{subdivision}, let $\varphi: \Xcal \rightarrow \Ucal$ be an \'etale morphism of admissible formal schemes over $\kdop^\circ$.
\begin{itemize}
\item[(a)] The base change $\varphi': \Xcal' \rightarrow \Ucal'$ of $\varphi$ to $\Ucal'$ is \'etale.
\item[(b)] The reduction $\tilde{\varphi}'$ maps an irreducible component $Y'$ of $\tilde{\Xcal}'$ dominantly to a unique irreducible component $Y$ of $\tilde{\Ucal}'$.
\item[(c)] Let $Y$ be the irreducible component of $\tilde{\Ucal}'$ associated to the vertex $\ub$ in $\Dcal$ and let $Z$ be the torus orbit of $\tilde{\Ucal}$ corresponding to the unique open face $\relint(\Delta)$, $\Delta \in \Ccal$, containing $\ub$ (see Proposition \ref{torus orbits global}). Then 
$$\sum_{Y'} [\tilde{\kdop}(Y'):\tilde{\kdop}(Y)] = \sum_{Z'} [\tilde{\kdop}(Z'):\tilde{\kdop}(Z)]$$
where $Y'$ (resp. $Z'$) ranges over all irreducible components of $(\tilde{\varphi}')^{-1}(Y)$ (resp. $\tilde{\varphi}^{-1}(Z)$).
\item[(d)] If $\dim(\Delta)=n$ in (c), then $\tilde{\varphi}$ maps all these $Y'$ isomorphically onto the toric variety $Y$. Then $Z$ is  a closed point and the number of such $Y'$ is equal to $\card(\tilde{\varphi}^{-1}(Z)) < \infty$.
\end{itemize}
Note that the right hand side in formula (c) and the cardinality in (d) depend only on $\varphi$ and not on the subdivision $\Dcal$.
\end{prop}
\proof Since $\varphi$ is \'etale, (a) and (b) are obvious and are true in much more generality. By Corollary \ref{u points}, the generic point $\tilde{\xi}_\ub$ of $Y$ maps to the generic point $\zeta$ of $Z$ with respect to $\tilde{\Ucal}' \rightarrow \Ucal$, hence formula (c) follows from the fact that the degree of the fibre $\tilde{\varphi}^{-1}(\zeta)$ is invariant under base change to $\tilde{\kdop}(Y)$.

In (d), Proposition \ref{torus orbits 1}(c) yields that $Z$ is a closed point of $\tilde{\Ucal}$. Since $\tilde{\varphi}$ is \'etale, $\tilde{\varphi}^{-1}(Z)$ is the disjoint open union of its points $Z'$. Similarly, $(\tilde{\varphi}')^{-1}(Y)$ is the disjoint open union of its components $Y'$. We conclude that $(\tilde{\varphi}')^{-1}(Y)$ is the base change of $\tilde{\varphi}^{-1}(Z)$ to $Y$  and hence there is exactly one $Y'$ over $Z'$. Moreover, it is isomorphic to $Y$. This proves (d). \qed

\section{Tropical analytic geometry}

We study the analytic analogue of tropical varieties for the  polytopal domains in $\Tor$ considered in the previous section.  We generalize the basic results of Einsiedler, Kapranov and Lind \cite{EKL} to analytic subvarieties of such a domain. In this analytic setting, the Bieri--Groves set from \cite{EKL} (see also \cite{BiGr}) may be strictly larger than the tropical variety and we have to use different methods from analytic and formal geometry.


\begin{art} \label{tropical analytical variety} \rm 
Let $\Delta$ be a $\Gamma$-rational polytope in $\rdop^n$ and
let $X$ be a closed analytic subvariety   of $U_\Delta$. By continuity of $\val$, we conclude that  $\val(X)$ is a compact subset of $\rdop^n$. Note that $X$ is given by an ideal $I$ of $\kdop \langle U_\Delta \rangle$ and it is endowed with the structure of an analytic space. If $X$ is connected, then $\val(X)$ is connected.

For a closed subscheme $X$ of $\Tor$ defined over $\kdop$, the closure of $\val(X(\kdop))$ is equal to $\val(X^{\rm an})$. This is clear by density, hence $\val(X^{\rm an})$ is equal to the usual definition of a tropical variety.
\end{art}  

\begin{prop} \label{polyhedron}
Let $X$ be a closed analytic subvariety of $U_\Delta$. Then $\val(X)$ is a $\Gamma$-rational  polytopal set in $\Delta$.
\end{prop}

This is a fundamental result and follows from  \cite{Ber5}, Corollary 6.2.2. Nonetheless, we give a proof under the assumption that $X$ has a semistable alteration. The proof in the special case is very instructive and will be used later.

\proof In fact, we assume only that there is a dominant morphism $f:X' \rightarrow X$, where $X'$ is the generic fibre of a quasicompact strictly semistable formal scheme $\Xcal'$ over $\kdop^\circ$. This assumption is satisfied in all cases relevant for our applications. Indeed, let $\kdop$ be the completion of the algebraic closure of a field $K$ with $v|_K$ a  complete discrete valuation. If $X$ is a smooth compact analytic space, then Hartl proved the existence of such an $f$ which is even \'etale and surjective (\cite{Ha}, Corollary 1.5). Together with de Jong's alteration theorem (\cite{dJ}, Theorem 4.1), we conclude that such an $f$ exists if $X$   is an analytic subdomain of an algebraic variety. 

By Proposition \ref{semistable properties}, $\Xcal'$ is covered by finitely many formal open affine subsets $\Ucal'$ allowing an \'etale morphism 
$$\psi: \Ucal' \longrightarrow \Scal :=\Spf \left( \kdop^\circ \langle x_0', \dots , x_m' \rangle /  \langle x_0' \cdots x_r' - \pi \rangle \right)$$
such that every irreducible component of the special fibre $\tilde{\Ucal}'$ passes through a closed point $\tilde{P} \in \tilde{\Ucal}'$ with $\tilde{\psi}(\tilde{P}) = \tilde{\mathbf 0}$. Moreover, for $j=1,\dots,n$, $f^*(x_j)$ is a unit on the generic fibre $U'$ of $\Ucal'$ and hence
\begin{equation} \label{semistable f-coord}
f^*(x_j)= \lambda_j u_j \psi^*(x_1')^{m_{j1}} \cdots \psi^*(x_r')^{m_{jr}}
\end{equation}
for suitable $\lambda_j \in \kdop^\times$, $u_j \in \Ocal(\Ucal')^\times$ and $\mb_j \in \zdop^r$. Let $\Sigma(r,\pi):=\{\ub' \in \rdop_+^r \mid u_1' + \dots + u_r' \leq v(\pi) \}$ be the ``standard simplex'' in $\rdop^r$ and let $f_{\rm aff}^{(0)}: \rdop^r \rightarrow \rdop^n$ be the affine map given by 
$$f_{\rm aff}^{(0)}(\ub')_j := \mb_j \cdot \ub' + v(\lambda_j) \quad (j=1, \dots, n).$$
By density of $f(X')$ in $X$, it is enough to prove 
\begin{equation} \label{simplex image}
\val(f(U')) = f_{\rm aff}^{(0)}(\Sigma(r,\pi)).
\end{equation}
Since $|u_j|=1$ on $U'$, the inclusion ``$\subset$'' is obvious. Proposition \ref{implicit function} yields $U_+'(\tilde{P}) \cong \Scal_+^{\rm an}(\tilde{\mathbf 0})$ for formal fibres. We conclude that
$$\val(f(U_+'(\tilde{P}))) = f_{\rm aff}^{(0)}(\relint(\Sigma(r,\pi))).$$
The right hand side is dense in $f_{\rm aff}^{(0)}(\Sigma(r,\pi))$ and we get equality in \eqref{simplex image}. \qed

\vspace{2mm}
For the following two results, we assume that $\Delta$ is an $n$-dimensional $\Gamma$-rational polytope in $\rdop^n$. The interior of $\Delta$ is denoted by $\Int(\Delta)$.
\begin{prop} \label{totally concave}
For a closed analytic subvariety $X$ of $U_\Delta$, the polytopal set $\val(X)$ is concave in all the points of  ${\Int(\Delta)}$.
\end{prop}

\proof Let $\tilde{X}$ be the reduction of $X$ (see \ref{Berkovich spectrum}). By \cite{BGR}, Theorem 6.3.4/2, the morphism $\tilde{X} \rightarrow \tilde{U}_\Delta$ is finite. By Proposition \ref{torus orbits 1}, the reduction of $\val^{-1}({\Int(\Delta)})$ in $\tilde{U}_\Delta$ is just one $\tilde{\kdop}$-rational point $\tilde{P}$ (the zero-dimensional torus orbit). The image of $X^\circ := X \cap \val^{-1}({\Int(\Delta)})$ under the reduction  map $X \rightarrow \tilde{X}$ is lying over $\tilde{P}$ with respect to the above finite morphism. 
We conclude that the reduction of $X^\circ$ consists of finitely many closed points $\tilde{x}_1, \dots, \tilde{x}_r$ in $\tilde{X}$. Since the inverse image of a closed point with respect to the reduction map is open in $X$, the topological space $X^\circ$ decomposes into disjoint open and closed subsets $V_j$ lying over $\tilde{x}_j$. 

We have to show that  $\val(X)$ is  concave in $\ub_0 \in \val(X^\circ)$. By Proposition \ref{polyhedron} and \ref{convexe geometry}, there is a $\Gamma$-rational polytopal decomposition of $\val(X)$. Let $\sigma$ be the polytope of this decomposition such that $\ub_0$ is contained in the relative interior of $\sigma$. All  $\ub \in \relint(\sigma)$ have the same local cone ${\rm LC}_{\ub}(\val(X))$. The points with coordinates in the value group $\Gamma$ are dense in $\relint(\sigma)$ and hence we may assume that $\ub_0$ is such a point, i.e. there is $\zb \in U_\Delta$ with coordinates in $\kdop$ and $\ub_0=\val(\zb)$. By  the coordinate transform $\xb':=\xb/{\mathbf z}$, we may assume that $\ub_0=\mathbf 0$.

Note that  ${\rm LC}_{\mathbf 0}(\val(X))$ is a finite union of $\qdop$-rational polyhedral cones centered at $\mathbf 0$. By convex geometry, the convex hull of ${\rm LC}_{\mathbf 0}(\val(X))$  is a finite intersection of half spaces  $\{\ub \cdot \mb \geq 0\}$ with $\mb \in \zdop^n$. To show concavity in $\mathbf 0$, we have to prove that the convex hull is a linear subspace. If ${\rm LC}_{\mathbf 0}(\val(X))$ is contained in a half space $\{\ub \cdot \mb \geq 0\}$  as above, then it is enough to show that ${\rm LC}_{\mathbf 0}(\val(X)) \subset \{\ub \cdot \mb=0\}$. 

By shrinking $\Delta$, we may assume that $\val(X)$ is contained in $\{\ub \cdot \mb \geq 0\}$ and that for $j=1, \dots, r$, there is $v_j \in V_j$ with $\val(v_j)=\mathbf 0$. Then $|\xb^\mb|$ takes its maximum $1$ in $v_1, \dots, v_r$. But every point in $V_j$ has the same reduction in $\tilde{X}$ as $v_j$ and hence $|\xb^\mb|=1$ on $X^\circ$. We conclude that $\val(X^\circ) \subset \{ \ub \cdot \mb =0 \}$. This proves ${\rm LC}_{\mathbf 0}(\val(X)) \subset \{ \ub \cdot \mb = 0\}$. \qed

\begin{prop} \label{pure dimension}
Let $X$ be a closed analytic subvariety of $U_\Delta$ such that $\val(X)\cap \Int(\Delta) \neq \emptyset$. If $X$ is of pure dimension $d$, then $\val(X) \cap {\Int(\Delta)}$ is also of pure dimension $d$.
\end{prop}

\proof We have seen in Proposition \ref{polyhedron} that $\val(X)$ is a $\Gamma$-rational polytopal set. Moreover, its proof or  \cite{Ber5}, Corollary 6.2.2, show that $\dim(\val(X)) \leq d$ holds even without considering interior points of $\Delta$. By subdivision of $\Delta$, it is enough to prove that $\val(X)$ is at least $d$-dimensional.

We proceed by induction on $N:= \dim(\val(X))$. We may assume that $X$ is irreducible. Then $\val(X)$ is connected (see \ref{tropical analytical variety}). If $N=0$, then $\val(X)$ is an interior point of $\Delta$. As in the proof of Proposition \ref{totally concave}, we conclude that $\tilde{X}$ is finite and hence $\dim(X)=0$.

Now assume $N>0$. 
By shrinking $\Delta$, we may assume that $\val(X)$ is of pure dimension. By density of $X(\kdop)$, there is $\ub \in {\Int(\Delta)} \cap \val(X(\kdop))$. By a change of coordinates as in the proof of Proposition \ref{totally concave}, we may assume that $\mathbf 1 \in X$ and $\ub = \mathbf 0$. There is $\mb \in \zdop^n \setminus \{\mathbf 0\}$ such that the hyperplane $\{\ub \cdot \mb =  0\}$ intersects $\val(X)$ transversally. The dimension of the closed analytic subvariety $X' := X \cap \{\xb^\mb = 1\}$ is $d-1$.  We have 
$$\val(X') \subset \val(X) \cap \{ \ub \cdot \mb =0\}$$
and hence $\dim(\val(X')) \leq N-1$.
By induction, we get $\dim(\val(X')) \geq d-1$.   We conclude $d \leq N$ proving the claim. \qed

\begin{art} \rm  \label{globalization2}
As in \ref{globalization}, let $\Ucal$ be the admissible formal scheme over $\kdop^\circ$ associated to the $\Gamma$-rational polytopal complex $\Ccal$ and let $U:=\Ucal^{\rm an}$. For a closed analytic subvariety $X$ of $U$, the set $\val(X)$ is called the {\it tropical variety associated to $X$}. We set $\Pi := \cup_{\Delta \in \Ccal} \Delta$. 
\end{art}

\begin{thm} \label{tropical main theorem} Under the assumptions in \ref{globalization2}, the following properties hold:
\begin{itemize}
\item[(a)] $\val(X)$ is a locally finite union of $\Gamma$-rational polytopes in $\rdop^n$.
\item[(b)] $\val(X) \cap \Int(\Pi)$ is totally concave.
\item[(c)] If $\val(X) \cap \Int(\Pi)$ is non-empty and if $X$ is of pure dimension $d$, then $\val(X) \cap \Int(\Pi)$ is of pure dimension $d$.
\end{itemize} 
\end{thm}

\proof Statement (a) is immediate from Proposition \ref{polyhedron}. Let $\ub \in \Int(\Pi)$. If $\ub \in \Int(\Delta)$ for some $\Delta \in \Ccal$, then (b) follows from Proposition \ref{totally concave}. If no such $\Delta$ is available, then one may easily adjust the polytopes in $\Ccal$ without changing their union $\Pi$ such that $\ub \in \Int(\Delta)$ for some $\Delta \in \Ccal$. Similarly, we deduce (c) from Proposition \ref{pure dimension}. \qed

\begin{rem} \label{algebraic situation} \rm 
This is most useful if $\Ccal$ is a polytopal decomposition of $\rdop^n$. Then no boundary points occur and (b), (c) hold for $\val(X)$. In particular, Theorem \ref{tropical main theorem} holds for a closed algebraic subvariety $X$ of $\Tor$ over $\kdop$ and hence implies the well known statements from tropical algebraic geometry (see \cite{EKL}, \S 2). The only thing which does not hold analytically is that $\val(X)$ is a finite union of polyhedrons.
\end{rem}

Now we are able to deduce a toric version of Theorem \ref{Theorem 2}.

\begin{thm} \label{toric dimension bound}
Let $\Xcal'$ be a quasicompact strictly semistable formal scheme over $\kdop^\circ$ with generic fibre $X'$ and let $X$ be a $d$-dimensional closed analytic subvariety of $(\Tor)_\kdop^{\rm an}$. If there is a dominant morphism $f:X' \rightarrow X$, then $\tilde{\Xcal}'$ has a point contained in at least $d+1$ irreducible components.
\end{thm}

\proof By Theorem \ref{tropical main theorem}, there is a $d$-dimensional $\Gamma$-rational polytope $\Delta$ in $\val(X)$. The quasicompact set $f^{-1}(U_\Delta)$ may be covered by finitely many sets $U':=(\Ucal')^{\rm an}$ of the same form as in the proof of Proposition \ref{polyhedron}. The same proof shows that $\Delta$ is a finite union of simplices $f_{\rm aff}^{(0)}(\Sigma(r,\pi))$. It follows that $r \geq d$ for at least one $\Ucal'$ proving the claim. \qed

\begin{art} \rm \label{closure}
In the remaining part of this section,  we consider  a $\Gamma$-rational polytopal decomposition $\Ccal$ of $\rdop^n$. The associated admissible formal scheme over  $\kdop^\circ$ is denoted by $\Ucal$. Let $X$ be a closed analytic subvariety of $(\Tor)_\kdop^{\rm an} = \Ucal^{\rm an}$ of pure dimension $d$. In the following, we relate $\Ccal$ to the {\it closure} $\Xcal$ of $X$ in $\Ucal$. The closure is the $\kdop^\circ$-model $\Xcal$ of $X$ locally defined by
$$\Ucal_\Delta \cap \Xcal := \Spf \left( \kdop \langle U_\Delta \rangle^\circ / \left(I_\Delta(X) \cap \kdop \langle U_\Delta \rangle^\circ \right)\right),$$
where $I_\Delta(X)$ is the ideal of vanishing on $U_\Delta$ (see \cite{Gu2}, Proposition 3.3). Note that $\tilde{\Xcal}$ is a closed subvariety of $\tilde{\Ucal}$ of pure dimension $d$.  
\end{art}

\begin{lem} \label{transverse lemma}
Let $\Delta \in \Ccal$ with $\codim(\Delta,\rdop^n)=d$. We assume that $\Delta \cap \val(X)$ is a non-empty finite subset of $\tau := \relint(\Delta)$. Then the toric variety $Y_\Delta$ in $\tilde{\Ucal}_\Delta$ (see Remark \ref{polytopal toric variety}) is an irreducible component of $\tilde{\Xcal}$.
\end{lem} 

\proof There is a $\Gamma$-rational half-space $H_+:=\{\mb \cdot \ub \geq c\}$ containing the finite set $\Delta \cap \val(X)$ such that the boundary hyperplane intersects $\Delta \cap \val(X)$ in a single point $\ub$. Then $\xb^\mb$ achieves its maximum absolute value $e^{-c}$  on $X \cap U_\Delta$ in every $x \in X \cap U_\Delta$ with $\val(x) = \ub$. The  Shilov boundary of the affinoid space $X \cap U_\Delta$ is given by the points reducing to the generic points of
$$(X \cap U_\Delta)\sptilde = \left( ( \Xcal \cap \Ucal_\Delta)^{\rm f-an} \right)\sptilde$$  
(see \cite{Ber}, Proposition 2.4.4). Hence there is an irreducible component $Z$ of $(\Xcal^{\rm f-an})\sptilde$ with $\val(\xi_Z) = \ub$. Let $Y$ be the image of $Z$ under the canonical finite surjective morphism $\tilde{\iota}: (\Xcal^{\rm f-an})\sptilde  \rightarrow \tilde{\Xcal}$. Note that $Y$ is $d$-dimensional and has the generic point $\tilde{\iota}(\tilde{\xi}_Z)$. By Proposition \ref{torus orbits 1},  $\tilde{\iota}(\tilde{\xi}_Z)$ is also contained in the $d$-dimensional torus orbit $Z_\tau$ and hence $Y= \overline{Z}_\tau$. \qed

\begin{thm} \label{transversal correspondence}
Under the hypothesis of \ref{closure}, we assume that $\Ccal$ is transversal to $\val(X)$ (see \ref{local cone and transversality}). Then there is a bijective correspondence between:
\begin{itemize}
\item[(a)] equivalence classes of transversal vertices of $\Ccal \cap \val(X)$ (see \ref{local cone and transversality}), 
\item[(b)] irreducible components $Y$ of $\tilde{\Xcal}$.
\end{itemize}
An equivalence class in (a) is contained in a unique $\Delta \in \Ccal$ of codimension $d$. The corresponding irreducible component $Y$ is the toric variety $Y_\Delta$ in $\tilde{\Ucal}$ from Remark \ref{polytopal toric variety}.
\end{thm}

\proof We have seen in Lemma \ref{transverse lemma} that $Y_\Delta$ is an irreducible component of $\tilde{\Xcal}$. Conversely, let $Y$ be an irreducible component of $\tilde{\Xcal}$. Using the notation from the proof of Lemma \ref{transverse lemma}, there is an irreducible component $Z$ of $(\Xcal^{\rm f-an})\sptilde $ with $\tilde{\iota}(Z)=Y$.

We claim that $\ub_Z:= \val(\xi_Z)$ is a transversal vertex of $\Ccal \cap \val(X)$. To see this, let $\Delta$ be the unique polytope from $\Ccal$ with $\ub_Z \in \tau := \relint(\Delta)$. By definition of a transversal vertex, we have to prove $\codim(\Delta,\rdop^n)=d$. 
By Proposition \ref{torus orbits 1} and Proposition \ref{torus orbits global}, we have $\codim(\Delta,\rdop^n) = \dim(Z_\tau)$ and $Z_\tau$ contains the generic point $\tilde{\iota}(\tilde{\xi}_Z)$ of $Y$, hence
$$\codim(\Delta,\rdop^n) \geq \dim(Y) =d.$$
Since $\val(X)$ is of pure dimension $d$ and $\ub_Z \in \val(X) \cap \tau$, transversality yields $\codim(\Delta, \rdop^n)=d$ proving that $\ub_Z$ is a transversal vertex of  $\Ccal \cap \val(X)$. Moreover, we see that $Y= \overline{Z}_\tau= Y_\Delta$. This shows that the map $Y \mapsto \Delta$ is independent of the choice of $Z$. By construction, it is inverse to the map $\Delta \mapsto Y_\Delta$ from the beginning. \qed

\section{Mumford's construction}

We review Mumford's construction of models $\Acal$ of a totally degenerate abelian variety $A$. For a closed subscheme $X$ of $A$, we  study the periodic tropical variety $\val(X)$ using the previous section. If we choose the polytopal decomposition for $\Acal$ transversally to $\val(X)$, then the irreducible components of the closure of $X$ in $\Acal$ turn out to be toric varieties.


In this section, $A$ denotes a {\it totally degenerate abelian variety} over $\kdop$, i.e. 
$A^{\rm an}$ is isomorphic to $(\Tor)^{\rm an}_\kdop / M$, where $M$ is a subgroup of $\Tor(\kdop)$ which maps isomorphically onto a complete lattice $\Lambda$ of $\rdop^n$ under the map $\val$. Such an $M$ is called a {\it lattice} of $(\Tor)^{\rm an}_\kdop$.
 
Let $\rdop^n \rightarrow \rdop^n/\Lambda, \, \ub \mapsto \overline \ub$, be the quotient map. Clearly, the map $\val$ from Section 4 descends to a continuous map $\overline{\val}: A^{\rm an} \rightarrow \rdop^n / \Lambda$.
First, we translate the notions of convex geometry introduced  in the appendix to the torus $\rdop^n/\Lambda$. 

\begin{art} \label{toric convex geometry} \rm 
A {\it polytope} $\overline{\Delta}$ in $\rtor$ is given by a polytope $\Delta$ in $\rdop^n$ such that $\Delta$ maps bijectively onto $\Deltabar$. We say that $\overline{\Delta}$ is $\Gamma$-rational if $\Delta$ is $\Gamma$-rational. A ($\Gamma$-rational) {\it polytopal set} $S$ in $\rtor$ is a finite union of ($\Gamma$-rational) polytopes in $\rtor$.  

A {\it polytopal decomposition} of $\rdop^n/\Lambda$ is a finite family $\overline{\Ccal}$ of polytopes in $\rdop^n/\Lambda$ induced by a $\Lambda$-periodic polytopal decomposition $\Ccal$ of $\rdop^n$. It is easy to see that $\rdop^n/\Lambda$ has a $\Gamma$-rational polytopal decomposition. The other notions from  the appendix transfer also to the periodic situation. 
\end{art}



\begin{art} \rm \label{Mumford's construction}
Let $\overline{\Ccal}$ be a $\Gamma$-rational polytopal decomposition of $\rdop^n/\Lambda$. By \ref{globalization}, $(U_\Delta)_{\Delta \in \Ccal}$ is a formal analytic atlas of $(\Tor)_\kdop^{\rm an}$. We may form the quotient by $M$ leading to a formal analytic variety over $\kdop$. The associated formal scheme $\Acal$ is a $\kdop^\circ$-model of 
$A= (\Tor)^{\rm an}_\kdop/M$ which has  a covering by formal open affine sets $\Ucal_{\Deltabar}$ obtained by gluing 
$\Ucal_{\Delta + \mathbf{\lambda}}$ for all $\mathbf{\lambda} \in \Lambda$.

The generic fibre of the formal torus  ${\mathbb T}_1$
acts naturally on $A$ and there is a unique extension to an action of ${\mathbb T}_1$ on $\Acal$. We get a torus action of $\tilde{{\mathbb T}}_1 = (\Tor)_{\tilde{\kdop}}$ on the special fibre $\tilde{\Acal}$. On $\tilde{\Ucal}_\Deltabar$, this action agrees with the action on $\Ucal_\Delta$ defined in \ref{torus action}. 

Using the $\Lambda$-periodic decomposition $\Ccal$ and passing to the quotient, we may transfer the results from \S 4 and \S 5 to $A$. By  Proposition \ref{torus orbits 1} and Proposition \ref{torus orbits global}, we get:
\end{art}

\begin{prop} \label{torus orbit 2}
For the formal $\kdop^\circ$-model $\Acal$ of $A$ associated to $\overline{\Ccal}$, we have:
\begin{itemize}
\item[(a)] There is a bijective order reversing correspondence between torus orbits $Z$ of $\tilde{\Acal}$ and open faces $\overline \tau$ of $\Ccalbar$, given by
$$\overline \tau = \overline{\val}(\pi^{-1}(Z)), \quad Z= \pi(\overline{\val}^{-1}(\overline{\tau})),$$
where $\pi: A \rightarrow \tilde{\Acal}$ is the reduction map. Moreover, we have $\dim(Z)+\dim(\tau)=n$.
\item[(b)]
The irreducible components $Y$ of $\tilde{\Acal}$ are toric varieties and correspond to the vertices $\ub$ of $\Ccalbar$ by $\overline{\ub}:= \valbar(\xi_Y)$.
\end{itemize}
\end{prop}

\begin{prop} \label{multiplication with m}
Let $\Ccalbar$ be a $\Gamma$-rational polytopal decomposition of $\rdop^n/\Lambda$ and let $m \in \zdop \setminus \{0\}$. Then 
$\overline{\frac{1}{m} \Ccal} := \left\{ \overline{\frac{1}{m} \Delta} \mid \Delta \in \Ccal \right \}$
is also a $\Gamma$-rational polytopal decomposition of $\rtor$. The associated $\kdop^\circ$-model $\Acal_m$ of $A$ has the following properties:
\begin{itemize}
\item[(a)] The morphism $[m]: A \rightarrow A, \, x \mapsto mx$, has a unique extension to a morphism $\varphi_m: \Acal_m \rightarrow \Acal_1$ of admissible formal schemes over $\kdop^\circ$.
\item[(b)] The morphism $\tilde{\varphi}_m$ is finite of degree $m^{2n}$.
\item[(c)] The behaviour of the reduction $\tilde{\varphi}_m$ with respect to the torus actions is given by
$$\tilde{\varphi}_m(\tb \cdot \zb) = \tb^m \cdot \tilde{\varphi}_m(\zb) \quad (\zb \in \tilde{\Acal}_m, \, \tb \in (\Tor)_{\tilde{\kdop}}).$$
\item[(d)] The inverse image of a $k$-dimensional torus orbit of $\tilde{\Acal}_1$ with respect to $\tilde{\varphi}_m$ is equal to the disjoint union of $m^n$ $k$-dimensional torus orbits of $\tilde{\Acal}_m$.
\end{itemize}
\end{prop}

\proof Obviously, $\overline{\frac{1}{m}\Ccal}$ is a $\Gamma$-rational polytopal decomposition. The extension of $[m]$ is constructed locally by $\Ucal_{\frac{1}{m}\Delta} \rightarrow \Ucal_\Delta, \, \xb \mapsto \xb^m$. Uniqueness is clear formal analytically and hence follows from \ref{admissible formal schemes}. This proves (a). By construction, we get immediately (b) and (c). Now (c) implies that the inverse image of a $k$-dimensional torus orbit $O$ is the disjoint union of $k$-dimensional torus orbits. By Proposition \ref{torus orbit 2}, $O$ corresponds to an  open face $\overline{\tau}$ of $\Ccalbar$ of dimension $n-k$. Since $\{ \overline{\ub} \in \rdop^n/\Lambda \mid \overline{m \ub} \in \overline{\tau}\}$ is the disjoint union of $m^n$ open faces, Proposition \ref{torus orbit 2} yields (d). \qed

\begin{art} \rm \label{toric line bundles}
We describe line bundles on $A = (\Tor)^{\rm an}_\kdop/M$ similarly as in the complex analytic situation (see \cite{FvdP}, Ch. VI, and \cite{BL2}, \S 2, for details). 

Let $L$ be a line bundle on $A$. The pull-back to $T:=(\Tor)^{\rm an}_\kdop$ with respect to the quotient morphism $p$ is trivial and will be identified with $T \times \kdop$. It is given by a cocycle $\gamma \mapsto Z_\gamma$ of $H^1(M, \Ocal(T)^\times)$ and $L= (T \times \kdop)/M$ where the quotient is with respect to the $M$-action
$$M \times \left( T \times \kdop \right)  \longrightarrow 
T \times \kdop, \quad (\gamma, (\xb, \alpha)) \mapsto (\gamma \cdot \xb, Z_\gamma(\xb)^{-1} \alpha).$$
The cocycle has the form $Z_\gamma(\xb)=d_\gamma \cdot \sigma_\gamma(\xb)$, where $\gamma \mapsto \sigma_\gamma$ is a homomorphism of $M$ to the character group $\check{T}$ and where $d_\gamma \in \kdop^\times$ satisfies
\begin{equation} \label{quadratic}
d_{\gamma \rho} \cdot d_\gamma^{-1} \cdot d_\rho^{-1} = \sigma_\rho(\gamma) \quad (\gamma, \rho \in M). 
\end{equation}
By the isomorphism $M \stackrel{\val}{\rightarrow} \Lambda$, we get a unique symmetric bilinear form $b$ on $\Lambda$ characterized by
$$b(\val(\gamma),\val(\rho)) = v\left( \sigma_\rho(\gamma) \right).$$
Then $b$ is positive definite on $\Lambda$ if and only if $L$ is ample. 
Note that the cocycle $Z_\gamma$ factors over $\rdop^n$, i.e. for every $\lambda = \val(\gamma) \in \Lambda$, there is a unique real function $z_\lambda$ on $\rdop^n$ such that 
$$z_{\lambda}(\val(\xb))= v (Z_\gamma(\xb)) \quad (\gamma \in M, \, \xb \in T).$$
The function $z_\lambda$ is affine with
\begin{equation} \label{cycle and bilinear}
z_\lambda(\ub)= z_\lambda(\mathbf 0) + b(\ub, \lambda) \quad (\lambda \in \Lambda, \, \ub  \in \rdop^n).
\end{equation}
\end{art}

\begin{prop} \label{Mumford's line bundles}
Let $L$ be a line bundle on $A$. Repeat that $L$ is given analytically by $L= (T \times \kdop)/M$ and by a cocycle $(Z_\gamma)_{\gamma \in M}$ leading to a family $(z_\lambda)_{\lambda   \in \Lambda}$ of real functions as above. Let $\Acal$ be the formal $\kdop^\circ$-model of $A$ associated to a given $\Gamma$-rational polytopal decomposition $\Ccalbar$ of $\rtor$. Then there is a bijective correspondence between isomorphism classes of formal $\kdop^\circ$-models $\Lcal$ of $L$ on $\Acal$ with trivialization $(\Ucal_\Deltabar)_{\Deltabar \in \Ccalbar}$ and continuous real functions $f$ on $\rdop^n$ satisfying the following two conditions:
\begin{itemize}
\item[(a)] For every $\Delta \in \Ccal$, there are $\mb_\Delta \in \zdop^n$ and $c_\Delta \in \Gamma$ with $f(\ub)= \mb_\Delta \cdot \ub + c_\Delta$ on $\Delta$.
\item[(b)] $f(\ub + \lambda) = f(\ub) + z_\lambda(\ub) \quad (\lambda \in \Lambda, \ub \in \rdop^n)$.
\end{itemize}
If $\metr_\Lcal$ denotes the formal metric on $L$ associated to $\Lcal$ (see \ref{formal metrics}), then the correspondence is given by $f_\Lcal \circ \val := - \log \circ p^* \|1\|_\Lcal$ on $T$.
\end{prop}

{\bf First step of proof: \/} {\it $g=\sum_{\nu \in \zdop^n} a_\nu \xb^\nu \in \kdop \langle \val^{-1}(\Delta) \rangle$ is a unit if and only if there is a $\nu_0 \in \zdop^n$ such that $|a_{\nu_0} \xb^{\nu_0}| > |a_\nu \xb^\nu|$ for all $\xb \in U_\Delta$ and all $\nu \neq \nu_0$.}

\vspace{2mm}
If $g$ has such a dominant term, then we may assume that $\nu_0 = \mathbf 0$ and $a_{\mathbf 0}=1$. Then we have $g:=1-h$ with $|h|_{\rm sup} < 1$ and $g^{-1}= \sum_{n=0}^\infty h^n \in \kdop \langle \val^{-1}(\Delta) \rangle$. Conversely, if $g$ has no dominant term, then there is $|a_{\nu_0} \xb^{\nu_0}|= |a_{\nu_1} \xb^{\nu_1}| = |g|_{\rm sup}$ for certain $\nu_0 \neq \nu_1$ and $\xb \in U_\Delta$. Note that $\val^{-1}(\ub)$ is isomorphic to the affinoid torus ${\mathbb T}_1^{\rm an}$ for $\ub := \val(\xb)$. Then the restriction of $g$  to ${\mathbb T}_1^{\rm an}$ has no dominant term as well and hence it is not invertible on ${\mathbb T}_1^{\rm an}$ (\cite{BGR}, Lemma 9.7.1/1). Since $\val^{-1}(\ub)$ is an analytic subdomain of $U_\Delta$, $g$ is not an unit of $\kdop\langle \val^{-1}(\Delta) \rangle$. 

{\bf Second step: \/} {\it $f_\Lcal$ is a continuous function satisfying (a) and (b).}

\vspace{2mm}
Continuity follows from the continuity of formal metrics on analytic spaces and (b) is by construction (see \ref{toric line bundles}). On a trivialization $U_{\Deltabar}$ of $\Lcal$, the section $1$ corresponds to a unit $g$ in $\kdop \langle U_\Delta \rangle$. By definition of formal metrics, we have $p^* \|1\|_\Lcal=|g|$ on $U_\Delta$. 
By the first step, $g$ has a dominant term $a_\Delta \xb^{\mb_\Delta}$ leading to (a) with $c_\Delta:=v(a_\Delta)$. This proves the second step.

\vspace{2mm}
Now the proof is quite easy. A continuous function $f$ on $\rdop^n$ gives rise to a metric $\metr'$ on $p^*(L)=T \times \kdop$ by $\|1\|'=e^{-f\circ \val}$. If $f$ satisfies (b), then $\metr'$ passes to the quotient modulo $M$, i.e. there is a unique metric $\metr_f$ on $L$ with 
$$f\circ \val = - \log \circ p^* \| 1 \|_f. $$
Now we assume that $f$ also satisfies (a). Since $\Gamma$ is the value group, there is $a_\Delta \in \kdop^\times$ with  $c_\Delta = v(a_\Delta)$. 
For every $\Delta \in \Ccal$, the unit $(a_\Delta \xb^{\mb_\Delta})^{-1}$ gives a frame of $p^*(L)=T \times \kdop$ over $U_\Delta$ leading to a trivialization of $L$ over $U_\Deltabar$. This trivialization extends to the trivial line bundle over $\Ucal_\Deltabar$ and induces a metric $\metr_\Delta$ on $L|_{U_\Deltabar}$. A priori, $\metr_\Delta$ depends on the choice of $\Delta$, but by construction and (a), we deduce that $\metr_\Delta$ agrees with $\metr_f$ over $U_{\Deltabar}$. Therefore the transition functions $g_{\Delta\Delta'}$ of the above trivializations have constant absolute value $1$ on $U_\Deltabar \cap U_{\Deltabar'}$. This means that they define a formal $\kdop^\circ$-model $\Lcal_f$ of $L$ on $\Acal$ with associated metric equal to $\metr_f$. 

We have to prove that  $f \mapsto \Lcal_f$ is inverse to $\Lcal \mapsto f_\Lcal$. For $\Mcal := \Lcal_f$, we have
$$f_\Mcal \circ \val = -\log \circ p^*\|1\|_\Mcal= -\log\circ p^*\|1\|_f = f \circ \val$$
and hence $f_\Mcal=f$. Conversely, it is clear that $\metr_\Lcal=\metr_{\Lcal_g}$ for given $\Lcal$ and $g:=f_\Lcal$. By considering trivializations as above or by using \cite{Gu2}, Proposition 5.5, we see that the formal $\kdop^\circ$-models $\Lcal$ and $\Lcal_g$ of $L$ are isomorphic.   \qed

\begin{cor} \label{ample reduction}
In the notation of Proposition \ref{Mumford's line bundles}, let $\tilde{\Lcal}$ be the reduction of $\Lcal$ on $\tilde{\Acal}$. Then $\tilde{\Lcal}$ is ample if and only if $f_\Lcal$ is a strongly polyhedral convex function with respect to $\Ccal$.
\end{cor}


\proof Note that $\tilde{\Lcal}$ is ample if and only if its restriction to every irreducible component $Y$ of $\tilde{\Acal}$ is ample. By Proposition \ref{torus orbit 2}, $Y$  is a toric variety corresponding to a vertex $\ubb$ of $\Ccalbar$. For simplicity of notation, we may assume that $\ub = \mathbf 0$. For $\Delta \in \Star(\ub)$ and $\tau := \relint(\Delta)$, we choose the equation $\tilde{\xb}^{\mb_\Delta}$ on the torus orbit $Z_\tau$. By Proposition \ref{Mumford's line bundles}(a), it is easy to see that we get a Cartier divisor $D$ on $Y$ of a meromorphic section of $\tilde{\Lcal}|_Y$. Let $\psi_D$ be the continuous function on the complete fan of $Y$ centered  in $\ub=\mathbf 0$ which is equal to $-\mb_\Delta \cdot \ub$ on $\rdop \Delta$. By \cite{Fu2}, 3.4, $-\psi_D$ is a strongly polyhedral convex function. Note that $f_\Lcal=-\psi_D$ on every $\Delta \in \Star(\ub)$. We deduce easily that $f_\Lcal$ is a strongly polyhedral convex function with respect to $\Ccal$. \qed

\begin{art} \label{periodic tropic assumption} \rm 
In the remaining part of this section, we consider a closed subscheme $X$ of $A$ of pure dimension $d$. The {\it tropical variety} $\valbar(X^{\rm an})$ is well-defined in $\rtor$. The following properties are easily deduced from Theorem \ref{tropical main theorem} and the continuity of $\valbar$.
\end{art}

\begin{thm} \label{periodic main theorem}
The tropical variety $\valbar(X^{\rm an})$ is a totally concave $\Gamma$-rational polytopal set in $\rtor$ of pure dimension $d$. If $X$ is connected, then $\valbar(X^{\rm an})$ is also connected. 
\end{thm}

Let $\Acal$ be the formal $\kdop^\circ$-model associated to the $\Gamma$-rational polytopal decomposition $\Ccalbar$ of $\rtor$. We denote by $\Xcal$ the closure of $X^{\rm an}$ in $\Acal$ (see \cite{Gu2}, Proposition 3.3). From Theorem \ref{transversal correspondence}, we deduce immediately:

\begin{thm} \label{periodic transverse} 
Unter the hypothesis above and assuming that $\Ccalbar$ is $\valbar(X^{\rm an})$-transversal, we have a bijective correspondence between irreducible components of $\tilde{\Xcal}$ and equivalence classes of transversal vertices of $\Ccalbar \cap \valbar(X^{\rm an})$. If $\ubb$ is a transversal vertex, then there is a unique  $\Deltabar \in \Ccalbar$ with $\ubb \in \Deltabar$ and $\codim(\Delta,\rdop^n)=d$.  The corresponding irreducible component of $\tilde{\Xcal}$ is the closure $Y_{\Deltabar}$ of the torus orbit of $\tilde{\Acal}$ associated to the open face $\relint(\Deltabar)$ (see Proposition \ref{torus orbit 2}). In particular, every irreducible component of $\tilde{\Xcal}$ is a toric variety.
\end{thm}

\section{Semistable alterations}

As in the previous section, we consider a totally degenerate abelian variety $A^{\rm an} =(\Tor)_\kdop^{\rm an}/M$ over $\kdop$ with Mumford model $\Acal$ over $\kdop^\circ$ associated to the $\Gamma$-rational polytopal decomposition $\Ccalbar$ of $\rtor$, where $\Lambda := \val(M)$. We fix also an irreducible  closed subvariety $X$ of $A$ with closure $\Xcal$ in $\Acal$. The goal of this section is to describe the multiplicities of an irreducible component $Y$ of $\tilde{\Xcal}$ using a strictly semistable alteration. Under a non-degeneracy condition for $Y$, this may be done in terms of convex geometry and will be used in the following section to prove the main result.

\begin{art} \rm \label{semistable situation}
Let $\Xcal'$ be a strictly semistable formal scheme over $\kdop^\circ$ with generic fibre $X'$ and let $f:X' \rightarrow A^{\rm an}$ be a morphism over $\kdop$. In the first paragraphes, we will show that the polytopal decomposition $\Ccalbar$ endows $X'$ with a canonical formal analytic structure $\frak X'$ and with a morphism $\varphi: {\frak X}' \rightarrow \Acal^{\rm f-an}$ of formal analytic varieties.

By Proposition \ref{semistable properties}, $\Xcal'$ is covered by formal open affine subsets $\Ucal'$ such that all irreducible components of $\tilde{\Ucal}'$ pass through  $\tilde{P} \in \tilde{\Ucal}'(\tilde{\kdop})$ and with an \'etale morphism
$$\psi: {\Ucal'} \longrightarrow {\Scal} :=\Spf \left( \kdop^\circ \langle x_0', \dots, x_d' \rangle / \langle x_0' \dots x_r' = \pi \rangle \right)$$
for suitable $\pi \in \kdop^{\circ \circ}$ such that $\tilde{\psi}^{-1}(\tilde{\mathbf 0}) = \{\tilde{P}\}$. The simplex $$\Delta(r,\pi):=\{\ub' \in \rdop_+^{r+1} \mid u_0'+ \dots + u_r' = v(\pi) \}$$ is canonically associated to $\Scal$. To apply \S 5, we represent $\Delta(r,\pi)$ by the standard simplex $$\Sigma(r,\pi)=\{\ub' \in \rdop_+^r \mid u_1'+ \dots u_r' \leq v(\pi)\}$$ omitting $u_0'$. Then we have
$$\Scal = \Spf \left( \kdop^\circ \langle x_0', \dots , x_r' \rangle / \langle x_0' \cdots x_r'= \pi \rangle \right) \times \Spf \left( \kdop^\circ \langle x_{r+1}', \dots , x_d' \rangle \right).$$
We denote the first factor by $\Scal_1$. By definition, the morphism $\phi: \Ucal' \rightarrow \Scal_1$, obtained by composition with the first projection, is smooth and we have $\tilde{\phi}(\tilde{P}) = \tilde{\mathbf 0}$. We have an isomorphism $\Scal_1 \stackrel{\sim}{\rightarrow} \Ucal_{\Sigma(r,\pi)}$ by omitting $x_0'$. By composition, we get a morphism $\phi_0:\Ucal' \rightarrow \Ucal_{\Sigma(r,\pi)}$ Note that we have maps $\val$ on $\Scal_1^{\rm an}$ and $\Ucal_{\Sigma(r,\pi)}^{\rm an}$ with images $\Delta(r,\pi)$ and $\Sigma(r,\pi)$. 
\end{art}

\begin{art} \rm \label{lift} 
The generic fibre of $\Ucal'$ is denoted by $U'$.
We claim that $f:U' \rightarrow A^{\rm an}$ has a lift $F:U' \rightarrow (\Tor)_\kdop^{\rm an}$, unique up to $M$-translation. For a Hausdorff analytic space $Y$, we consider the cohomology group $H^1(Y,\zdop)$ of the underlying topological space. Note that it is the same as 
$H^1(Y_{\rm rig}, \zdop)$ for the underlying rigid space $Y_{\rm rig}$ (see \cite{Ber2}, 1.6). By \cite{Ber4}, Theorem 5.2, the generic fibre of a strictly semistable formal scheme over $\kdop^\circ$ is contractible to the skeleton. The skeleton of $\Ucal'$ is $\Delta(r,\pi)$ and therefore $H^1(U',\zdop)=0$. Since $(\Tor)_\kdop^{\rm an}$ is the universal covering space of $A^{\rm an}$, we get a lift $F$ as desired (see \cite{BL2}, Theorem 1.2). 
\end{art} 

\begin{art} \rm \label{affine map}
By the proof of Proposition \ref{polyhedron}, there is a unique map $F_{\rm aff}: \Delta(r,\pi) \rightarrow \rdop^n$ with
$$F_{\rm aff} \circ \val \circ \phi = \val \circ F$$
on $U'$. Let $f_{\rm aff}$ be such an $F_{\rm aff}$ (without fixing $F$), it is determined up to $\Lambda$-translation. Then $\overline{f}_{\rm aff}: \Delta(r,\pi) \rightarrow \rdop^n/\Lambda$ is uniquely characterized by 
\begin{equation} \label{affine characterization}
\overline{f}_{\rm aff} \circ \val \circ \phi = \overline{\val} \circ f
\end{equation}
on $U'$ Note that uniqueness always follows from $\val\circ \phi (U')=\Delta(r,\pi)$. This was also part of the proof of Proposition \ref{polyhedron}, where we have considered the affine map $f_{\rm aff}^{(0)}:\Sigma(r,\pi) \rightarrow \rdop^n$ determined by 
$$f_{\rm aff}^{(0)}(u_1', \dots, u_r') = f_{\rm aff}(u_0', \dots, u_r')$$
for $\ub' \in \Delta(r,\pi)$. By \eqref{semistable f-coord}, there are $\mb_i \in \zdop^r$ and $\lambda_i \in \kdop^\times$ $(i=1, \dots, n)$ such that
\begin{equation} \label{explicit affinization}
f_{\rm aff}^{(0)}(\ub') = \left( \mb_i \cdot \ub' + v(\lambda_i) \right)_{i=1, \dots, n}
\end{equation}
for every $\ub' \in \Sigma(r,\pi)$.
\end{art}

\begin{art} \label{formal refinement} \rm 
Now we are ready to describe the formal analytic structure $\Xfrak'$ on $X'$ induced by $\Ccalbar$. It is given by the atlas $U' \cap f^{-1}(U_{\Deltabar})$, where $U'$ ranges over all formal open affinoids as in \ref{semistable situation} and where $\Deltabar \in \Ccalbar$. Note that $U' \cap f^{-1}(U_{\Deltabar}) = U' \cap F^{-1}(U_\Delta)$ is a Weierstrass domain in $U'$. We have  unique morphisms $\iota: \Xfrak' \rightarrow (\Xcal')^{\rm f-an}$ and $\varphi:\Xfrak' \rightarrow \Acal^{\rm f-an}$ extending the identity on $X'$ and $f:X' \rightarrow A$, respectively.
\end{art}

Our next goal is to relate this formal analytic structure on a subset $U'$ from \ref{semistable situation} to the simplex $\Sigma(r,\pi)$. Let $\Ufrak'$ be the formal analytic variety on $U'$ induced by $\Xfrak'$. Note that
$$\left(f_{\rm aff}^{(0)}\right) ^{-1} (\Ccal):= \left(\left(f_{\rm aff}^{(0)} \right)^{-1} (\Delta) \right)_{\Delta \in \Ccal}$$
is a $\Gamma$-rational polytopal decomposition of $\Sigma(r,\pi)$. We denote the associated formal scheme (see \ref{globalization}) by $\Tcal$ coming with a canonical morphism $i:\Tcal \rightarrow \Ucal_{\Sigma(r,\pi)}$ extending the identity.

\begin{prop} \label{convex formal analytic}
Under the hypothesis above, the following properties hold:
\begin{itemize}
\item[(a)] $\Ufrak'$ is given by the atlas $\phi^{-1}_0(U_\sigma)$, $\sigma \in \left(f_{\rm aff}^{(0)}\right)^{-1}(\Ccal)$. 
\item[(b)] There is a unique morphism $\phi_0': \Ufrak' \rightarrow \Tcal^{\rm f-an}$ with $i^{\rm f-an} \circ \phi_0' = \phi_0^{\rm f-an} \circ \iota$.
\item[(c)] If $r=d$, then every irreducible component $Z$ of $\tilde{\Ufrak}'$ with $\val \circ \phi_0(\xi_Z) \in \relint(\Sigma(r,\pi))$ is a toric variety.
\end{itemize}
\end{prop}

\proof Clearly, (a) follows from \eqref{affine characterization} and hence $\Ufrak'$ is obtained from $(\Ufrak')^{f-an}$ by base change to $\Tcal^{\rm f-an}$ proving also (b). Finally, (c) is a consequence of Proposition \ref{toric and etale}. \qed

\begin{rem} \rm \label{additional remark}
We assume $r=d$ and hence $\Scal \cong \Ucal_{\Sigma(d,\pi)}$. The irreducible components $Z$ in (c) are the irreducible components of $\tilde{\Xfrak}'$ contracting to the distinguished singularity $\tilde{P}$ from \ref{semistable situation} under the canonical morphism $\tilde{\iota}:\tilde{\Xfrak}' \rightarrow \tilde{\Xcal}'$. Moreover, $Z$ is isomorphic to the toric variety $Y_{\ub'}$ in $\tilde{\Tcal}$ associated to the vertex $\ub':= \val \circ \phi_0(\xi_Z)$ of $(f_{\rm aff}^{(0)})^{-1}(\Ccal)$. By Proposition \ref{toric and etale} again, we get a bijective correspondence between the above $Z$ and  vertices of $(f_{\rm aff}^{(0)})^{-1}(\Ccal)$ contained in $\relint(\Sigma(d,\pi))$. The behaviour of $\tilde{\varphi}:Z \rightarrow \tilde{\Acal}$ with respect to the torus actions is 
\begin{equation} \label{torus equivariance}
\tilde{\varphi}(\tb \cdot z) = (\tb^{\mb_1}, \dots , \tb^{\mb_n}) \cdot \tilde{\varphi}(z)
\end{equation}
for $z \in Z$ and $\tb \in ({\mathbb G}_m^d)_{\tilde{\kdop}}$. This follows from the description \eqref{semistable f-coord} of $F$ and $\tilde{u}_j|_Z \equiv \tilde{u}_j(\tilde{P})$ for $u_j \in \Ocal(\Ucal')^\times$ (use $U_+'(\tilde{P}) \cong \Scal^{\rm an}_+(\tilde{\mathbf 0})$ and the proof of Proposition \ref{semistable properties}).
\end{rem}

\begin{art} \rm \label{strata}
Since $\tilde{\Xcal}'$ is a strictly semistable formal scheme, its special fibre $\tilde{\Xcal}'$ has a canoncial stratification: Let $\tilde{\Xcal}'(i)$ be the closed subvariety of points in $\tilde{\Xcal}'$ which are contained in at least $i+1$ irreducible components of $\tilde{\Xcal}'$.  Then the irreducible components of the disjoint sets $\tilde{\Xcal}'(i) \setminus \tilde{\Xcal}'(i+1)$ are called the {\it strata} of $\tilde{\Xcal}'$.
\end{art}

For $\tilde{P} \in \tilde{\Xcal}'(\tilde{\kdop})$, let $\Ucal'$ be a formal neighbourhood in $\Xcal'$ as in \ref{semistable situation} leading to an affine map $f_{\rm aff}: \Delta(r_i,\pi_i) \rightarrow \rdop^n$.

\begin{prop} \label{affinization and strata}
The map $\overline{f}_{\rm aff}$ is determined by the stratum containing $\tilde{P}$ up to permutation of the coordinates on $\Delta(r,\pi)$. 
\end{prop}

\proof We consider $\tilde{P}_1, \tilde{P}_2 \in \tilde{\Xcal}'(\tilde{\kdop})$ in the same stratum with corresponding affine maps $f_{i, \rm aff}: \Delta(r_i,\pi_i) \rightarrow \rdop^n$. Note that $r_i+1$ is the number of irreducible components in $\tilde{\Xcal}'$ passing through $\tilde{P}_i$, hence $r_1=r_2$. By interchanging a suitable  $\tilde{P}_3$, we may assume $\Ucal_1' = \Ucal_2'$. After a permutation, Proposition \ref{semistable properties}(d) yields that $\tilde{\phi}^*_1(\tilde{x}_j')$ is equal to $\tilde{\phi}^*_2(\tilde{x}_j')$ up to a unit on $\tilde{\Ucal}'_1$. The latter lifts to a unit on $\Ucal_1'$ and hence $\val \circ \phi_1=\val \circ \phi_2$ on $U_1'=U_2'$. By \eqref{affine characterization} and $\val \circ \phi_i(U_i')= \Delta(r,\pi)$, we deduce $\overline{f}_{1,\rm aff}=\overline{f}_{2,\rm aff}$. \qed

\begin{art} \label{alteration} \rm  
Up to now, we assume that $f$ is an {\it alteration} of $X^{\rm an}$, i.e. $f$ is a proper surjective morphism $X' \rightarrow X^{\rm an}$ for an irreducible variety $X'$ of dimension $d:=\dim(X)$. By the GAGA-principle (see \cite{Ber}, 3.4.7), everything may be formulated algebraically. We assume that $X'$ has a strictly semistable formal $\kdop^\circ$-model $\Xcal'$.

From \ref{formal refinement}, we get a morphism $\varphi:\Xfrak' \rightarrow \Xcal^{f-an}$ of formal analytic varieties with $\varphi^{\rm an}=f$, where $\Xcal$ is the closure of $X$ in $\Acal$. Since $f$ is proper,  $\tilde{\varphi}$ is also proper (\cite{Gu2}, Remark 3.14). 
\end{art}

\begin{art} \label{non-degenerate alteration} \rm 
For every stratum $S$ of $\tilde{\Xcal}'$, we get a map $\overline{f}_{\rm aff}:\Delta(r_S,\pi_S) \rightarrow \rdop/\Lambda$, determined up to permutation. For $\overline{\rho}_S:=\overline{f}_{\rm aff}(\Delta(r_S,\pi_S))$, the proof of Proposition \ref{polyhedron} and Proposition \ref{affinization and strata} show 
$$ \valbar({X^{\rm an}}) = \bigcup_S \overline{\rho}_S.$$
Moreover, by Theorem \ref{periodic main theorem}, we may restrict $S$ to the strata with $\dim(\rho_S)=d$. 

We call $\ubb \in \valbar({X^{\rm an}})$ {\it non-degenerate} with respect to $f$ if 
$\ub \not \in \rho_S$ for all simplices $\rho_S$ of dimension $<d$. Note that $\ub$ and $\rho_S$ are only determined up to $\Lambda$-translation.
\end{art}

\begin{art} \rm \label{list of singularities}
Since $\Xcal'$ is a strictly semistable formal $\kdop^\circ$-model of the $d$-dimensional irreducible variety $X'$,  the set $\tilde{\Xcal}'(d)$ of strata introduced in \ref{strata} is zero-dimensional. Let $\tilde{P}_1, \dots , \tilde{P}_R$ be the points  of 
$\tilde{\Xcal}'(d)$. We denote the affine map $\Delta(d,\pi_j) \rightarrow \rdop^n$ corresponding to the stratum $\tilde{P}_j$ by $f_{j,\rm aff}$. The image of $f_{j,\rm aff}$ is a simplex in $\rdop^n$ which we denote by $\rho_j$. After renumbering, we may assume that $f_{j,\rm aff}$ is one-to-one exactly for $j=1, \dots, N$. By \ref{non-degenerate alteration}, we have the decomposition
$$\valbar({X^{\rm an}}) = \bigcup_{j=1}^N \overline{\rho}_j.$$ 
The lower dimensional simplices $\rho_{N+1}, \dots, \rho_R$ will play only a minor role in the following. For $j=1, \dots, N$, the bijective projection $\Delta(d,\pi_j) \rightarrow \Sigma(d,\pi_j)$ and $f_{j,\rm aff}$ induce  $f_{j, \rm aff}^{(0)}: \Sigma(d,\pi_j) \rightarrow \rdop^n$ (see \ref{affine map}) which extends canonically to an injective affine map $\rdop^d \rightarrow \rdop^n$ also denoted by $f_{j, \rm aff}^{(0)}$.
\end{art}

\begin{art} \rm \label{transverse and non-degenerate assumption}
We consider a polytope $\overline{\Delta} \in \Ccalbar$ of codimension $d$ with relative interior $\overline{\tau}$. We assume that $\Deltabar \cap \valbar({X^{\rm an}})$ is a non-empty finite subset of $\overline{\tau}$. We suppose also that the points of $\Deltabar \cap \valbar({X^{\rm an}})$ are non-degenerate with respect to $f$.
\end{art}

\begin{art} \rm \label{index definition}
For $j \in \{1, \dots ,N\}$ with $\overline{\rho}_j \cap \Deltabar \neq \emptyset$, we are going to define an index of $\Deltabar$ relative to $\overline{f}_{j, \rm aff}$. Note that $j \leq N$ means that the simplex $\rho_j$ from \ref{list of singularities} is $d$-dimensional and hence $\overline{\rho}_j \cap \Deltabar$ is a transversal intersection. The index will depend only on the $d$-codimensional linear subspace $\ldop_\Delta$ of $\rdop^n$ with $\Delta \subset \ub + \ldop_\Delta$ and on the injective linear map $\ell_{j}^{(0)}:=f_{j, \rm aff}^{(0)} - f_{j,\rm aff}^{(0)}({\mathbf 0}): \rdop^d \rightarrow \rdop^n$. 

Note that $\ldop_\Delta$ is defined over $\qdop$ and hence $N_\Delta:= \ldop_\Delta \cap \zdop^n$ is a complete lattice in $\ldop_\Delta$. Let $q_\Delta$ be the quotient map $\rdop^n \rightarrow \rdop^n/\ldop_\Delta$. Since $\Deltabar$ is transversal to $\overline{\rho}_j$ and since $\ell_{j}^{(0)}$ is injective, $q_\Delta \circ \ell_{j}^{(0)}$ is also injective on $\rdop^d$. Using that $\ell_{j}^{(0)}$ is defined over $\zdop$, we conclude that $(q_\Delta \circ \ell_{j}^{(0)})(\zdop^d)$ is a complete lattice in $\rdop^n/\ldop_\Delta$ of finite index in $\zdop^n/N_\Delta$. This crucial index will be denoted by
$$\ind(\Deltabar,\overline{f}_{j, \rm aff}):= \left[ \zdop^n/N_\Delta : (q_\Delta \circ \ell_{j}^{(0)})(\zdop^d)\right].$$ 
It is important to note that $\ind(\Deltabar, \overline{f}_{j, \rm aff})$ depends only on $\ldop_\Delta$ and on $\ell_{j}^{(0)}$.   
The index may be more canonically described in terms of the linear map $\ell_{j} := f_{j, \rm aff}-f_{j, \rm aff}(\mathbf 0)$ defined on the hyperplane $\{\ub' \in \rdop^{d+1} \mid u_0' + \dots + u_d' = 0\}$. 

In the above definition of the index, not all assumptions in \ref{transverse and non-degenerate assumption} are needed. In fact, $\ind(\Deltabar, \overline{f}_{j, \rm aff})$ is defined for all $j \in \{1, \dots, N\}$ and every $\Deltabar \in \Ccalbar$ of codimension $d$ with $\ldop_\Delta \cap \ldop_{\rho_j} = \{0\}$. 
\end{art}

\begin{prop} \label{multiplicity formula}
Recall that $A$ is a totally degenerate abelian variety over $\kdop$ with Mumford model $\Acal$ associated to $\Ccalbar$ and that $X$ is a $d$-dimensional irreducible closed subvariety of $A$. Let $f: X' \rightarrow X$ be an alteration with strictly semistable formal $\kdop^\circ$-model $\Xcal'$ of $X'$.  
Let $\Deltabar$ be a polytope from $\Ccalbar$ satisfying the transversality assumption  \ref{transverse and non-degenerate assumption}. 
\begin{itemize}
\item[(a)] There is a unique formal analytic structure $\Xfrak'$ on $X'$ which refines $(\Xcal')^{\rm f-an}$ such that $f$ extends to a morphism $\varphi:\Xfrak' \rightarrow \Acal^{\rm f-an}$.
\item[(b)] The toric variety $Y:=Y_{\Deltabar}$ associated to $\Deltabar$ (see Remark \ref{polytopal toric variety}) is an irreducible component of $\tilde{\Xcal}$ for the closure $\Xcal$ of $X$ in $\Acal$.
\item[(c)] Let $\tilde{P}_1, \dots , \tilde{P}_N$ be the zero-dimensional strata of $\tilde{\Xcal}'$ such that the associated affine maps $f_{j, \rm aff}: \Delta(d,\pi_j) \rightarrow \rdop^n$ are one-to-one (see \ref{list of singularities}). An irreducible component $Z$ of $\tilde{\Xfrak}'$ with $\tilde{\varphi}(Z)=Y$ contracts to a unique $\tilde{P}_j$ with respect to $\tilde{\Xfrak}' \stackrel{\tilde{\iota}}{\rightarrow} \tilde{\Xcal}'$ for some $j \in \{1, \dots ,N\}$.
\item[(d)] This gives a bijective correspondence between such $Z$ and $J:=\{j \in \{1, \dots, N\} \mid \Deltabar \cap \overline{f}_{j, \rm aff}(\Delta(d, \pi_j)) \neq \emptyset \}$.
\item[(e)] For $j=1, \dots, N$, let $f_{j, \rm aff}^{(0)}: \Sigma(d,\pi_j) \rightarrow \rdop^n$ be the affine map induced from $f_{j, \rm aff}$ as in \ref{affine map}. We choose a formal affine neighbourhood $\Ucal_j'$ of $\tilde{P}_j$ in $\Xcal'$  and an \'etale morphism $\phi_{j,0}:\Ucal_j' \rightarrow \Ucal_{\Sigma(d,\pi_j)}$ as in \ref{semistable situation} (for $\tilde{P} = \tilde{P}_j$). 

If $Z$ corresponds to $j$ by (d), then
$\ub' := \val \circ \phi_{j,0}(\xi_Z)$ is a vertex of $\left(f_{j, \rm aff}^{(0)} \right)^{-1} ( \Ccal)$ and $Z$ is isomorphic to the toric variety $Y_{\ub'}$  associated to $\ub'$ (see Proposition \ref{torus orbits global}).
\item[(f)] If $Z$ is as in (c), then $[Z:Y] = \ind(\Deltabar, \overline{f}_{j, \rm aff})$.
\item[(g)] The multiplicity of $Y$ in $\tilde{\Xcal}$ satisfies
$m(Y,\tilde{\Xcal})= \frac{1}{[X':X]} \sum_{j \in J} \ind(\Deltabar,\overline{f}_{j,\rm aff})$.
\end{itemize}
\end{prop}

\proof Obviously, (a) is a reformulation of \ref{formal refinement}. The transversality assumption on $\Deltabar$ yields (b) by Lemma \ref{transverse lemma}. Let $Z$ be an irreducible component of $\tilde{\Xfrak}'$ with $\tilde{\varphi}(Z)=Y$. Then $\tilde{\varphi}(\tilde{\xi}_Z)$ is the generic point of $Y$, hence $f(\xi_Z)$ is in the open dense orbit $Z_{\overline{\tau}}$ of $Y$ for $\overline{\tau}:=\relint(\Deltabar)$. We conclude that $\ubb:=\valbar(f(\xi_Z))$ is contained in the finite set 
$$\Deltabar \cap \valbar({X^{\rm an}}) = \overline{\tau} \cap \valbar({X^{\rm an}}).$$
There is a formal affine open subset $\Ucal'$ of $\Xcal'$ as in \ref{semistable situation} with $\xi_Z $ contained in the generic fibre $U'$. The image of the corresponding affine map $f_{\rm aff}$ is a simplex $\rho$ with $\ubb \in \overline{\rho}$ by \eqref{affine characterization}, hence non-degeneracy of $\Deltabar \cap \valbar({X^{\rm an}})$ with respect to $f$ yields that $\rho$ is $d$-dimensional. By Proposition \ref{affinization and strata}, $\tilde{P}_j \in \Ucal'$ for some $j \in \{1, \dots, N\}$ and so we may assume $\Ucal'=\Ucal_j'$, $f_{\rm aff}= f_{j, \rm aff}$ (see \ref{list of singularities}). Moreover, \eqref{affine characterization} and non-degeneracy show that
$$\ub' := \val \circ \phi_{j,0}(\xi_Z) \in \relint(\Sigma(d, \pi_j))$$
and $\ubb = \overline{f}_{j, \rm aff}^{(0)}(\ub')$. For the polytopal decompositon $\Ccal_j:=(f_{j, \rm aff}^{(0)})^{-1}(\Ccal)$ of $\Sigma(d,\pi_j)$, we have seen in Remark  \ref{additional remark} that there is a bijective correspondence between irreducible components of $\tilde{\Xcal}'$ contracting to $\tilde{P}_j$ and vertices of  $\Ccal_j$ contained in $\relint(\Sigma(d,\pi))$. This leads to (c) and (e). 

To prove (d), we have to show  for $j \in J$ that there is a unique irreducible component $Z$ of $\tilde{\Xfrak}'$ contracting to $\tilde{P}_j$ with $\tilde{\varphi}(Z)=Y$. The $d$-dimensional simplex $\rho_j := f_{j, \rm aff}(\Delta(d,\pi_j))$ satisfies $\overline{\rho}_j \subset \valbar({X^{\rm an}})$ by \ref{list of singularities}. The assumption $\Deltabar \cap \overline{\rho}_j \neq \emptyset$ leads to a lift $\Delta \subset \rdop^n$ with $\Delta \cap \rho_j \neq \emptyset$. By transversality, $\Delta \cap \rho_j = \tau \cap \rho_j$ consists of a single point $\ub$. Since $f_{j, \rm aff}^{(0)}$ is injective, there is a unique $\ub' \in \Sigma(d,\pi_j)$ with $f_{j, \rm aff}^{(0)}(\ub') = \ub$. We note that $\ub'$ is a vertex of $\Ccal_j$. 
Since $\ubb$ is non-degenerate with respect to $f$, we have $\ub' \in \relint(\Sigma(d, \pi_j))$. As we have seen above, this vertex $\ub'$ corresponds to a unique irreducible component $Z$ of $\tilde{\Xcal}'$ contracting to $\tilde{P}_j$. 
By $\val(\phi_{j,0}(\xi_Z))=\ub'$ and \eqref{affine characterization}, we get $\valbar \circ f(\xi_Z) = \ubb$. Proposition \ref{torus orbit 2} proves
$$\tilde{\varphi}(\tilde{\xi}_Z)=(f(\xi_Z))\sptilde \in Z_{\overline{\tau}}$$
and hence $\tilde{\varphi}(Z) \subset Y$. The above application of Remark \ref{additional remark} shows more precisely that $Z$ is isomorphic to $Y_{\ub'}$ and therefore $Z$ is a toric variety with respect to the induced $({\mathbb G}_m^d)_{\tilde{\kdop}}$-action. On the other hand, $Y$ is a toric variety with respect   to the torus $T$ over $\tilde{\kdop}$ given by 
$$T(\tilde{\kdop})=\left(\zdop^n/N_\Delta \right) \otimes_\zdop \tilde{\kdop}^\times$$
(see Remark \ref{polytopal toric variety}). Here and in the following, we use the notation from \ref{index definition}.  For $i=1, \dots, n$, there is $\mb_i \in \zdop^d$ and $\lambda_i \in \kdop^\times$ such that the $i$-th coordinate of $f_{j, \rm aff}^{(0)}(\ub'')$ is equal to $\mb_i \cdot \ub'' + v(\lambda_i)$ (see \eqref{explicit affinization}).   Let $\nu:({\mathbb G}_m^d)_{\tilde{\kdop}} \rightarrow T$ be the composition of $\tb \mapsto (\tb^{\mb_1}, \dots, \tb^{\mb_n}) \in (\Tor)_{\tilde{\kdop}}$  with the quotient homomorphism $(\Tor)_{\tilde{\kdop}} \rightarrow T$. Then the homomorphism $\nu$ is induced by the linear map $q_\Delta \circ \ell_{j}^{(0)}:  \zdop^d \rightarrow \zdop^n/N_\Delta$. Since $q_\Delta \circ \ell_{j}^{(0)}$  is one-to-one (see \ref{index definition}), we deduce that $\nu$ is a finite surjective homomorphism of degree $\ind(\Deltabar, \overline{f}_{j, \rm aff})$. Now \eqref{torus equivariance} yields
\begin{equation} \label{equation of equivarence}
\tilde{\varphi}(\tb \cdot z)= \nu(\tb) \tilde{\varphi}(z) \quad (\tb \in ({\mathbb G}_m^d)_{\tilde{\kdop}}, z \in Z)
\end{equation}
and we conclude that $\tilde{\varphi}(Z)=Y$. Uniqueness of $Z$ is clear from the construction and hence we  get (d). Moreover, (f) follows easily from our des\-cription of $\nu$ and \eqref{equation of equivarence}. Finally, (g) is a consequence of (d), (f) and the projection formula
$$m(Y,\tilde{\Xcal})= \frac{1}{[X':X]} \sum_Z [Z:Y],$$
where $Z$ ranges over all irreducible components of $\tilde{\Xfrak}'$ with $\tilde{\varphi}(Z)=Y$. 
\qed

\section{Canonical measures}

In this section, $K$ is a field with a discrete valuation $v$. The completion of the algebraic closure of the completion of $K$ with respect to $v$ is an algebraically closed field denoted by $\kdop$ (\cite{BGR}, Proposition 3.4.1/3). The unique extension of $v$ to a valuation of $\kdop$ is also denoted by $v$ with corresponding absolute value $|\phantom{a}| := e^{-v}$. The value group $\Gamma = v(\kdop^\times)$ is equal to $\qdop$.

Let $A$ be an abelian variety over $K$ which is totally degenerate over $\kdop$, i.e. $A^{\rm an}= (\Tor)_\kdop^{\rm an}/M$ for a lattice $M$ isomorphic to $\Lambda := \val(M)$ in $\rdop^n$. Let $X$ be a geometrically integral $d$-dimensional closed subvariety of $A$. 

We will show first that a generic rational polytopal decomposition $\Ccalbar$ of $\rtor$ is transversal to $\valbar({{\Xan}})$. If $\Acal$ denotes the associated formal $\kdop^\circ$-model of $A$ and if $\Xcal$ is the closure of ${\Xan}$ in $\Acal$, then transversality allows us to handle the special fibre of $\Xcal$ by the theory of toric varieties. By rationality, the decomposition $\overline{\frac{1}{m} \Ccal} := \{\overline{\frac{1}{m} \Delta} \mid \Delta \in \Ccal \}$ of $\rtor$ can not be transversal  to $\valbar(\Xan)$ for all $m \in \ndop \setminus \{0\}$ simultaneously. However, this may be achieved over a sufficiently large base extension $\kdop'$ with value group $\Gamma'$ by a ``completely irrational'' construction which is also suitable for extending an ample line bundle on $A$ to a  positive formal $(\kdop')^\circ$-model on $\Acal$. This will be used to prove Theorem \ref{Theorem 3}.
Moreover, we will get an explicit formula for the canonical measures on $\val(\Xan)$ given in Theorem \ref{explicit expression}.

There is no restriction of generality assuming that $X$ is geometrically integral as we may proceed by  base change and linearity to get the canonical measures in the general case. 

\begin{art} \rm \label{generic decompositions}
Fix some $\Lambda$-periodic set $\Sigma$ of polytopes in $\rdop^n$ and set $\overline{\Sigma} := \{ \overline{\sigma} \subset \rtor \mid \sigma \in \Sigma \}$. We assume that $\overline{\Sigma}$  is a finite set. If a polytope is in $\Sigma$, then we require that all its closed faces are also in $\Sigma$. For a polytope $\sigma \in \rdop^n$, $\adop_\sigma$ denotes the affine space in $\rdop^n$ generated by $\sigma$.

The polytopal decomposition $\Ccalbar$ of $\rtor$ is said to be {\it generic} with respect to $\overline{\Sigma}$ if the following conditions hold for every $\sigma \in \Sigma$,  $\Delta \in \Ccal$:
\begin{itemize}
\item[(a)] $\dim(\adop_\sigma \cap \adop_\Delta) = D$ if $D:=\dim(\sigma)+\dim(\Delta)-n \geq 0$,
\item[(b)] $\adop_\sigma \cap \adop_\Delta = \emptyset$ if $D <0$.
\end{itemize}
A polytopal decomposition $\Ccalbar$ of $\rtor$ is called {\it $\Sigmabar$-transversal} if $\Delta \cap \sigma$ is either empty or of dimension $\dim(\Delta) + \dim(\sigma) - n$ for all $\Delta \in \Ccal$, $\sigma \in \Sigma$. If the  union $\overline{S}$ of all polytopes in $\Sigmabar$ is pure dimensional, then a $\Sigmabar$-transversal $\Ccalbar$ is transversal to $\overline{S}$ in the sense of \ref{local cone and transversality}. 
\end{art}

\begin{prop} \label{generic is transversal}
Every $\Sigmabar$-generic polytopal decomposition of $\rtor$ is $\Sigmabar$-transversal.
\end{prop}

\proof Let $\Delta$, $\Sigma$ be polytopes in $\rdop^n$ with $\Delta \cap \Sigma \neq \emptyset$ such that all closed faces $\sigma'$ of $\sigma$ and $\Delta'$ of $\Delta$ satisfy (a) and (b). It is enough to show $\dim(\Delta \cap \sigma) = \dim(\Delta) + \dim(\sigma) - n$. If $\relint(\Delta) \cap \relint(\sigma) \neq \emptyset$, then this is obvious from (a). So we may assume $\relint(\Delta) \cap \relint(\sigma) = \emptyset$ which will lead to a contradiction. By symmetry and passing to closed faces if necessary, we may assume that there is a closed face $\sigma'$ of codimension $1$ in $\sigma$ such that $\Delta \cap \sigma' = \Delta \cap \sigma$ and 
$\relint(\Delta) \cap \relint(\sigma') \neq \emptyset$. 
Note that $\adop_{\sigma'}$ divides $\adop_\sigma$ into half spaces, one is containing $\sigma$. Since $\relint(\Delta) \cap \relint(\sigma') \neq \emptyset$, we conclude that $\adop_\Delta \cap \adop_\sigma$ is contained in $\adop_{\sigma'}$. Thus $\adop_\Delta \cap \adop_{\sigma'}=\adop_\Delta \cap \adop_\sigma$ contradicts (a) and (b). \qed

\begin{art} \rm \label{rationality} 
Starting with an arbitrary rational triangulation of $\rtor$ and varying the vertices a little bit in $\qdop^n$, we get a rational triangulation $\Ccalbar$ of $\rtor$ which is $\Sigmabar$-generic. 

Up to now, we assume that $\Sigma$ is rational. For the proof of Theorem \ref{Theorem 3}, we need that $\overline{\frac{1}{m} \Ccal}$ is $\Sigmabar$-generic for all non-zero $m \in \ndop$ simultaneously which is not possible for a rational $\Ccalbar$. Instead, we are working with an infinite dimensional $\qdop$-subspace $\Gamma'$ of $\rdop$ containing $\qdop$. By \cite{Bou}, Ch. VI, $n^\circ$ 10, Prop. 1, there is an algebraically closed field $\kdop'$, complete with respect to an absolute value $|\phantom{a}|'$ extending $|\phantom{a}|$ such that $v'((\kdop')^\times) = \Gamma'$. Now we will see that a $\Gamma'$-rational $\Ccalbar$ with the desired property and which behaves well with respect to the extension of an ample line bundle can be obtained by a ``completely irrational'' construction. Since the lemma will be important for the sequel, we give a rather detailed proof.
\end{art}

\begin{lem} \label{irrational construction}
Let $L$ be an ample line bundle on $A$. Then there is a $\Gamma'$-rational polytopal decomposition $\Ccalbar$ of $\rtor$ with the following properties:
\begin{itemize}
\item[(a)] $\overline{\frac{1}{m} \Ccal}$ is $\overline{\Sigma}$-generic and hence $\Sigmabar$-transversal for all $m \in \ndop \setminus \{0\}$.
\item[(b)] If $\Acal$ denotes the formal $(\kdop')^\circ$-model of $A_{\kdop'}^{\rm an}$ associated to $\Ccalbar$, then there are $N \in \ndop \setminus \{0\}$ and a formal $(\kdop')^\circ$-model $\Lcal$ of $L^{\otimes N}$ on $\Acal$ with $\tilde{\Lcal}$ ample. 
\end{itemize}
\end{lem}

\proof In the terminology of \ref{toric line bundles} and of Proposition \ref{Mumford's line bundles},  $L$ induces $z_\lambda$ and a bilinear form $b$ on $\Lambda$. 
By \ref{toric line bundles}, we deduce that $z_\lambda(\mathbf 0)$ is a quadratic function on $\Lambda$ and therefore 
$$z_\lambda(\mathbf 0) = q(\lambda) + \ell(\lambda)$$
for the quadratic form $q(\lambda):=\frac{1}{2}b(\lambda,\lambda)$ and a linear form $\ell$ on $\Lambda$ (see \cite{BG}, 8.6.5).  Both extend to corresponding forms on $\rdop^n$ also denoted by $q$ and $\ell$. Since $L$ is ample, $q$ is positive definite on $\Lambda$ (see \ref{toric line bundles}) and hence its extension to $\rdop^n$ is also positive definite (using $\Gamma =\qdop$).  We conclude that $f:=q+\ell$ is a strictly convex function on $\rdop^n$ (see \ref{strongly polyhedral function}). Formula \eqref{cycle and bilinear} yields
\begin{equation} \label{automorphy}
f(\ub + \lambda) = f(\ub) + z_\lambda(\ub) \quad (\lambda \in \Lambda, \, \ub \in \rdop^n).
\end{equation}

Our goal is to construct a $\Gamma'$-rational polytopal decomposition $\Ccalbar$ of $\rtor$ with (a) and a strongly polyhedral convex function $g: \rdop^n \rightarrow \rdop$  with respect to $\Ccal$ such that \eqref{automorphy} holds for $g$ and such that for every $\Delta \in \Ccal$, there are $\mb_\Delta \in \qdop^n$, $c_\Delta \in \Gamma'$ with 
\begin{equation} \label{affine equation}
 g(\ub) = \mb_\Delta \cdot \ub + c_\Delta \quad (\ub \in \Delta).
\end{equation}
We show first that this implies the lemma. By \ref{toric line bundles}, there is $\mb_\lambda \in \zdop^n$, additive in $\lambda \in \Lambda$, such that $b(\lambda, \ub)=\mb_\lambda \cdot \ub$. Now \eqref{automorphy}, \eqref{affine equation} and \eqref{cycle and bilinear} yield
\begin{equation} \label{m-transform}
 \mb_{\Delta+\lambda} = \mb_\Delta+\mb_\lambda
\end{equation}
and hence there is a common denominator $N$ of all $\mb_\Delta$, $\Delta \in \Ccal$. By Proposition \ref{Mumford's line bundles}, there is a formal $\kdop^\circ$-model $\Lcal$ of $L^{\otimes N}$ on $\Acal$ with $f_{\Lcal} = N g$. Since $g$ is a strongly polyhedral convex function with respect to $\Ccal$, Corollary \ref{ample reduction} yields that $\tilde{\Lcal}$ is ample. 

Before we start the construction, we note that we may assume $\ell=0$: This corresponds to the replacement of $L$ by $L \otimes [-1]^*(L)$. If $\widehat{g}$ is a solution for the latter, then $g:=\frac{1}{2} \widehat{g} + \ell$ is a solution for the original problem (use \ref{toric line bundles}). So we may assume $f=q$. 

Let $e_1', \dots, e_n'$ be a basis of $\Lambda$ with fundamental domain $F_\Lambda:= \sum_{i=1}^n [0,1) e_i'$. We number the $r=2^n$ points
$$\theta_1 e_1' + \dots + \theta_n e_n' \quad (\theta_k \in \{0, \frac{1}{2}\})$$
by $\ub_1, \dots , \ub_r$.  They form the set $F_\Lambda \cap \frac{1}{2} \Lambda$. We have the affine approximation
\begin{equation} \label{Aij}
 A_i(\ub) := b(\ub,\ub_i)-q(\ub_i)
\end{equation}
of $q$ in $\ub_i$. We have $A_i(\ub_i)=q(\ub_i)$ and $A_i(\ub) < q(\ub)$ for $\ub \neq \ub_i$ by strict convexity of $q$. By periodicity, we extend these definitions to $\rdop^n$, i.e. we set for $\lambda \in \Lambda$:
\begin{equation} \label{uil}
\ub_{i, \lambda}:= \ub_i + \lambda, \quad A_{i, \lambda}(\ub) := b(\ub, \ub_{i,\lambda}) - q(\ub_{i, \lambda}).
\end{equation}
Then the $\ub_{i,\lambda}$ form just $\frac{1}{2} \Lambda$ and $A_{i, \lambda}$ is the affine approximation of $q$ in $\ub_{i, \lambda}$. The strongly polyhedral convex function 
\begin{equation} \label{gj}
 g:= \max_{i=1, \dots, r; \; \lambda \in \Lambda} A_{i, \lambda}
\end{equation}
is a  lower bound of $q$ and satisfies \eqref{automorphy}. Moreover, $g$ is affine on the rational polytopes 
\begin{equation} \label{Dil}
 \Delta_{i,\lambda} := \{ \ub \in \rdop^n \mid A_{i, \lambda}(\ub)=g(\ub)\}
= \bigcap_{(h,\mu) \neq (i,\lambda)} \{A_{i,\lambda}(\ub) \geq A_{h,\mu}(\ub)\}.
\end{equation}
They are the maximal polytopes  where $g$ is affine  and  they are the $d$-dimensional polytopes of the Voronoi decomposition $\Ccal$ of $\frac{1}{2} \Lambda$ with respect to the euclidean metric $q$, i.e.
$$\Delta_{i, \lambda} = \{\ub \in \rdop^n \mid q(\ub- \ub_{i,\lambda}) \leq q(\ub - \ub_{h,\mu}) \,\; \forall h=1, \dots, r, \, \; \forall \mu \in \Lambda\}.$$
The polytopal decomposition $\Ccal$ is  $\frac{1}{2}\Lambda$-periodic and induces a rational polytopal decomposition $\Ccalbar$ of $\rtor$. By construction, \eqref{automorphy} and \eqref{affine equation} are satisfied but $\Ccalbar$ does certainly not satisfy (a). 

To achieve this,  we modify the construction by a small perturbation. The union $G_0$ of the Voronoi cells $\Delta_{i,\mathbf 0}$ $(1 \leq i \leq r)$ is the closure of a fundamental domain for $\rtor$. Then
$$G_1:= \bigcup_{(i,\lambda) \in T_1} \Delta_{i,\lambda}, \quad T_1:=\{(i,\lambda) \in \{1, \dots, r\} \times \Lambda \mid \Delta_{i,\lambda} \cap G_0 \neq \emptyset\},$$
is the set of neighbours of $G_0$. We approximate  the gradient $\nabla q(\ub_i)$ by $\mb_i \in \qdop^n$ and $-q(\ub_i)$ by $c_i \in (\Gamma')^n$. Then we define the affine function
$$A_i(\ub) := \mb_i \cdot \ub + c_i$$
which is very close to the old definition. We may still assume that $g$ is an upper bound of $A_i$. With the $\mb_\lambda \in \zdop^n$ introduced at the beginning of the proof, we define affine functions 
$$A_{i,\lambda}(\ub) := \mb_{i,\lambda} \cdot \ub + c_{i,\lambda},
\quad \mb_{i,\lambda} := \mb_i + \mb_\lambda \in \qdop^n, \quad c_{i,\lambda}:=c_i - q(\lambda)- \mb_i \cdot \lambda.$$
By construction, they are very close to the old approximations, $A_{i,\mathbf 0} = A_i$ and \eqref{cycle and bilinear} yields 
\begin{equation} \label{A-transform}
A_{h,\rho+\lambda}(\ub + \lambda) = A_{h,\rho}(\ub) + z_\lambda(\ub) \quad (\lambda, \mu, \rho \in \Lambda).
\end{equation}

We assume that the approximations $\mb_{i,\lambda}$, $c_{i,\lambda}$ of $\nabla q(\ub_i)$, $-q(\ub_i)$ satisfy the conditions:
\begin{itemize}
\item[(c)] $1, (c_i)_{i=1,\dots,r}$ are $\qdop$-linearly independent in $\Gamma'$.
\item[(d)] For all $\sigma \in \Sigma$, given by linearly independent equations $\lb_1 \cdot \ub =a_1, \dots , \lb_c \cdot \ub = a_c$, with $\lb_k \in \zdop^n$ and $a_k \in \zdop$ (possible by rationality of $\Sigma$), and for every $S \subset \{1, \dots, r\}$ with $\card(S)  \leq n-c+1$, $\lambda:S \rightarrow \Lambda \cap (G_0-G_0)$, the vectors 
$$\left(\mb_{i,\lambda(i)} - \mb_{i_0,\lambda(i_0)} \right)_{i \in S \setminus\{i_0\}} , \lb_1, \dots ,\lb_c$$
are linearly independent, where $i_0$ is the minimal member of $S$.
\end{itemize}

The existence of $c_i$  follows from $[\Gamma': \qdop] = \infty$ and the construction of the approximations $\mb_i$ is by induction where in each step, we have to omit finitely many hyperplanes which is possible in every neighbourhood of $\nabla q(\ub_i)$.

The function $g$ defined by \eqref{gj} is a strongly polyhedral convex function with respect to $\Ccal$ which is equal to $A_{i,\lambda}$ on the $\Gamma'$-rational polytope $\Delta_{i,\lambda}$ from \eqref{Dil}. The latter are again the $d$-dimensional polytopes of a $\Gamma'$-rational decomposition $\Ccal$ of $\rtor$. $\Ccal$ is very close to the Voronoi decompositon considered above in the sense that the boundary of the new $\Delta_{i,\lambda}$ is near to the boundary of the corresponding Voronoi cell. Moreover, if two cells intersect, then also the corresponding Voronoi cells intersect.  By \eqref{A-transform}, $\Ccal$ is $\Lambda$-periodic and $g$ satisfies \eqref{automorphy}. Note that \eqref{affine equation} is clear by construction. We get a polytopal decomposition $\Ccalbar$ of $\rtor$.

It remains to see that $\Ccalbar$ satisfies (a): Let $\sigma \in \Sigma$, $m \in \ndop \setminus \{0\}$ and $\Delta \in \Ccal$. We may represent $\Delta$ as a closed face of $\Delta_{i_0,\lambda_0}$ given by the hyperplanes
$$A_{i,\lambda}(\ub) = A_{i_0,\lambda_0}(\ub), \quad (i, \lambda) \in S \subset \{1,\dots,r\}\times \Lambda .$$
We may assume that every hyperplane is generated by a face of  $\Delta_{i_0,\lambda_0}$ and that 
\begin{equation} \label{Ad}
\adop_\Delta = \bigcap_{(i,\lambda) \in S} \left\{ A_{i,\lambda}(\ub)=A_{i_0,\lambda_0}(\ub) \right\} 
\end{equation}
is a transversal intersection. By $\Lambda$-periodicity, we may assume that 
$\Delta \subset G_1$. 
Since $\Ccal$ is very close to the Voronoi decomposition of $\frac{1}{2}\Lambda$ with respect to $q$, it is clear that for given $i$, there is at most one $\lambda \in \Lambda$ involved in \eqref{Ad}. Note that such a $\lambda$ is contained in $G_0-G_0$. By \eqref{Ad}, $\adop_\sigma \cap \adop_\Delta$ is the solution space of the $c+\codim(\Delta)$ linear equations
\begin{equation} \label{leq}
 \lb_1 \cdot \ub = a_1, \dots, \lb_c \cdot \ub = a_c
\end{equation}
(in the notation borrowed from (d)) and 
\begin{equation} \label{meq}
\left(\mb_{i,\lambda} - \mb_{i_0,\lambda_0} \right) \cdot \ub = c_{i_0,\lambda_0} - c_{i,\lambda} \quad \left((i,\lambda) \in S \right). 
\end{equation}
If $D:= \dim(\sigma) + \dim(\Delta) - n  \geq 0$, then the assumption (d) yields that the homogeneous linear equations associated to \eqref{leq} and \eqref{meq} are linearly independent and hence $\dim(\adop_\sigma \cap \adop_\Delta) = D$. If $D<0$, then we may express $\lb_1 \cdot \ub$ as a $\qdop$-linear combination of the left hand sides of the other equations in \eqref{leq} and \eqref{meq}. By (c), $a_1$ can not be the same linear combination of the right hand sides. This proves $\adop_\sigma \cap \adop_\Delta = \emptyset$ and therefore $\Ccalbar$ is $\Sigmabar$-generic.

If we replace $\Delta$ by $\frac{1}{m}\Delta$, then the right hand side of \eqref{meq} is multiplied by $\frac{1}{m}$. This does not change the above argument and we conclude that $\overline{\frac{1}{m}\Ccal}$ is $\Sigmabar$-generic. \qed

\begin{art} \rm \label{setting}
Let $\overline{L}_1, \dots , \overline{L}_d$ be ample  line bundles on $A$ endowed with canonical metrics. Our goal is to give an explicit formula for the canonical measure
$$\mu := (\valbar)_*\left(c_1(\overline{L}_1|_X) \wedge \cdots \wedge c_1(\overline{L}_d|_X)\right)$$ 
on the tropical variety $\valbar(X_\kdop^{\rm an})$. 
By de Jong's alteration theorem (\cite{dJ}, Theorem 6.5), there is an alteration $f:X' \rightarrow X_\kdop^{\rm an}$ and a strictly semistable formal $\kdop^\circ$-model $\Xcal'$ of $X'$ (see \ref{alteration}). 
Let $\tilde{P}_1, \dots , \tilde{P}_N$ be the zero-dimensional strata of $\tilde{\Xcal}'$ such that the associated affine maps $f_{j, \rm aff}: \Delta(d,\pi_j) \rightarrow \rdop^n$ are one-to-one (see \ref{list of singularities}). The image of $f_{j, \rm aff}$ is a $d$-dimensional simplex denoted by $\rho_j$. The sets $\overline{\rho}_1, \dots , \overline{\rho}_N$ cover $\valbar(\Xan)$ (see \ref{list of singularities}). 

For simplicity of notation, we may assume that $\overline{\rho}_1, \dots , \overline{\rho}_N$ are simplices in $\rtor$, i.e. the projection $\rho_j \rightarrow 
\overline{\rho}_j$ is bijective. In general, a subdivision of the $\rho_j$ may be needed but this does not change the description of the  measure $\mu$ in Theorem 8.6.


Let $\overline{\sigma}$ be an atom of the covering $\bigcup_j \overline{\rho}_j=\valbar(\Xan)$, i.e. $\overline{\sigma}$ is the closure of $\bigcap_j \overline{T}_j$, where $\overline{T}_j$ is either $\relint(\overline{\rho}_j)$ or $\valbar(\Xan) \setminus \overline{\rho}_j$. We omit $\sigmabar = \emptyset$. Then $\overline{\sigma}$ is a finite union of $d$-dimensional polytopes in $\rtor$. Moreover, the sets $\overline{\sigma}$ form a finite covering of $\valbar(\Xan)$ with overlappings of dimension $<d$. To get polytopes in $\rtor$, a subdivision would be needed, but this is irrelevant for the description of the measure $\mu$. 

Since $\sigmabar \neq \emptyset$, the set $J(\sigmabar):=\{j=1,\dots,N \mid \sigmabar \subset \overline{\rho}_j\}$ is non-empty. Let $\adop_\sigma$ be the affine space generated by $\rho_j$ for some $j \in J(\sigmabar)$. Up to $\Lambda$-translation, it is independent of the choice of $j \in J(\sigmabar)$. We may lift a measurable subset $\Omegabar$ of $\sigmabar$ to  $\Omega \subset \rho_j$. We conclude that the relative Lebesgue measure on $\adop_\sigma$ induces a canonical measure on $\sigmabar$ which we denote by $\vol(\Omegabar)$. 

Using $\Gamma = \qdop$, we deduce that the stabilizer $\Lambda(\adop_\sigma):= \{\lambda \in \Lambda \mid \adop_\sigma + \lambda \subset \adop_\sigma \}$ of $\adop_\sigma$ is a complete lattice in the linear space $\ldop_\sigma$ parallel to $\adop_\sigma$.
For $j \in J(\sigmabar)$, let $f_{j, \rm aff}^{(0)}: \Sigma(d,\pi_j) \rightarrow \rdop^d$ be the affine map induced from $f_{j, \rm aff}$ as in \ref{affine map}. Then the linear map  $\ell_{j}^{(0)} := f_{j, \rm aff}^{(0)} - f_{j, \rm aff}^{(0)}(\mathbf 0)     :\rdop^d \rightarrow \rdop^n$ is one-to-one, defined over $\zdop$ and has image $\ldop_\sigma$. We get a dual map $\hat{\ell}_{j}^{(0)}: \ldop_\sigma^* \rightarrow (\rdop^d)^*=\rdop^d$ which is bijective and defined over $\zdop$. 

Let $b_j$ be the bilinear form associated to $L_j$ (see \ref{toric line bundles}). Since $b_j$ is positive definite on $\rdop^n$,
$$\Lambda(\adop_\sigma)^{L_j}:=\{b_j(\cdot,\lambda) \in \ldop_\sigma^* \mid \lambda \in \Lambda(\adop_\sigma)\}$$
is a complete lattice in $\ldop_\sigma^*$. We will also use the dual lattice of $\Lambda(\adop_\sigma)$, given by
$$\Lambda(\adop_\sigma)^* := \{ \ell \in \ldop_\sigma^* \mid \ell(\Lambda(\adop_\sigma)) \subset \zdop\}.$$

We denote  by $\vol(\Lambda(\adop_\sigma))$ the volume of a fundamental domain of the lattice $\Lambda(\adop_\sigma)$ with respect to the relative Lebesgue measure on $\ldop_\sigma$. Let $\vol(\Lambda(\adop_\sigma)^{L_1}, \dots, \Lambda(\adop_\sigma)^{L_d})$  be the mixed volume of the corresponding fundamental domains in $\ldop_\sigma^*$ (see \ref{mixed volume} for  definition and properties). Since the mixed volume is positive, the following formula implies Theorem \ref{Theorem 3}:
\end{art}

\begin{thm} \label{explicit expression}
For a measurable subset $\Omegabar$ of the atom $\sigmabar$ as in  \ref{setting}, we have 
$$\mu(\Omegabar) = \frac{d!}{[X':\Xan]} \sum_{j \in J(\sigmabar)}
[\zdop^d : \hat{\ell}_{j}^{(0)}(\Lambda(\adop_\sigma)^*)] \cdot 
\frac{\vol(\Lambda(\adop_\sigma)^{L_1}, \dots, \Lambda(\adop_\sigma)^{L_d}) \cdot\vol(\Omegabar)}{\vol(\Lambda(\adop_\sigma)^*)\cdot\vol(\Lambda(\adop_\sigma))},$$
 \end{thm}



\proof To prove the theorem, we may assume that $\Omegabar$ is a $d$-dimensional polytope contained in $\sigmabar$ (using monotone convergence). Note that the odd part of $L_j$ does not influence the bilinear form $b_j$ and hence we may suppose, using \ref{can measure}, that every $L_j$ is symmetric. 
By multilinearity of $\mu$ and the mixed volume, we may assume $\overline{L}:=\overline{L}_1= \cdots = \overline{L}_d$. In the next paragraph, we will fix a Mumford model $\Acal$ of $A$ associated to a ``generic'' choice of a polytopal decomposition $\Ccalbar$. It will be crucial for the proof that this choice is as generic as possible. In particular, $\Ccalbar$ is only $\Gamma'$-rational for an infinite dimensional $\qdop$-subspace $\Gamma'$ of $\rdop$ equal to the value group of a complete algebraically closed extension $\kdop'$ of $\kdop$  and hence $\Acal$ is only defined over the valuation ring of $\kdop'$. By Remark \ref{base change} and since $\valbar(\Xan)$ is invariant under base change,  we are allowed to make the analytic calculations for $\mu$ over $\kdop'$.

Let $\tilde{P}_1, \dots, \tilde{P}_R$ be the zero-dimensional strata of $\tilde{\Xcal}'$. As in \ref{list of singularities}, they induce affine maps $f_{j, \rm aff}: \Delta(d,\pi_j) \rightarrow \rdop^n$ and we denote the image simplex by $\rho_j$. According to \ref{setting}, we assume that $\rho_1,\dots,\rho_N$ are $d$-dimensional and that $\rho_{N+1},\dots,\rho_R$ are lower dimensional. Let $\Sigma$ be the collection of these simplices $\rho_j$ in $\rdop^n$ together with all their faces and their $\Lambda$-translates. 
We choose $\Ccalbar$, $\Acal$ and  $\Lcal$ as in Lemma \ref{irrational construction}.  
By multilinearity, we may assume that $\Lcal$ is a $(\kdop')^\circ$-model of $L$.

We will proof Theorem \ref{explicit expression} now in four steps. We give first an {\it outline} of the plan:  

We fix a rigidification on $L$ such that the associated canonical metric $\metr_{ {\rm can}}$ is the metric of $\overline{L}$ (see \ref{canonical metrics}). The rigidification and the theorem of the cube yield an identification $[m]^*L=L^{\otimes m^2}$ for $m \in \zdop$. Let $\metr$ be the formal metric associated to $\Lcal$ (see \ref{formal metrics}). A variant of Tate's limit argument (\cite{BG}, proof of Theorem 9.5.4) says
\begin{equation} \label{Tate's limit argument}
\metr_{ {\rm can}} = \lim_{m \to \infty} \left( [m]^* \metr \right)^{1/m^2}.
\end{equation}
Let  $\Acal_m$ be the Mumford model of $A$ associated to $\overline{\Ccal_m}:=\overline{\frac{1}{m}\Ccal}$ and let $\Xcal_m$ be the closure of $X$ in $\Acal_m$. Using the very definition of the measure $\mu$ and \eqref{Tate's limit argument}, we will show in a first step that $\mu$ is a weak limit of a sum of Dirac measures $\valbar_*(\delta_Z)$, where $Z$ ranges over the irreducible components of $(\Xcal_m^{\rm f-an})\sptilde$ and where  $\xi_Z$ is the 
point of the Berkovich space $\Xan$ corresponding to $Z$ (see \ref{formal analytic varieties}). In a second step, we will replace the $Z$ by the irreducible components $Y$ of the special fibre $\tilde{\Xcal}_m$. Since the reduction $(\Xcal_m^{\rm f-an})\sptilde$ is a finite covering of $\tilde{\Xcal}_m$, this will be a consequence of the projection formula.  
 
After the first two steps we have $\mu(\Omegabar)=\lim_{m \to \infty} a_m$, where  $a_m$ depends on the multiplicities and the degrees of the $Y$'s. In the third step, we transform the limit into a multiple of $\vol(\Omegabar)$. To make this plausible, note that the multiplication by $m$ on $A$ extends uniquely to a morphism $\varphi_m:\Acal_m \rightarrow \Acal_1=\Acal$ (see \ref{multiplication with m}). Applying projection formula to $\tilde{\varphi}_m$, we will relate the degree of $Y$ to $\deg_{\tilde{\Lcal}}(Y_\ub)$ for a vertex $\ub$ of $\Ccal \cap (m \adop_\sigma)$ with $\ubb \in \Omegabar$. Here, $Y_\ub:=Y_{\Delta(\ub)}$ is the toric variety in $\tilde{\Acal}$ associated to the unique $\Delta(\ub) \in \Ccal$ with $\ub \in \relint(\Delta(\ub))$ (see Remark \ref{polytopal toric variety}). We will use the alteration  $f:X' \rightarrow \Xan$ from \ref{setting} and Section 7 to express the multiplicity of $Y$ in terms of indices of lattices. For $m \to \infty$, we will get $\mu(\Omegabar)=s \cdot\vol(\Omegabar)$, where $s$ is a linear combination of $\deg_{\tilde{\Lcal}}(Y_\ub)$ with $\ub$ ranging over the vertices of $\Ccal \cap \ldop_\sigma$ modulo the stabilizer $\Lambda(\adop_\sigma)$ from \ref{setting}. In the theory of toric varieties, there is a formula for $\deg_{\tilde{\Lcal}}(Y_\ub)$ as the volume of a polytope associated to the fan of $Y_\ub$. We apply this in the fourth step to calculate $s$ in terms  of the dual complex of $\Ccal \cap \ldop_\sigma$ using  the appendix on convex geometry.  This will prove Theorem \ref{explicit expression}.

{\it Step 1:  $\mu$ is a weak limit of discrete measures related to the irreducible components of $(\Xcal_m^{\rm f-an})\sptilde$.}

For $m \geq 1$, we have seen in the outline that a unique morphism $\varphi_m:\Acal_m \rightarrow \Acal_1=\Acal$ extends multiplication by $m$ on $A$. Recall that $\metr$ is the formal metric associated to the $(\kdop')^\circ$-model $\Lcal$ of $L$ on $\Acal$. Clearly, $[m]^*\metr$ is the formal metric associated to $\varphi_m^*(\Lcal)$.  The composition of the canonical finite morphism $\tilde{\iota}:(\Xcal_m^{\rm f-an})\sptilde \rightarrow \tilde{\Xcal}_m$ from \ref{admissible formal schemes} with $\tilde{\varphi}_m$ is denoted by $\tilde{\phi}_m$. By \eqref{Tate's limit argument}, \ref{Dirac decomposition} and \ref{Chern continuity}, we have the following weak limit of regular Borel measures on $\valbar(\Xan)$:
\begin{equation} \label{toric measure 2}
\mu = \lim_{m \to \infty} m^{-2d} \sum_Z \deg_{\tilde{\phi}_m^*(\tilde{\Lcal})}(Z) \valbar_*(\delta_{\xi_Z}),
\end{equation}
where $Z$ ranges over all irreducible components of $(\Xcal_m^{\rm f-an})\sptilde$. For our polytope $\Omegabar \subset \sigmabar$, this yields 
\begin{equation} \label{toric measure 3}
\mu(\Omegabar)=\lim_{m \to \infty} m^{-2d} \sum_Z \deg_{\tilde{\phi}_m^*(\tilde{\Lcal})}(Z), 
\end{equation}
where $Z$ ranges over all irreducible components of $(\Xcal_m^{\rm f-an})\sptilde$ with $\valbar(\xi_Z) \in \Omegabar$. 

{\it Step 2: We replace the $Z$'s in \eqref{toric measure 3} by the irreducible components $Y$ of $\tilde{\Xcal}_m$.}

Since the morphism $\tilde{\iota}$ is finite and surjective, the set of irreducible components of $(\Xcal_m^{\rm f-an})\sptilde$ is mapped onto the set of irreducible components of $\tilde{\Xcal}_m$. By our choice of $\Ccalbar$ from Lemma \ref{irrational construction}, we note that  $\overline{\Ccal_m}$ is transversal to $\valbar(\Xan)$ and hence every irreducible component $Y$ of $\tilde{\Xcal}_m$ corresponds to an equivalence class $\overline{\Delta} \cap \valbar(\Xan)$ of transversal vertices of $\overline{\Ccal_m} \cap \valbar(\Xan)$ for a unique $d$-codimensional polytope $\Deltabar \in \overline{\Ccal_m}$ (see Theorem \ref{periodic transverse}). Since $\overline{\Ccal_m}$ is $\Sigmabar$-generic, it is clear that $\Deltabar$ intersects the $d$-dimensional atom $\sigmabar$ in at most one point. If a $Z$ from \eqref{toric measure 3} is lying over $Y$, then the proof of Theorem \ref{transversal correspondence} shows that $\valbar(\xi_Z) \in \Deltabar \cap {\valbar(\Xan)}$, i.e. we have a transversal vertex corresponding to $Y$ which is contained in $\Omegabar$. We say that $Y$ is {\it $\Omegabar$-inner} if
$$\Deltabar \cap \valbar(\Xan) = \Deltabar \cap \Omegabar.$$
If $Y$ is not $\Omegabar$-inner, then the corresponding $\Deltabar$ intersects also an atom $\sigmabar' \neq \sigmabar$. This $\Deltabar$ has to be a face of an $n$-dimensional polytope in $\Ccalbar$ intersecting the boundary of $\sigmabar$. By an easy argument covering the boundary, we conclude that the number of such $Y$ is of order $O(m^{d-1})$. We will use this later to show that the $Z$ lying over non-$\Omegabar$-inner $Y$ may be neglected in \eqref{toric measure 3}.

We consider now an $\Omegabar$-inner $Y$. Since $\Omegabar$ is a polytope contained in $\sigmabar$, we conclude that $\Deltabar \cap \Omegabar$ is just a point $\ubb$. For an irreducible component $Z$ of $(\Xcal_m^{\rm f-an})\sptilde$ lying over $Y$, we have seen above that $\valbar(\xi_Z) \in \Deltabar \cap {\valbar(\Xan)}$ and hence  $\valbar(\xi_Z) = \ubb \in \Omegabar$.
Since $(\Xcal_m^{\rm f-an})\sptilde$ is reduced, we have
$$\sum_Z \deg_{\tilde{\phi}_m^*(\tilde{\Lcal})}(Z)
=  \deg_{\tilde{\phi}_m^*(\tilde{\Lcal})}((\Xcal_m^{\rm f-an})\sptilde),$$
where  $Z$ ranges over all irreducible components of $(\Xcal_m^{\rm f-an})\sptilde$. Since $\Xcal_m^{\rm f-an}$ and $\Xcal_m$ have the same generic fibre, the projection formula (\cite{Gu2}, Proposition 4.5) shows that  $\tilde{\iota}_*((\Xcal_m^{\rm f-an})\sptilde)$ is equal to the cycle associated to $\tilde{\Xcal}_m$. For the multiplicity $m(Y, \tilde{\Xcal}_m)$, projection formula yields
\begin{equation} \label{toric measure 4}
\sum_{\tilde{\iota}(Z)=Y}  \deg_{\tilde{\phi}_m^*(\tilde{\Lcal})}(Z) 
= m(Y, \tilde{\Xcal}_m) \deg_{\tilde{\varphi}_m^*(\tilde{\Lcal})}(Y). 
\end{equation}
This is true for any irreducible component $Y$ of $\tilde{\Xcal}_m$. First, we apply \eqref{toric measure 4} on the right hand side of \eqref{toric measure 3} to transform the contribution of all $Z$ lying over an $\Omegabar$-inner $Y$ into the sum 
\begin{equation} \label{toric measure 5}
S_m:=m^{-2d} \sum_Y m(Y, \tilde{\Xcal}_m) \deg_{\tilde{\varphi}_m^*(\tilde{\Lcal})}(Y),
\end{equation}
where $Y$ ranges over all $\Omegabar$-inner irreducible components of $\tilde{\Xcal}_m$. If $Y$ is an irreducible component of $\tilde{\Xcal}_m$ which is not $\Omegabar$-inner, then it is possible that there are irreducible components $Z$, $Z'$ of $(\Xcal_m^{\rm f-an})\sptilde$ lying over $Y$ with $\valbar(\xi_Z) \in \Omegabar$ and $\valbar(\xi_{Z'}) \not \in \Omegabar$.  Then
\eqref{toric measure 5} is just an upper bound for the contribution of those $Z$ in \eqref{toric measure 3} lying over such a $Y$.

{\it Step 3: Transformation of the limit in \eqref{toric measure 3} into a multiple of $\vol(\Omegabar)$.}

Let $Y$ be an irreducible component of $\tilde{\Xcal}_m$. By transversality of $\overline{\Ccal_m}$ and Theorem \ref{periodic transverse}, there is a unique $d$-codimensional polytope $\Deltabar \in \overline{\Ccal_m}$ such that $Y$ is the toric variety given as the closure of the torus orbit associated to $\relint(\Deltabar)$. 

Now we assume that $Y$ is $\Omegabar$-inner. As we have seen in step 2, $\Deltabar \cap {\valbar(\Xan)}$ is a transversal vertex $\ubb$ of $\overline{\Ccal_m} \cap \valbar(\Xan)$. Since $\ubb \in \Omegabar$, there is a unique lift $\ub$ to the affine space $\adop_\sigma$ from \ref{setting}. We conclude that $m \ub$ is a vertex of $\Ccal \cap (m \adop_\sigma)$. 

Every $\ub' \in \rdop^n$ is contained in  the relative interior of a unique $\Delta(\ub') \in \Ccal$.  Let $Y_{\ub'}:= Y_{\Delta(\ub')}$ be the associated toric variety in $\tilde{\Acal}$ (see Remark \ref{polytopal toric variety} and Proposition \ref{torus orbit 2}). Note that $Y_{\ub' +\lambda}=Y_{\ub'}$ for $\lambda \in \Lambda$. In the situation above, we have $m \Delta = \Delta(m \ub)$ and $\tilde{\varphi}_m(Y) = Y_{m \ub}$. Applying projection formula to the morphism $\tilde{\varphi}_m: \tilde{\Acal}_m \rightarrow \tilde{\Acal}$ and using Proposition \ref{multiplication with m}, we get
\begin{equation} \label{toric measure 5'}
\deg_{\tilde{\varphi}_m^*(\tilde{\Lcal})}(Y) = m^d \deg_{\tilde{\Lcal}}(Y_{m \ub}).
\end{equation}
Now we use the alteration $f:X' \rightarrow \Xan$ and the associated affine maps 
$f_{j, \rm aff}: \Delta(d,\pi_j) \rightarrow \rdop^n$  with $d$-dimensional images $\overline{\rho}_j$, $j=1,\dots, N$, from \ref{setting}. 
Since $\Ccalbar$ is $\Sigmabar$-generic, the final remark in \ref{index definition} shows that $\ind(\overline{ \Delta'},\overline{f}_{j, {\rm aff}})$ is well-defined for all $d$-codimensional $\overline{\Delta'} \in \Ccalbar$ and depends only on the linear space $\ldop_{\Delta'}$ for $j=1, \dots, N$. We consider the multiplicity
$$\vartheta(\overline{ \Delta'}, \sigmabar):= \frac{1}{[X':\Xan]} \sum_{j \in J(\sigmabar)} \ind(\overline{ \Delta'}, \overline{f}_{j,{\rm aff}})$$
of $\overline{ \Delta'}$ relative to the atom $\sigmabar$, where $J(\sigmabar):=\{j=1,\dots,N \mid \sigmabar \subset \overline{\rho}_j\}$.  Proposition \ref{multiplicity formula} yields
\begin{equation} \label{multiplicity and m}
m(Y,\tilde{\Xcal}_m)= \frac{1}{[X':\Xan]} \sum_{j \in J} \ind(\Deltabar,\overline{f}_{j,\rm aff}),
\end{equation}
where $J:=\{j=1, \dots, N \mid \Deltabar \cap \overline{\rho}_j \neq \emptyset\}$. Since $\overline{\Ccal_m}$ is $\Sigmabar$-transversal, we note that $\ubb$ is an interior point of the atom $\sigmabar$. By $\Deltabar \cap {\valbar(\Xan)} = \{\ubb\}$, we deduce that $j \in J$ if and only if $\Deltabar \cap   \overline{\rho}_j  = \{\ubb\}$ and this is equivalent to $\sigmabar \subset  \overline{\rho}_j $. This means $J=J(\sigmabar)$. We have $\ldop_{m \Delta}= \ldop_\Delta$ and 
hence  $\ind(\overline{m \Delta}, \overline{f}_{j,\rm aff})= \ind(\Deltabar,\overline{f}_{j, \rm aff})$. Now 
\eqref{multiplicity and m} yields
\begin{equation} \label{toric measure 5''}
\ind(\overline{m \Delta}, \overline{f}_{j,\rm aff})= \ind(\Deltabar,\overline{f}_{j, \rm aff}), 
\quad m(Y, \tilde{\Xcal}_m) = \vartheta(\Deltabar(m \ub), \sigmabar).
\end{equation}
In \eqref{toric measure 5'} and \eqref{toric measure 5''}, we have related the degree and the multiplicity of $Y$ (or equivalently an $\Omegabar$-inner transversal vertex of $\overline{\Ccal_m} \cap {\valbar(\Xan)}$ as we have seen in step 2) to the vertex $m\ub$ of $\Ccal \cap (m \adop_\sigma)$.  Both formulas are $\Lambda$-periodic in $m\ub$. Note that the atom $\sigmabar$ of $\Sigmabar$ is rational and hence there is $m_0 \in \ndop$ such that $m_0 \adop_\sigma = \lambda_0 + \ldop_\sigma$ for some $\lambda_0 \in \Lambda$. Up to now, we assume that $m \in \ndop m_0$. By periodicity, we may express \eqref{toric measure 5'} and \eqref{toric measure 5''} in terms of the vertex $\ub':= m \ub - \frac{m}{m_0}\lambda_0$ of $\Ccal \cap \ldop_\sigma$. Note that $\Ccal \cap \ldop_\sigma$ is  
$\Lambda(\adop_\sigma)$-periodic with fundamental domain denoted by $F$. For $m$ sufficiently large, the number of $\frac{1}{m}\Lambda(\adop_\sigma)$-translates of $\frac{1}{m}F$ contained in the lift $\Omega$ of $\Omegabar$ to $\adop_\sigma$ is $m^d\vol(\Omega)/\vol(F) + O(m^{d-1})$. By \eqref{toric measure 5'} and \eqref{toric measure 5''} inserted in \eqref{toric measure 5}, we get
\begin{equation} \label{toric measure 6}
S_m=\frac{\vol(\Omegabar)}{\vol(F)} \cdot  \sum_{\ub'}\vartheta(\Deltabar(\ub'), \sigmabar) \deg_{\tilde{\Lcal}}(Y_{\ub'}) +O(\frac{1}{m}),
\end{equation}
where $\ub'$ ranges over the  vertices of $\Ccal \cap \ldop_\sigma$ contained in $F$. We claim that
\begin{equation} \label{toric measure 1}
\mu(\Omegabar) = \frac{\vol(\Omegabar)}{\vol(\Lambda(\adop_\sigma))} \cdot \sum_{\ub'} \vartheta(\Deltabar(\ub'), \overline{\sigma}) \deg_{\tilde{\Lcal}}(Y_{\ub'}),
\end{equation}
where $\ub'$ ranges over all vertices of $\Ccal \cap \ldop_\sigma$ modulo $\Lambda(\adop_\sigma)$. By \eqref{toric measure 3}, \eqref{toric measure 5} and \eqref{toric measure 6},  it remains to show that the $Z$ in \eqref{toric measure 3} lying over non-$\Omegabar$-inner $Y$ may be neglected. We have seen in step 2 that the number of such $Y$ is $O(m^{d-1})$. Formula \eqref{multiplicity and m} holds also for non-$\Omegabar$-inner $Y$ and proves
$$m(Y,\tilde{\Xcal}_m) \leq \frac{1}{[X':\Xan]} \sum_{j=1}^N \ind(\Deltabar, \overline{f}_{j, \rm aff}),$$
where $\Deltabar$ is the $d$-codimensional polytope in $\overline{\frac{1}{m}\Ccal}$ corresponding to $Y$. We have
$$\ind(\Deltabar,\overline{f}_{j, \rm aff}) \leq  \max_{\overline{ \Delta'}} \left(\ind(\overline{ \Delta'},\overline{f}_{j, \rm aff})\right),$$
where $\overline{ \Delta'}$ is ranging over all $d$-codimensional simplices in $\Ccalbar$. This leads to a bound of $m(Y,\tilde{\Xcal}_m)$ independent of $Y$ and $m$. There is also such a bound for  $m^{-d}\deg_{\tilde{\varphi}_m^*(\tilde{\Lcal})}(Y)$. Indeed, this follows from projection formula as in \eqref{toric measure 5'} replacing $Y_{m \ub}$ by  the closure of a $d$-dimensional torus orbit of $\tilde{\Acal}$. Using  the final remark in step 2, the contribution of the $Z$ in  \eqref{toric measure 3} lying over non-$\Omegabar$-inner $Y$ is bounded by $O(\frac{1}{m})$ and therefore may be neglected in \eqref{toric measure 3}. This proves  \eqref{toric measure 1}.

{\it Step 4: Calculation of the sum in \eqref{toric measure 1} in terms of the dual complex of $\Ccal \cap \ldop_\sigma$.}

We have chosen $\Ccalbar$ and the $(\kdop')^\circ$-model $\Lcal$ of $L$ from Lemma \ref{irrational construction}. Recall from its proof that $\Lcal$ was constructed by a strongly polyhedral convex function $f_\Lcal$ with respect to $\Ccal$ using Proposition \ref{Mumford's line bundles} and Corollary \ref{ample reduction}. In particular, $f_\Lcal(\ub)=\mb_\Delta \cdot \ub + c_\Delta$ on $\Delta \in \Ccal$ with $m_\Delta \in \zdop^n$. Let $\ub'$ be a vertex of $\Ccal \cap \ldop_\sigma$ and let $\Delta \in \Star_n(\Delta(\ub'))$ with peg $\mb_\Delta$ (see \ref{strongly polyhedral function} and \ref{dual complex}). The theory of toric varieties (see \cite{Fu2}, 5.3) and \ref{alternative description of dual complex} show that
\begin{equation} \label{toric measure 7}
\deg_{\tilde{\Lcal}}(Y_{\ub'}) = d! \cdot \vol\left(\Delta({\ub'})^{f_\Lcal}  \right) \cdot \vol\left(\zdop^n \cap \Delta({\ub'})^\bot \right)^{-1},
\end{equation} 
where $\Delta({\ub'})^\bot$ denotes the orthogonal complement of $\Delta({\ub'})$.  Note that $\Delta({\ub'})^{f_\Lcal} $ is a $d$-di\-men\-sional polytope of the dual complex $\Ccal^{f_\Lcal}$  contained in  $\mb_\Delta + \Delta({\ub'})^\bot$. We may also consider the complex $\Ccal \cap \ldop_\sigma$ in $\ldop_\sigma$. Clearly, $f_\Lcal$ restricts to a strongly polyhedral convex function $g$ on $\ldop_\sigma$ and we may form the dual complex $(\Ccal \cap \ldop_\sigma)^g$ in $\ldop_\sigma^*$ which  is $\Lambda(\adop_\sigma)^L$-periodic and hence covers $\ldop_\sigma^*$ (use \cite{McM}, Theorem 3.1). Let $P$ be the dual map of $\ldop_\sigma \hookrightarrow \rdop^n$. Then $\{P(\mb_{\tau}) \mid \tau \in \Ccal, \, \dim(\tau)=n, \, \tau \cap \ldop_\sigma \neq \emptyset\}$ are the pegs of $\Ccal \cap \ldop_\sigma$. Since $\Ccal$ is transversal to $\ldop_\sigma$, we easily deduce that 
$$P\left(\Delta({\ub'})^{f_\Lcal} \right) =\{{\ub'}\}^g$$
and hence \eqref{toric measure 7} yields
\begin{equation} \label{toric measure 8}
\deg_{\tilde{\Lcal}}(Y_{\ub'}) = d! \cdot \vol(\{{\ub'}\}^g) 
\cdot \vol\left(P\left(\zdop^n \cap \Delta({\ub'})^\bot \right)\right)^{-1}.
\end{equation}
Duality of lattices and the definitions of the multiplicity $\vartheta(\Deltabar,\sigmabar)$ and the index in \ref{index definition} imply
\begin{align*} 
\vartheta(\Deltabar({\ub'}), \sigmabar) & = 
\frac{1}{[X':\Xan]} \sum_{j \in J(\sigmabar)}
[\zdop^d : \hat{\ell}_{j}^{(0)}(P(\zdop^n \cap \Delta({\ub'})^\bot ))]\\
  & =\frac{1}{[X':\Xan]} \sum_{j \in J(\sigmabar)}
[\zdop^d : \hat{\ell}_{j}^{(0)}(\Lambda(\adop_\sigma)^*)] \cdot  
\frac{\vol(P(\zdop^n \cap \Delta({\ub'})^\bot ))}{\vol(\Lambda(\adop_\sigma)^*)}.
\end{align*}
Using this and \eqref{toric measure 8}  in \eqref{toric measure 1}, we get
$$\mu(\Omegabar)=\frac{d!}{[X':\Xan]}  \sum_{\ub'} \vol(\{{\ub'}\}^g)
\sum_{j \in J(\sigmabar)}[\zdop^d : \hat{\ell}_{j}^{(0)}(\Lambda(\adop_\sigma)^*)]  \cdot 
\frac{\vol(\Omegabar)}{\vol(\Lambda(\adop_\sigma)^*)\vol(\Lambda(\adop_\sigma))},$$
where ${\ub'}$ ranges over the vertices of $\Ccal \cap \ldop_\sigma$ modulo $\Lambda(\adop_\sigma)$. The theorem follows now from the $\Lambda(\adop_\sigma)^L$-periodicity of the covering $(\Ccal \cap \ldop_\sigma)^g$ of $\ldop_\sigma^*$.
\qed

\begin{rem} \rm \label{arbitrary line bundles}
For arbitrary line bundles $\overline{L}_1, \dots , \overline{L}_d$ on $A$ endowed with canonical metrics, Theorem \ref{Theorem 3} yields that $\mu$ is a piecewise Haar measure on $\valbar(\Xan)$.  Indeed, every  line bundle on $A$ is isomorphic to the ``difference'' of two ample  line bundles and multilinearity yields the claim.
\end{rem}

\section{Generalizations}

First, we relate tropical varieties to the skeletion of a strictly semistable $\kdop^\circ$-model. We use it to prove Theorem \ref{Theorem 2}. Then we describe the canonical measures on a closed subvariety of a totally degenerate abelian variety generalizing Theorem \ref{explicit expression}. 

\begin{art} \label{skeleton} \rm 
Berkovich has shown in \cite{Ber5}, \S 5, that a strongly non-degenerate polystable fibration over $\kdop^\circ$ has a canonical polytopal subset called the skeleton. For simplicity, we restrict its description to the case of a strictly semistable formal scheme $\Xcal'$ over $\kdop^\circ$:

By Proposition \ref{semistable properties}, every irreducible component $Y$ of $\tilde{\Xcal}'$ induces a Cartier divisor $D_Y$ on $\Xcal'$ with $\cyc(D_Y)=v(\pi)Y$ just by lifting the local equations $\tilde{\gamma}$ for $Y$ in $\tilde{\Xcal}'$, where $\pi \in K^{\circ\circ}$ is from \ref{semistable}. If $X':=(\Xcal')^{\rm an}$ is connected, then we may choose $\pi$ independent of $Y$. For simplicity of notation, this will be assumed in the following.

By \ref{formal metrics} and \ref{pseudo-divisors},  there is a global section $s_Y$ of $O_{X'}$ and a formal metric $\metr_Y$ on $O_{X'}$ such that $D_Y=\widehat{\Div}(s_Y)$. Let $I$ be the set of irreducible components of $\tilde{\Xcal}'$. Then we get an analogue of tropical geometry by considering
$$\val:X' \longrightarrow \rdop^I, \quad x' \mapsto \left( - \log ||s_Y(x')|| \right)_{Y\in I}.$$
We have seen in Proposition \ref{affinization and strata} that every stratum $S$ of $\tilde{\Xcal}'$ gives rise to a canonical simplex
$$\Delta_S:= \{ \ub \in \rdop^I \mid \text{$u_Y = 0$ if $S \cap Y = \emptyset$}, \, \sum_{Y \cap S \neq \emptyset} u_Y = v(\pi) \}.$$
As in the proof of Proposition \ref{polyhedron}, we deduce that $\val(X')$ is covered by $\Delta_S$, where $S$ is ranging over the set ${\rm str}(\tilde{\Xcal}')$ of strata. The strata of dimension $d:=\dim(X')$ are in one-to-one correspondence with $I$ and hence with the vertices $(0, \dots, v(\pi), \dots, 0)$, where just the $Y$-th entry is non-zero. Note however that different lower dimensional strata may induce the same canonical simplex. To omit this, we define the  abstract metrized polytopal set 
$${\mathbb D}(\Xcal'):= \left(\coprod_{S \in {\rm str}(\tilde{\Xcal}')} (\Delta_S,S) \right)/\sim.$$
Here, $(\Delta_S,S)$ is isometrically identified with $\Delta_S$ and $(\ub_1,S_1) \sim (\ub_2,S_2)$ if there is $S_3 \in {\rm str}(\tilde{\Xcal}')$ such that the closure of $S_3$ contains $S_1 \cup S_2$ and if $\ub_1=\ub_2 \in \Delta_{S_3}$. The set ${\mathbb D}(\Xcal')$ reflects  the incidence of strata closures as there is a bijective correspondence between ${\rm str}(\tilde{\Xcal}')$ and the simplices $(\Delta_S,S)$ of ${\mathbb D}(\Xcal')$. Note that we may lift $\val$ to a continuous surjective map 
$$\Val: X' \longrightarrow {\mathbb D}(\Xcal'), \quad x' \mapsto (\val(x'),S(\tilde{x}')),$$
where $S(\tilde{x}')$ is the stratum containing $\tilde{x}'$. 

Berkovich introduces a partial ordering $\preceq$ on $X'$ (depending on $\Xcal'$) by $x' \preceq y'$ if, for every affinoid algebra $\Acal$ with  \'etale morphism $\varphi:\Xcal = \Spf(\Acal^\circ) \rightarrow \Xcal'$ and $x \in (\varphi^{\rm an})^{-1}(x')$, there is $y \in (\varphi^{\rm an})^{-1}(y')$ with $|f(x)| \leq |f(y)|$ for all $f \in \Acal$. The set of maximal points with respect to this ordering is called the {\it skeleton} $S(\Xcal')$. By \cite{Ber5}, Theorem 5.1.1, the map $\Val$ restricts to a homeomorphism of $S(\Xcal')$ onto ${\mathbb D}(\Xcal')$. We use it to endow $S(\Xcal')$ with the structure of a metrized polytopal set, i.e. we identify $S(\Xcal')$ with ${\mathbb D}(\Xcal')$. 
\end{art}

\begin{prop} \label{affinization on skeleton}
Let $\Xcal'$ be a strictly semistable formal scheme over $\kdop^\circ$ with generic fibre $X'$ and let $A$ be a totally degenerate abelian variety over $\kdop$, i.e. $A^{\rm an} = (\Tor)_\kdop^{\rm an}/M$ for a lattice $M$. For a morphism $f:X' \rightarrow A$ and $\Lambda := \val(M)$, there is a unique continuous map $\overline{f}_{\rm aff}: S(\Xcal') \rightarrow \rtor$ with
\begin{equation} \label{Val and aff}
\overline{f}_{\rm aff} \circ \Val = \valbar \circ f.
\end{equation}
Moreover, $f_{\rm aff}$ is an affine map on every simplex $(\Delta_S,S)$, $S \in {\rm str}(\tilde{\Xcal}')$. 
\end{prop}

\proof By \ref{affine map} and Proposition \ref{affinization and strata}, we get affine maps satisfying \eqref{Val and aff} on every simplex $(\Delta_S,S)$ of $S(\Xcal')$. They fit to define $\overline{f}_{\rm aff}$ on $S(\Xcal')$. Uniqueness follows from surjectivity of $\Val$. \qed 

\vspace{3mm}  
{\bf Proof of Theorem \ref{Theorem 2}: \/} 
We may assume that $X'$ is irreducible. By Chevalley's theorem, $f(X')$ contains an open dense subset of a $d$-dimensional closed subvariety $X$ of $A$. Since $\valbar$ is continuous, we conclude that
$$\valbar(X^{\rm an}) = \valbar(f(X')^{\rm an}).$$
By Proposition \ref{affinization on skeleton}, we have $\valbar(X^{\rm an})=\overline{f}_{\rm aff}(S(\Xcal'))$. The tropical variety $\valbar(X^{\rm an})$ is $d$-dimensional (Theorem \ref{periodic main theorem}) and hence $S(\Xcal')$ contains a simplex $(\Delta_S,S)$ of dimension at least $d$. The vertices of $\Delta_S$ correspond to irreducible components of $\tilde{\Xcal}'$ containing $S$. \qed

\begin{rem} \rm \label{analytic generalization of dimension theorem}
Theorem \ref{Theorem 2} is also true analytically in the style of Theorem \ref{toric dimension bound} with $A$ replacing $\Tor$. It wouldn't be difficult to deduce Theorem \ref{Theorem 2} directly from Theorem \ref{toric dimension bound}. 
\end{rem}

\begin{thm} \label{generalization of dimension theorem for arbitrary abelian var}
If $A$ is an arbitrary abelian variety over $\kdop$ with $t$-dimensional formal abelian scheme $\Bcal$ in the Raynaud extension $E$ of $A$ (see \cite{BL2}, \S 1), then Theorem \ref{Theorem 2} holds with $1-t+\dim f(X')$ replacing $1+\dim f(X')$. 
\end{thm}

\proof We define $\val$ on $E$ as in \cite{BL2}, p. \!656. Note that $E$ is locally trivial over $\Bcal^{\rm an}$ and hence we deduce easily form Theorem \ref{periodic main theorem} that $\valbar(X)$ is at least of dimension $\dim(X)-t$. Since Proposition \ref{affinization on skeleton} generalizes to this context, we can follow the proof of Theorem \ref{Theorem 2} to get the claim.  \qed

\begin{art} \rm \label{new setting}
For the remaining part of this section, we consider a field $K$ with a discrete valuation $v$ and we assume that $\kdop$ is the completion of the algebraic closure of the completion of $K$ with respect to $v$. The unique extension of $v$ to a valuation of $\kdop$ is also denoted by $v$ with corresponding absolute value $|\phantom{a}| := e^{-v}$. 

Let $A$ be an abelian variety over $K$ which is totally degenerate over $\kdop$, i.e. $A_\kdop^{\rm an} = (\Tor)_\kdop^{\rm an}/M$ for a lattice $M$. Again, let $\Lambda := \val(M)$ be the corresponding complete lattice in $\rdop^n$. Let $X$ be a closed geometrically integral $d$-dimensional closed subvariety of $A$.  

For canonically metrized ample line bundles $\overline{L}_1,\dots, \overline{L}_d$ on $A$, we want to describe the canonical measure $c_1(\overline{L}_1|_X) \wedge \dots \wedge c_1(\overline{L}_d|_X)$ on $X^{\rm an}$. By de Jong's alteration theorem (\cite{dJ}, Theorem 6.5), there is an alteration $f:X' \rightarrow X_\kdop^{\rm an}$ and a strictly semistable formal $\kdop^\circ$-model $\Xcal'$ of $X'$.

A $d$-dimensional simplex $(\Delta_S,S)$ of $S(\Xcal')$, $S \in {\rm str}(\tilde{\Xcal}')$, is called {\it non-degenerate} with respect to $f$ if $f_{\rm aff}(\Delta_S)$ is also $d$-dimensional. Then $S$ is a $\tilde{\kdop}$-rational point contained in $d+1$-irreducible components $Y_0, \dots ,Y_d$ of $\tilde{\Xcal}'$. If $u_0, \dots ,u_d$ denote the corresponding coordinates, then $\Delta_S$ is given by $\{u_0+ \dots + u_d = v(\pi)\}$. We consider the $d$-dimensional standard simplex 
$$\Sigma_S:= \{\ub \in \rdop_+^d \mid u_1 + \dots + u_d \leq v(\pi)\}$$
and the affine map $f_{\rm aff}^{(0)}:\Sigma_S \rightarrow \rdop^n$ given by 
$$ f_{\rm aff}^{(0)}(u_1, \dots ,u_d)= f_{\rm aff}(u_0, \dots, u_d), \quad (u_0, \dots, u_d) \in \Delta_S.$$
If we extend $f_{\rm aff}^{(0)}- f_{\rm aff}^{(0)}({\mathbf 0})$, then we get an associated injective linear map
$\ell_S^{(0)}:\rdop^d \rightarrow \rdop^n$ as in \ref{index definition}. It is defined over $\zdop$ and hence $\Lambda_S:=(\ell_S^{(0)})^{-1}(\Lambda)$ is a complete rational lattice in $\rdop^d$. The positive definite bilinear form $b_j$ associated to $L_j$ induces a complete lattice
$$\Lambda_S^{L_j} := \{b_j(\ell_S^{(0)}(\cdot),\lambda) \mid \lambda \in \Lambda\}$$
on $(\rdop^d)^*=\rdop^d$. We denote by $\vol$ the (mixed) volume with respect to the  Lebesgue measure on $\rdop^d$. By Corollary \ref{Chern properties}, the following result  describes the canonical measure on $\Xan$:
\end{art}

\begin{thm} \label{explicit generalization}
Under the hypothesis of \ref{new setting}, the measure $\mu := c_1(f^*(\overline{L}_1)) \wedge \dots \wedge c_1(f^*(\overline{L}_d))$ is supported in the union of the non-degenerated simplices $(\Delta_S,S)$ of $S(\Xcal')$. For a measurable subset $\Omega$ of such a simplex, we have
$$\mu(\Omega) = d! \cdot \frac{\vol(\Lambda_S^{L_1}, \dots, \Lambda_S^{L_d})}{\vol(\Lambda_S)} \cdot \vol(\Omega).$$
\end{thm}

\proof We follow the steps of the proof of Theorem \ref{explicit expression}. We may assume that $\overline{L}:=\overline{L}_1= \dots = \overline{L}_d$ and $L$ symmetric. We use the same $\Sigmabar$, $\Ccalbar$, $\Acal$, $\Lcal$, $\Acal_m$ and $\Xcal_m$. The main difference is that the role of $\Xcal_m^{\rm f-an}$ is replaced by the minimal formal analytic structure $\Xfrak_m'$ on $X'$ which refines $(\Xcal')^{\rm f-an}$ such that $f$ extends to a morphism $\phi_m:\Xfrak_m' \rightarrow \Acal_m^{\rm f-an}$. Note that we obtain $\Xfrak_m'$  and $\phi_m$ by applying \ref{formal refinement} to $\Ccal_m:=\frac{1}{m}\Ccal$ instead of $\Ccal$. In this sense, we may use in the following the description and the properties of $\Xfrak_m'$ from Section 7. Similarly as in step 1, we have
the following weak limit of regular Borel measures on $X'$:
\begin{equation} \label{can measure 1}
\mu = \lim_{ m \to \infty} m^{-2d} \sum_Z \deg_{\tilde{\phi}_m^*(\tilde{\Lcal})}(Z) \delta_{\xi_Z},
\end{equation}
where $Z$ ranges over all irreducible components of $\tilde{\Xfrak}_m'$. Note that $\Val(\xi_Z)$ is in the relative interior of a  simplex $\Delta:=(\Delta_S,S)$ of $S(\Xcal')$ for a unique $S \in {\rm str}(\tilde{\Xcal}')$. If $\Delta$ is degenerate with respect to $f$, then we claim
\begin{equation} \label{degenerate degree}
\deg_{\tilde{\phi}_m(\tilde{\Lcal})}(Z) = 0. 
\end{equation}
By definition, the simplex $\rho:=f_{\rm aff}(\Delta)$ is contained in $\Sigma$ and has dimension $<d$. Formula \eqref{Val and aff} yields $\valbar(f(\xi_Z)) \in \overline{\rho}$. 
By projection formula, it is enough to show that $Y:=\tilde{\phi}_m(Z)$ is no irreducible component of $\tilde{\Xcal}_m$. We argue by contradiction. We apply Theorem \ref{periodic transverse} to the irreducible component $Y$ using that $\overline{\Ccal_m}$ is transversal to $ \valbar(\Xan)$. We conclude that $Y$ corresponds to an equivalence class of transversal vertices, i.e. there is a unique $d$-codimensional $\Deltabar_m \in \overline{\Ccal_m}$ such that the torus orbit in $\tilde{\Acal}_m$ associated to $\relint(\Deltabar_m)$ is dense in $Y$.  Since $ \overline{\Ccal_m}$ is $\Sigmabar$-generic, we have $\overline{\rho} \cap \Deltabar_m = \emptyset$. But $(f(\xi_Z))\sptilde = \tilde{\phi}_m(\tilde{\xi}_Z)$ is the generic point of $Y$ and hence contained in the above open dense torus orbit. This means $\valbar(f(\xi_Z))  \in \relint(\Deltabar_m)$ (see Proposition \ref{torus orbit 2}) leading to a contradiction and proving \eqref{degenerate degree}.


Let $\Delta=(\Delta_S,S)$ be a canonical simplex of $S(\Xcal')$ which is non-degenerate with respect to $f$. We have seen in \ref{new setting} that $\Delta_S$ is $d$-dimensional and that $S$ is a $\tilde{\kdop}$-rational point contained in $d+1$ irreducible components of $\tilde{\Xcal}'$. Now we use Section 7, with $\Delta_S, \Sigma_S$ playing the role of $\Delta(d,\pi)$ and $\Sigma(d,\pi)$, to express properties of $\Xfrak_m'$ in terms of the polytopal decomposition $\Dcal_m:=(f_{\rm aff}^{(0)})^{-1}(\Ccal_m)$ of $\Sigma_S$ (see \ref{new setting}). By Remark \ref{additional remark}, the irreducible components $Z$  of $\tilde{\Xfrak}_m'$ with $\Val(\xi_Z) \in \relint(\Delta)$ correspond bijectively to the vertices $\ub'$ of $\Dcal_m$ contained in $\relint(\Sigma_S)$. For such  a $Z$, the corresponding   $\ub'$ is given by omitting the coordinate $u_0$ of $\ub:= \Val(\xi_Z) \in \relint(\Delta_S)$. Moreover, $Z$ is isomorphic to the toric variety $Y_{\ub'}$ associated to the vertex $\ub'$ and we will identify them later. This replaces step 2. 

Note that the point $\xi_{\ub'}$ of Corollary \ref{u points} is in the skeleton of the formal scheme $\Ucal_{\Sigma_S}$ associated to the standard simplex $\Sigma_S$  
(see \cite{Ber5}, 4.2 and Theorem 4.3.1). This skeleton consists of $\Sigma_S$ itself. We have seen in Remark \ref{additional remark} that there is a formal affine neighbourhood $\Ucal'$ of $S$ in $\Xcal'$ and an \'etale morphism $\phi_0:\Ucal' \rightarrow \Ucal_{\Sigma_S}$ with $\ub'=\val(\phi_0(\xi_Z))$. This proves $\phi_0(\xi_Z)=\xi_{\ub'}$ and hence  $\xi_Z \in S(\Xcal')$ (see \cite{Ber5}, Corollary 4.3.2). By the identification in \ref{skeleton}, $\Val$ is the identity on $S(\Xcal')$ and hence $\xi_Z \in \Delta$. By \eqref{can measure 1} and \eqref{degenerate degree}, we conclude  that the support of $\mu$ is contained in the union of the non-degenerate simplices with respect to $f$.

To prove the remaining formula, we may assume that $\Omega$ is a polytope contained in such a $\Delta=(\Delta_S,S)$. Note that $f_{\rm aff}^{(0)}$ from \ref{new setting} extends to an affine map $f_0:\rdop^d \rightarrow \rdop^n$ which is also one-to-one. The polytopal decomposition $\Dcal:=f_0^{-1}(\Ccal)$ of $\rdop^d$ is periodic with respect to the lattice $\Lambda_S$ from \ref{new setting}. Similarly as in step 3, we deduce from \eqref{can measure 1} the formula
\begin{equation} \label{can measure 2}
\mu(\Omega) = \frac{\vol(\Omega)}{\vol(\Lambda_S)} \cdot \sum_{\ub'} \deg_{\tilde{\phi}^*(\tilde{\Lcal})}(Y_{\ub'}),
\end{equation}
where $\ub'$ ranges over the vertices of $\Dcal$ modulo $\Lambda_S$. Since no multiplicities occur, the argument is easier here and will be omitted. We have seen in step 4 that $f_\Lcal$ is a strongly polyhedral convex function with respect to $\Ccal$ (Corollary \ref{ample reduction}) and hence $g:=f_\Lcal \circ f_0$ is a strongly polyhedral convex function with respect to $\Dcal$. As in \eqref{toric measure 7}, we conclude that
\begin{equation} \label{can measure 3}
\deg_{\tilde{\phi}^*(\tilde{\Lcal})}(Y_{\ub'}) = d! \cdot \vol(\{{\ub'}\}^g).
\end{equation}
If $\ub'$ ranges over the vertices of $\Dcal$, the rational polytopes $\{{\ub'}\}^g$ are the $d$-dimensional polytopes of the dual polytopal decomposition $\Dcal^g$ of $\rdop^d$. Since $\Dcal^g$ is $\Lambda_S^L$-periodic, the formula in the claim follows from \eqref{can measure 2} and \eqref{can measure 3}. \qed

\begin{cor} \label{non ample generalization}
If we do not require that $L_0, \dots, L_d$ are ample in Theorem \ref{explicit generalization}, then $\mu$ is still supported in the union of the non-degenerate simplices $(\Delta_S,S)$ of $S(\Xcal')$ and the restriction of $\mu$ to such a simplex is still a Haar measure.
\end{cor}

\proof This follows from Theorem \ref{explicit generalization} by multilinearity as in Remark \ref{arbitrary line bundles}. \qed

\begin{ex} \label{special case} \rm 
We consider the special case $X=A$ in Theorem \ref{explicit generalization}. Using $A_\kdop^{\rm an} = (\Tor)_\kdop^{\rm an}/M$, the points $\xi_\ub$ from Corollary \ref{u points} form a canonical subset $S(A)$ of $A_\kdop^{\rm an}$ which we call the {\it skeleton} of $A$. By \cite{Ber}, Example 5.2.12 and Theorem 6.5.1, this is a closed subset of $A_\kdop^{\rm an}$ and $\valbar$ restricts to a homeomorphism from $S(A)$ onto $\rdop^n/\Lambda$ which we use for identification.

By a combinatorial result of Knudson and Mumford (\cite{KKMS}, Chapter III), there is a rational triangulation $\Ccalbar$ of $\rtor$ (even refining any given rational polytopal decomposition) and $m \in \ndop \setminus \{0\}$ such that for every  $\Delta \in \Ccal$, the simplex $m \Delta$ is ${\rm GL}(n,\zdop)$-isomorphic to a $\zdop^n$-translate of the standard simplex $\{\ub \in \rdop_+^n \mid u_1 + \dots + u_n \leq 1 \}$. Then the formal $\kdop^\circ$-model $\Acal$ of $A$ associated to $\Ccalbar$ is strictly semistable. By the way, K\"unnemann generalized this construction to prove the existence of projective strictly semistable $\kdop^\circ$-models for arbitrary abelian varieties  (see \cite{Ku1} and also the erratum in \cite{Ku2}, 5.8). By the second step in the proof of Theorem \ref{explicit generalization}, the skeleton of $\Acal$ agrees with $S(A)$. We get a triangulation of $S(A)$ corresponding to $\Ccalbar$. 

We apply Theorem \ref{explicit generalization} with $X'=X=A$ and $\Xcal'=\Acal$. Note that the non-degenerate simplices $(\Delta_S,S)$ correspond to the $n$-dimensional simplices  of $\Ccalbar$ and hence the lattice $\Lambda_S^{L_j}$ does not depend on the choice of the stratum $S$. We conclude that the measure $\mu$ from Theorem \ref{explicit generalization} is supported in $S(A)$ and corresponds to a Haar measure on $\rtor$. By Proposition \ref{Chern continuity}, it has total measure $\deg_{L_1, \dots, L_n}(A)$. Using multilinearity for non-ample line bundles, this proves the following result:
\end{ex}

\begin{cor} \label{Haar measure on A}
Let $\overline{L}_1, \dots , \overline{L}_n$ be canonically metrized line bundles on the totally degenerate  abelian variety $A$ from above. Then $c_1(\overline{L}_1) \wedge \cdots \wedge c_1(\overline{L}_n)$ is supported in the skeleton $S(A)$ and corresponds to the Haar measure on $\rtor$ with total measure $\deg_{L_1, \dots, L_n}(A)$.
\end{cor}

\appendix
\section{Convex geometry}

In this appendix, we gather notions and results from convex geometry.

\begin{art} \rm \label{convexe geometry}
A {\it polyhedron} $\Delta$  in $\rdop^n$ is an intersection of finitely many closed half-spaces $\{\ub \in \rdop^n \mid \mb_i \cdot \ub \geq c_i\}$. We say that $\Delta$ is {\it $\Gamma$-rational} if we may choose all $\mb_i \in \zdop^n$ and all $c_i \in \Gamma$. A {\it closed face} of $\Delta$ is either $\Delta$ itself or has the form $H \cap \Delta$ where $H$ is the boundary of a closed half-space containing $\Delta$. An {\it open face} of $\Delta$ is a closed face  without all its properly contained closed faces. We denote by ${\Int(\Delta)}$ the topological interior of $\Delta$ in $\rdop^n$ and by $\relint(\Delta)$ the unique open face of $\Delta$ which is dense in $\Delta$. 

A bounded  polyhedron is called a {\it  polytope}. By linear algebra, a polytope is $\Gamma$-rational if and only if all vertices are in $\Gamma^n$ and the edges have rational slopes.
A ($\Gamma$-rational) {\it polytopal set} $S$ in $\rdop^n$ is a finite union of ($\Gamma$-rational) polytopes in $\rdop^n$. $S$ is said to have {\it pure dimension $d$} if all maximal polytopes of $S$ have dimension $d$. 

A {\it polytopal complex} $\cal C$ in $\rdop^n$ is a locally finite set of polytopes such that 
\begin{itemize}
\item[(a)] $\Delta \in {\cal C} \text{ $\Rightarrow$ all closed faces of $\Delta$ are in $\cal C$}$;
\item[(b)] $\Delta, \sigma \in {\cal C} \text{ $\Rightarrow$ $\Delta \cap \sigma$ is either empty or a closed face of $\Delta$ and $\sigma$}$.
\end{itemize}
The polytopal complex is called {\it $\Gamma$-rational} if every $\Delta \in \cal C$ is $\Gamma$-rational. A {\it polytopal decomposition} of $S \subset \rdop^n$ is a polytopal complex with $S= \cup_{\cal C} \Delta$. It is easy to see that every $\Gamma$-rational polytopal set has a finite $\Gamma$-rational polytopal decomposition. A {\it triangulation} of $S$ is a polytopal decomposition consisting only of simplices.  A polytopal complex $\Dcal$ {\it subdivides} $\Ccal$ if every polytope $\Delta$ in $\Ccal$ has a polytopal decomposition in $\Dcal$. 

A {\it cone} $\sigma$ in $\rdop^n$ is centered at $\mathbf 0$, i.e. it is characterized by $\rdop_+ \sigma = \sigma$. Its {\it dual} is defined by 
$$\check{\sigma}:= \{\ub' \in \rdop^n \mid \ub \cdot \ub' \geq 0 \, \; \forall \ub \in \sigma\}.$$
\end{art}

\begin{art} \rm \label{local cone and transversality}
Let $S$ be a locally finite union of polytopes in $\rdop^n$. The {\it local cone} ${\rm LC}_\ub(S)$ at $\ub$ is defined by 
$${\rm LC}_\ub(S):= \{ {\mathbf w} + \ub \mid 
\text{${\mathbf w} \in \rdop^n, \, [0,\ve){\mathbf w} +\ub \subset S$ for some $\ve >0$}\}.
$$
Then $S$  is said to be {\it concave in $\ub \in S$} if the convex hull of ${\rm LC}_\ub(S)$ is an affine subspace of $\rdop^n$. The set $S$ is called {\it totally concave} if it is concave in all $\ub \in S$.

Let $S$ be a locally finite union of $d$-dimensional polytopes in $\rdop^n$. A polytopal decomposition $\Ccal$ of $\rdop^n$ is said to be {\it transversal} to $S$ if the polytopal set $\Delta \cap S$ is either empty or of pure dimension $d-\codim(\Delta)$ for every $\Delta \in \Ccal$. 

If $\Ccal$ is transversal to $S$ and $\Delta \in \Ccal$ is of codimension $d$, then $\Delta \cap S$ consists of finitely many points. Such points are called {\it transversal vertices} of $\Ccal \cap S$. Two transversal vertices are called {\it equivalent} if they are contained in the same open face of $\Ccal$. 
\end{art}

\begin{art} \rm \label{strongly polyhedral function}
Let $\Ccal$ be a polytopal decomposition of $\rdop^n$. A {\it strongly polyhedral convex  function} $f$ with respect to $\Ccal$ is a  convex function $f: \rdop^n \rightarrow \rdop$ such that the $n$-dimensional $\Delta \in \Ccal$ are the maximal subsets of $\rdop^n$ where $f$ is affine, i.e. there are $\mb_\Delta \in \rdop^n$, $c_\Delta \in \rdop$ with
$$f(\ub) = \mb_\Delta \cdot \ub + c_\Delta$$
for every $\ub \in \Delta$. The vector $\mb_\Delta$ is called the {\it peg} of $\Delta$. 
Recall that $f$ is {\it convex} if 
\begin{equation} \label{convex}
f(rx+sy) \leq r f(x) +s f(y)
\end{equation}
for $x,y \in \rdop^n$ and $r,s \in [0,1]$ with $r+s=1$. {\bf Warning:} In the theory of toric varieties, convex functions are defined the opposite way! Here, we follow the terminology from analysis and we call a convex function $f$ {\it strictly convex} if we have $<$ in \eqref{convex} for $x \neq y$ and $0<r<1$. 
\end{art}

\begin{art} \label{dual complex} \rm 
There is a {\it dual complex} $\Ccal^f$ of $\Ccal$ realized in $\rdop^n$ with respect to $f$: For $\sigma \in \Ccal$, let 
$$\Star(\sigma):= \{\Delta \in \Ccal \mid \sigma \subset \Delta\}, \quad \Star_n(\sigma):= \{\Delta \in \Ccal \mid \sigma \subset \Delta, \, \dim(\sigma) = n\}.$$
The convex hull of $\{\mb_\Delta \mid \Delta \in \Star(\sigma)\}$ is a polytope denoted by $\sigma^f$. These polytopes form the dual complex $\Ccal^f$ (\cite{McM}, Theorem 3.1). It may happen that $\Ccal^f$ does not cover $\rdop^n$. For more details and biduality, we refer to \cite{McM}. 
\end{art}

\begin{art} \rm \label{alternative description of dual complex}
For a vertex $\ub_0$ of $\Ccal$, the polytope $\{\ub_0\}^f$ depends only on the local cones ${\rm LC}_{\ub_0}(\Delta)$, $\Delta \in \Ccal$, and the values of $f$ in a small neighbourhood of $\ub_0$. We have
$$\{\ub_0\}^f = \{\omega  \in \rdop^n \mid \ub \in \Delta \in \Star_n(\Delta) \Rightarrow \omega \cdot (\ub - \ub_0) \leq \mb_\Delta \cdot (\ub - \ub_0) \}.$$
This follows from \cite{Oda}, A.3 and Lemma 2.12. Indeed, for a $d$-dimensional polytope $\sigma$ with vertex $\ub_0$,  $\{\ub_0\}^f$  is the polytope associated to the fan in $\ub_0$ by the theory of toric varieties and 
$$\sigma^f = \{\ub_0\}^f \cap \left(\mb_\Delta + \sigma^\bot\right), \quad \dim(\sigma) + \dim(\sigma^f) = n,$$
where $\Delta$ is any element of $\Star_n(\sigma)$ (see \cite{Oda}, Corollary A.19).  \end{art}

\begin{art} \label{mixed volume} \rm
For compact convex subsets $P$ and $Q$ of $\rdop^n$, we have the {\it Minkowski sum} $$P+Q:=\{\ub + \ub' \mid \ub \in P, \ub' \in Q\}.$$ This is again a compact convex set. Similarly or by associativity, we define the Minkowski sum for more than two summands. For  a non-negative real number $\lambda$, we use $\lambda P:= \{\lambda \ub \mid \ub \in P\}$. There is a unique symmetric real function $V(P_1, \dots, P_n)$ on the set of compact convex subsets of $\rdop^n$ which is multilinear with respect to the above operations and which satisfies
$$V(P, \dots, P) = \vol(P).$$
The number $V(P_1, \dots, P_n)$ is called the {\it mixed volume} of $P_1, \dots, P_n$. The mixed volume is monotone increasing with respect to inclusion and hence it is non-negative. Moreover, it follows from translation invariance that $V(P_1, \dots, P_n)>0$ if all the $P_j$ have non-empty interiors. For proofs and more details, we refer to \cite{BZ}, Chapter 4.
\end{art}

\small

{\small Walter Gubler, Fachbereich Mathematik, Universit\"at Dortmund,
 D-44221 Dortmund, walter.gubler@mathematik.uni-dortmund.de}

\begin{thebibliography}{EGA III}      




\bibitem[Ber1]{Ber}{V.G. Berkovich: Spectral theory and analytic geometry over non-archimedean fields.  Mathematical Surveys and Monographs, 33. Providence, RI:  AMS. ix, 169 p. (1990).}

\bibitem[Ber2]{Ber2}{V.G. Berkovich:  \'Etale cohomology for non-archimedean analytic spaces. Publ. Math. IHES 78, 5--161 (1993).}

\bibitem[Ber3]{Ber3}{V.G. Berkovich: Vanishing cycles for formal schemes. Invent. Math.  115, No.3, 539--571 (1994).}

\bibitem[Ber4]{Ber4}{V.G. Berkovich:  Smooth $p$-adic analytic spaces are locally contractible. Invent. Math.  137, No.1, 1--84 (1999).}

\bibitem[Ber5]{Ber5}{V.G. Berkovich: Smooth $p$-adic analytic spaces are locally contractible. II. Adolphson, Alan (ed.) et al., Geometric aspects of Dwork theory. Vol. I. Berlin: de Gruyter. 293--370 (2004).}

\bibitem[BiGr]{BiGr}{R. Bieri, J.R.J. Groves: The geometry of the set of characters induced by valuations. J. Reine Angew. Math. 347, 168--195 (1984).}

\bibitem[BoGu]{BG}{E. Bombieri, W.Gubler: Heights in Diophantine geometry. Cambridge University Press. xvi, 652 p. (2006).}




\bibitem[Bo]{Bo}{S. Bosch: Zur Kohomologietheorie rigid analytischer R\"aume. Manuscr. Math. 20, 1--27 (1977).}

\bibitem[BGR]{BGR}{S. Bosch, U. G\"untzer, R. Remmert: Non-Archimedean analysis. A systematic approach to rigid analytic geometry. Grundl. Math. Wiss., 261. Berlin etc.: Springer Verlag. xii, 436 p. (1984).}


\bibitem[BL1]{BL1}{S. Bosch, W. L\"utkebohmert: N\'eron models from the rigid analytic viewpoint. J. Reine Angew. Math. 364, 69--84 (1986).}

\bibitem[BL2]{BL2}{S. Bosch, W. L\"utkebohmert: Degenerating abelian varieties. Topology 30, No.4, 653--698 (1991).}

\bibitem[BL3]{BL3}{S. Bosch, W. L\"utkebohmert: Formal and rigid geometry. I: Rigid spaces. Math. Ann. 295, No.2, 291--317 (1993).}

\bibitem[BL4]{BL4}{S. Bosch, W. L\"utkebohmert: Formal and rigid geometry. II: Flattening techniques. Math. Ann. 296, No.3, 403--429 (1993).}



\bibitem[Bou]{Bou}{N. Bourbaki: \'El\'ements de Math\'ematique. 
Fasc. XXX. Alg\`ebre commutative. Chap. 5: Entiers. Chap. 6: 
Valuations. Actualit\'es scientifiques et industrielles. 1308.  
Paris: Hermann, 207 (1964).}

\bibitem[BZ]{BZ}{Y. D. Burago, V. A. Zalgaller: Geometric inequalities. Translated from the Russian by A. B. Sossinsky. Grundlehren der Mathematischen Wissenschaften, 285. Springer Series in Soviet Mathematics. Berlin: Springer Verlag. xiv, 331 pp. (1988).}


\bibitem[Ch]{Ch}{A. Chambert-Loir: Mesure et \'equidistribution sur les espaces de Berkovich. J. Reine Angew. Math. 595, 215--235 (2006).} 


\bibitem[EKL]{EKL}{M. Einsiedler, M. Kapranov, D. Lind: Non-archimedean amoebas and tropical varieties. arXiv:math.AG/0408311, v2, 1--19.} 





\bibitem[FvdP]{FvdP}{J. Fresnel, M. van der Put: Rigid analytic geometry and its applications. Progress in Mathematics, 218. Boston, MA: Birkh\"auser. xii, 296 p. (2004).}

\bibitem[Fu1]{Fu1}{W. Fulton: Intersection theory. Ergebnisse der  Mathematik und ihrer Grenzgebiete. 3. Folge, Bd. 2. Berlin etc.: Springer-Verlag. xi, 470 p. (1984).}

\bibitem[Fu2]{Fu2}{W. Fulton: Introduction to toric varieties. The 1989 W. H. Roever lectures in geometry. Annals of Math. Studies, 131.  Princeton, NJ: Pinceton University Press. xi, 157 p. (1993).}











\bibitem[EGA IV]{EGA IV}{A. Grothendieck, J. Dieudonn\'e: \'El\'ements de g\'eometrie alg\'ebrique. IV: \'Etude locale des sch\'emas et des morphismes de sch\'emas (Quatrieme partie). Publ. Math. IHES, 1--361 (1967).}


\bibitem[Gu1]{Gu1}{W. Gubler: Heights of subvarieties over $M$-fields. F. Catanese (ed.), Arithmetic geometry. Proceedings of a symposium, Cortona, 1994. Cambridge: Cambridge University Press.  Symp. Math.  37, 190--227 (1997).}

\bibitem[Gu2]{Gu2}{W. Gubler: Local heights of subvarieties over non-archimedean fields. J. Reine Angew. Math. 498, 61--113 (1998).} 

\bibitem[Gu3]{Gu3}{W. Gubler: Local and canonical heights of subvarieties. 
Ann. Sc. Norm. Super. Pisa, Cl. Sci. (Ser. V), 2, No.4., 711--760 (2003).}

\bibitem[Gu4]{Gu4}{W. Gubler: The Bogomolov conjecture for totally degenerate abelian varieties. arXiv:math.NT/0609387v2, 1--21.}

\bibitem[Ha]{Ha}{U. Hartl: Semi-stable models for rigid-analytic spaces. Manuscr. Math. 110, No.3, 365--380 (2003).}




\bibitem[dJ]{dJ}{A. J. de Jong: Smoothness, semi-stability and alterations. Publ. Math. IHES 83, 51--93 (1996).}

\bibitem[KKMS]{KKMS}{G. Kempf, F. Knudsen, D. Mumford, B. Saint-Donat: Toroidal embeddings. I. LNM 339. Berlin etc.: Springer-Verlag. viii, 209 p. (1973).}

\bibitem[Ku1]{Ku1}{K. K\"unnemann: Projective regular models for abelian varieties, semistable reduction, and the height pairing. Duke Math. J. 95, No. 1, 161--212 (1998).}

\bibitem[Ku2]{Ku2}{K. K\"unnemann: Height pairings for algebraic cycles on abelian varieties. Ann. Sci. \'Ec. Norm. Sup\'er. (4) 34, No. 4, 503--523 (2001).}




\bibitem[Kl]{Kl}{S. L. Kleiman: Toward a numerical theory of ampleness. Ann. Math. (2)  84, 293--344 (1966).}




\bibitem[McM]{McM}{P. McMullen: Duality, sections and projections of certain Euclidean tilings. Geom. Dedicata 49, No.2, 183--202 (1994).}

\bibitem[Mi]{Mi}{G. Mikhalkin: Amoebas of algebraic varieties and tropical geometry. S. Donaldson (ed.) et al., Different faces of geometry. New York, NY: Kluwer Academic/Plenum Publishers. Int. Math. Ser. 3,  257--300 (2004).}






\bibitem[Mu]{Mu}{D. Mumford: An analytic construction of degenerating abelian varieties over complete rings. Compositio Math. 24, 239--272 (1972).}






\bibitem[Oda]{Oda}{T. Oda: Convex bodies and algebraic geometry. An introduction to the theory of toric varieties. Ergebnisse der Mathematik und ihrer Grenzgebiete. 3. Folge, Bd. 15. Berlin etc.: Springer-Verlag. viii, 212 p. (1988).}

\bibitem[Ru1]{Ru1}{W. Rudin: Functional Analysis. McGraw-Hill 
Series in Higher Mathematics. New York etc.: McGraw-Hill. xiii, 397 p. (1973).}

\bibitem[Ru2]{Ru2}{W. Rudin: Real and complex analysis.
3rd. edition. New York, NY: McGraw-Hill. xiv, 416 p. (1987).}




\bibitem[Ul]{Ul}{P. Ullrich:  The direct image theorem in formal and rigid geometry. Math. Ann. 301, No.1, 69--104 (1995).}








\end{thebibliography}
\end{document}